\documentclass{amsart}
\usepackage[margin=2cm]{geometry}
\usepackage[svgnames]{xcolor}
\usepackage{enumitem}

\usepackage{mathtools}
\usepackage{tensor}
\usepackage{amsmath,amssymb}
\usepackage{tikz}
\usetikzlibrary{decorations.pathmorphing, calc,shapes.geometric,patterns} 
\pgfdeclarelayer{background layer}
\pgfdeclarelayer{foreground layer}
\pgfsetlayers{background layer,main,foreground layer}
\usetikzlibrary{automata}
\usepackage{multimedia}
\usepackage{multicol}
\usepackage[T1]{fontenc} 
\usepackage{times}
\usepackage{ytableau}
\usepackage{graphicx}
\usepackage[normalem]{ulem}
\usepackage{chngcntr}
\counterwithin{figure}{section} 

\newcommand{\omitt}[1]{}
\newcommand\mc[1]{\mathcal{#1}}

\newcommand{\Z}{\mathbb Z}

\newcommand{\R}{\mathbb R}

\newcommand{\len}[1]{\ell(#1)}
\newcommand{\n}{n}

\newcommand{\brac}[2]{{\langle #1\mid#2\rangle}}

\DeclareMathOperator{\shi}{shi} 
\DeclareMathOperator{\ish}{ish} 
\DeclareMathOperator{\rank}{rank} 
\DeclareMathOperator{\area}{area} 
\DeclareMathOperator{\bounce}{bounce}
\DeclareMathOperator{\ad}{ad}
\DeclareMathOperator{\im}{Im}

\def\highestroot{highest}
\def\hecke{\mathfrak H}
\def\KLRing{\mc A}
\def\Lang{\mc L}
\newcommand{\Region}{{\mathfrak A}_\m}
\newcommand{\ldes}[1]{{\mathcal L}(#1)} 
\newcommand{\rdes}[1]{{\mathcal R}(#1)} 
\newcommand{\m}{{m}}  

\newcommand{\A}{A}
\newcommand{\Afund}{\A_0}
\newcommand{\domCham}{D}
\newcommand{\Hy}{{H}}
\newcommand{\Hak}[2]{\Hy_{#1, #2}} 
\newcommand{\Hp}[2]{{{\Hak{#1}{#2}}^+}}

\newcommand{\Hn}[2]{{{\Hak{#1}{#2}}^-}}
\newcommand{\Shi}[1]{{\mathcal Shi}_{#1}}
\newcommand{\Ish}[1]{{\mathcal Ish}_{#1}}
\newcommand{\Cox}[1]{{\mathcal Cox}_{#1}}
\newcommand{\Cat}[1]{{\mathcal Cat}_{#1}}
\newcommand{\mCat}[2]{{\Cat{#1}^{#2}}}
\newcommand{\mShi}[2]{{\Shi{#1}^{#2}}}
\newcommand{\rec}[1]{{\mc{R}(#1)}}

\newcommand{\roots}{{\Delta}}
\newcommand{\simples}{{\Pi}}
\newcommand{\pos}{{\roots^+}}

\newcommand{\e}{{\varepsilon}}  

\newcommand\acoord[2]{\tensor[^{}_{#1}]a{^{}_{#2}}}
\newcommand\aprimecoord[2]{\tensor[^{}_{#1}]{\tilde{a}}{^{}_{#2}}}
\def\regNum{{\mathfrak r}}
\def\bddRegNum{{\mathfrak b}}
\def\S{{{\mathcal S}}}
\def\F{{\mathbb F}}
\newcommand{\poin}[1]{P_{#1}}
\newcommand\Sym[1]{\mathfrak{S}_{#1}}
\def\Sn{\Sym{n}}
\newcommand\affS[1]{{\widehat{\mathfrak{S}}}_{#1}}
\def\affSn{\affS{n}}
\def\affA[1]{{\widehat{\mathcal{A}}}_n}
\def\affA{{\widehat{\mathcal{A}}}_n}

\newcommand\Triangle[3]{\node at (1,1) {#1};\node at (0,0) {#2};\node at (2,0) {#3};}
\newcommand\Lusztig[1]{\textcolor{black}{#1}}
\newcommand\Levear[1]{\textcolor{black}{#1}}
\newcommand\Christos[1]{\textcolor{black}{#1}}
\newcommand\DGdO[1]{\textcolor{black}{#1}}
\newcommand\Zaslavsky[1]{\textcolor{black}{#1}}

\numberwithin{equation}{section}
\theoremstyle{definition}
  \newtheorem{theorem}{Theorem}[section]
  
  \newtheorem{proposition}[theorem]{Proposition}
  \newtheorem{lemma}[theorem]{Lemma}

  \newtheorem{definition}[theorem]{Definition}
  \newtheorem{example}[theorem]{Example}
\theoremstyle{remark}

\usepackage[colorinlistoftodos]{todonotes}
\newcommand{\Susanna}[1]{\todo[size=\tiny,inline,color=red!5]{#1
      \\ \hfill --- Susanna}}

\newcommand{\sftodo}[1]{\Susanna{#1}}
\newsavebox{\smlmat}
\savebox{\smlmat}{$\left(\begin{smallmatrix}1&2&3&4&5&6\\2&4&5&3&6&1\end{smallmatrix}\right)$}

\newsavebox{\smlmatInv}
\savebox{\smlmatInv}{$\left(\begin{smallmatrix}1&2&3&4&5&6\\6&1&4&2&3&5\end{smallmatrix}\right)$}

\title{A survey of the Shi arrangement}
\author{Susanna Fishel}
\address{School of Mathematical and Statistical Sciences, Arizona State University, P.O. Box 871804, Tempe, AZ 85287-1804, USA} \email{sfishel1@asu.edu}

\begin{document}
\maketitle

In \cite{lusztig1983}, Lusztig defined a
map $\sigma$ from affine permutations of $n$ to partitions of $n$. He conjectured that
for any partition $\lambda$ of $n$, $\sigma^{-1}(\lambda)$ is a two-sided cell.
Shi \cite{shiBook} proved \Lusztig{part of }this conjecture. As a byproduct, Shi introduced the {\em
  Shi arrangement} of hyperplanes and found a few of its remarkable
properties. The Shi arrangement has since become a central object in
algebraic combinatorics. This article is intended to be a fairly
gentle introduction to the Shi arrangement, intended for readers with
a modest background in combinatorics, algebra, and Euclidean
geometry. After background material in Section~\ref{sec:background},
this introduction to the arrangement will be by way of a discussion in
Section~\ref{sec:origin} of how it arose, some of its marvelous
enumerative properties in Section~\ref{sec:enumeration}, and some of
its surprising connections to algebra in
Section~\ref{sec:connections}. For some brief comments on recent
extensions, see Section~\ref{sec:further} and for an incomplete list
of topics we left out, see Section~\ref{sec:themes}.


\section{Background}\label{sec:background} In this section, we will give very brief introductions to some of the ingredients needed to define the Shi arrangement.

\subsection{Root systems and Coxeter group notation}
\label{sec:cox} Let $V$ be a finite dimensional real vector space with fixed inner product $\brac{}{}$.
We'll use $\roots$ to denote a {\em root system}: a finite set of
vectors in $V$ which satisfies

\begin{enumerate}
\item $\roots\cap \R\alpha=\{\alpha,-\alpha\}$ for all $\alpha\in\roots$ and
\item $s_{\alpha}\roots=\roots$ for all $\alpha\in\roots,$
\end{enumerate}
where $s_{\alpha}$ is the reflection about the hyperplane with normal
$\alpha$.  We use $\pos$ to denote a choice of positive roots of
$\roots$, so that $\roots=\pos\cup-\pos$, and $\simples$ to denote the simple roots, which are a basis
for the $\R$-span of $\roots.$ The reflections
$S=\{s_{\alpha}\}_{\alpha\in\simples}$ generate a finite reflection
group $W$. The {\em rank} of the system and of $W$ is the dimension of the space spanned by $\roots$. 

\begin{figure}


\begin{tikzpicture}[scale=.33]
\tikzstyle{every node}=[font=\footnotesize]

\path [<->,draw = gray,  draw opacity = 1](150:6) -- (330:6) node[black,above] {$\alpha_2$};

\path [<->,draw = gray, draw opacity = 1](270:6) -- (90:6) node[black,above] {$\alpha_1$};

\path [<->,draw = gray, draw opacity = 1](210:6) -- (30:6) node[black,above] {$\theta$};

\path [draw = red, very thick, draw opacity = 1] (240:8) -- (60:8) node[black,above] {$H_{\alpha_2,0}$};

\path [draw = red,very thick, draw opacity = 1] (-8,0) -- (8,0) node[black,right] {$H_{\alpha_1,0}$};

\path [draw = red,very thick,draw opacity = 1] (0,0) -- (0,0);

\path [draw = black,draw opacity = 1] (120:8) -- (300:8) node[black,below] {$H_{\theta,0}$};
\path [draw = red ,very thick,draw opacity = 1,shift={(2.5,0)}] (120:8) -- (300:8) node[black,below] {$H_{\theta,1}$};

\end{tikzpicture}

\caption{The roots and reflecting hyperplanes of affine type $A_2$. The
  reflection $s_1$ (resp. $s_2$) flips the plane over the hyperplane
  $\Hak{\alpha_1}{0}$ ($\Hak{\alpha_2}{0}$). The reflection $s_0$ reflects over $\Hak{\theta}{1}$.}
\label{fig:A2}
\end{figure}
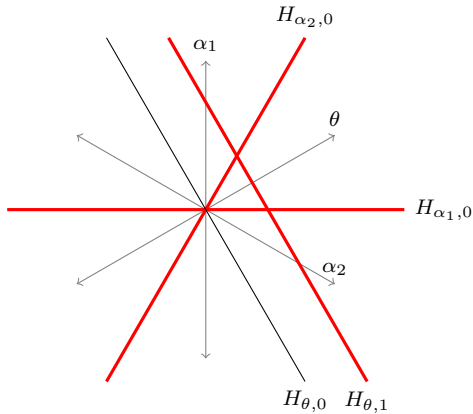

Coxeter groups generalize finite reflection groups.  Let $W$ be a
group with a set of generators $S\subset W$. Let $m_{st}$ be the order
of the element $st$, with $s,t\in S$. If there is no relation between
$s$ and $t$, we set $m_{st}=\infty$. \omitt{If $(W,S)$ satisfies} If $W$ has a presentation such that 

\omitt{A group $W$ is a {\em Coxeter group} if it has a set of
  generators $S\subset W$ and a set $\{m_{st}\}_{s,t\in S}$ which
  satisfy}

\begin{enumerate}
\item $m_{ss}=1$
  \item for $s,t\in S,s\neq t$,
    $1<m_{st}\leq\infty$,
\end{enumerate}
then $W$ is a {\em Coxeter group}. We refer to $(W,S)$ as a {\em
  Coxeter system}.
If $m_{st}\in\{2,3,4,6\}$ when $s\neq t$, then the Coxeter group is
called {\em crystallographic} and, if finite, is a {\em Weyl} group. It
is also a reflection group.  The product in any order of all the elements in $S$ is called a {\em Coxeter element}; all Coxeter elements for a given $W$ are conjugate and their order is the {\em Coxeter number} of $W$.

The expression for the reflection $s_{\alpha}$, $\alpha\in\roots$,
is $$s_{\alpha}(v)=v-2\frac{\brac{v}{\alpha}}{\brac{\alpha}{\alpha}}\alpha$$
for $v\in V$. For any $k\in\R$, we can define an {\em affine
  reflection} $s_{\alpha,k}$ by
$$s_{\alpha,k}(v)=v-2\frac{\brac{v}{\alpha}-k}{\brac{\alpha}{\alpha}}\alpha.$$
We define the {\em affine Weyl group} to be the group generated by all
affine reflections $s_{\alpha,k}$ for $\alpha\in\roots$ and
$k\in\Z$. It is also a Coxeter group. Its simple reflections are the
simple reflections $S$ of the finite Weyl group, plus an extra
reflection, $s_0$, about a translate of a certain other hyperplane in
the arrangement. See Figure~\ref{fig:A2}. Its root system, for us, is
the root system for the corresponding finite Weyl group. Our proofs
will not be detailed enough to need the full set of affine roots, and
we will not define them. Given a root system $\roots$, we write
$W_{\roots}$ for the corresponding finite group. Please see
\cite{humphreys1990} or \cite{kac1990} for more information.

We put a partial order on any root system. The {\em root poset} of
$\roots$ is the set of positive roots $\pos$, partially ordered by
setting $\alpha\leq\beta$ if $\beta-\alpha$ is a nonnegative linear
combination of simple roots.  If the root system is irreducible, then
there is a unique highest root relative to this ordering. We will
denote this root by $\theta$. See Figure~\ref{fig:catObjs} for a
picture of the root poset for type $A_4$.

Let $W$ be a Coxeter group.  Every $w\in W$ has an expression as a
product of elements of $S$: $w=s_{i_1}\cdots s_{i_k}$. If $k$ is minimal
among all expressions for $w$, then $k$ is the {\em length}
$\len{w}$. Any expression for $w$ of length $\len{w}$ is a {\em
  reduced expression}.

We will often refer to the {\em type} of a group or root system,
particularly ``type $A$'' and ``affine type $A$,'' which are the symmetric group or affine symmetric group if we are referring to groups. Please see
\cite{humphreys1990} for more information on the classification of
finite reflection groups and Coxeter groups. Humphreys and
\cite{bjorner-brenti2005} are good sources for the definitions of
irreducible, Bruhat order, and other material omitted
here.

\omitt{Now for a warning. Although one of the charms of the Shi arrangement is
how well it generalizes, in our examples and discussions we are mainly
sticking with type $A$. This is where it began and there is plenty to
discuss here. That is, here usually $W=\affSn$ and
$S=\{s_0,\ldots,s_{n-1}\}$, defined in Section~\ref{subsec:affSn}.|}

\subsubsection{Type $A$}
\label{subsec:affSn}
We will be seeing type $A$ often, so we'll be a little more concrete.
For $A_{n-1}$, the vector space $V$ is $\{ (a_1, \ldots, a_\n ) \in
\R^n \mid a_1+\cdots+a_n =0 \}$. Let $\{ \e_1, \ldots, \e_\n \}$ be
the standard basis of $\R^\n$ and $\brac{\,}{\,}$ be the bilinear form
for which this is an orthonormal basis. The set of roots is $\roots =
\{\e_i - \e_j \mid i \neq j\}$ and a root $\alpha\in\roots$ is
positive, written $\alpha > 0$, if $\alpha \in \pos = \{\e_i - \e_j
\mid i < j\}$. Set $\alpha_i = \e_i - \e_{i+1}$; the simple roots are
$\{\alpha_1, \ldots, \alpha_{\n-1} \}$.  The set $\simples$ is a basis of $V$.

The Coxeter group that Shi studied was the affine symmetric group
$\affSn$, and we review that here.  There are several possible
descriptions, here we give one due to Lusztig \cite{lusztig1983}. It
is the set of permutations $w$ of $\Z$ such that \begin{enumerate}
\item $w(i+n)=w(i)+n$ for all $i\in\Z$
\item $\sum_{i=1}^nw(i)=\binom{n}{2}$
\end{enumerate}


It's a Coxeter group: for any $i$, $0\leq i<n$, $s_i$ corresponds to
the permutation
$$t\mapsto\begin{cases}t&\text{if $t\bmod n\neq i$ and $t\bmod n\neq
  i+1$}\\t-1 &\text{if $t\bmod n=i+1$}\\t+1
&\text{if $t\bmod n=i $}\end{cases}$$

The set of reflections $S$ is $\{s_1, \ldots, s_{\n-1}, s_0\}$ and

$$\affSn = \langle s_1, \ldots, s_{\n-1}, s_0 \rangle.$$ 
\omitt{For $n=2$, ${\affS{2}} = \langle s_1, s_0 \mid s_i^2 = 1 \rangle$.}


The affine symmetric group contains the symmetric group $\Sn$ as a
subgroup.  $\Sn$ is the subgroup generated by the $s_i$, $0<i<n$.  We
identify $\Sn$ as permutations of $\{1, \ldots, \n\}$ by identifying
$s_i$ with the simple transposition $(i,i+1)$. We act on the right, as did Shi.

\subsection{A taste of Coxeter combinatorics, type $A$}
The number of parking functions and the Catalan numbers appear in
every discussion of the Shi arrangement. We'll define the parking
functions when we first see them, in Section~\ref{sec:stanEnum}, but
we'll collect some facts on the Catalan objects here, mostly type
$A$. There are an awful lot of Catalan objects, but only a few of them
will appear in this survey.

A {\em partition} is a finite sequence $\lambda = (\lambda_1, \ldots,
\lambda_r)$ of positive integers in decreasing order:
$\lambda_1\geq\lambda_2\geq\ldots\geq\lambda_r >0$. We identify a
partition with its {\em Young diagram}, the left-justified array of
boxes where the $i^{\text{th}}$ row from the top has $\lambda_i$
boxes. A box is called {\em removable} (respectively {\em addable}) if
we can remove (respectively add) it and still have a diagram of a
partition. We use $| \lambda | = \lambda_1+\cdots+\lambda_r$ and
$\len{\lambda}=r$.

\subsubsection{Set partitions} 
\begin{figure}[ht]
\begin{tikzpicture}

\begin{scope}[thick,font=\small]
\def\a{3}
\def\b{-2}
\def\c{.4}
\begin{scope}[shift={(0*\a,0*\b)},scale=\c]
  \node[below,blue] (1) at (0,0) {};
  \node[below,blue] at (0,-\c) {1};
\fill[red] ($(1)$) circle (2pt);
\node[below,blue] (2) at (1,0) {};
\node[below,blue]  at (1,-\c) {2};
\fill[red] ($(2)$) circle (2pt);
\node[below,blue] (3) at (2,0) {};
\node[below,blue]  at (2,-\c) {3};
\fill[red] ($(3)$) circle (2pt);
\node[below,blue] (4) at (3,0) {};
\node[below,blue]  at (3,-\c) {4};
\fill[red] ($(4)$) circle (2pt);
\node[below,blue] (5) at (4,0) {};
\node[below,blue]  at (4,-\c) {5};
\fill[red] ($(5)$) circle (2pt);

\draw (1) .. controls (1,.5) .. (3);
\draw (4) .. controls (3.5,.5) .. (5);
\draw (2) .. controls (2,.5) .. (4);

\end{scope}
\end{scope}

  \begin{scope}[shift={(2.4,-2)}]
\begin{scope}[font=\small, scale=.9]
\def\a{1}
\def\b{1}
\node (1) at (0*\a,0*\b) {$\e_1-\e_2$};
\node (2) at (2*\a,0*\b) {$\e_2-\e_3$};
\node (3) at (4*\a,0*\b) {$\e_3-\e_4$};
\node [draw=green,shape=ellipse](4) at (6*\a,0*\b) {$\e_4-\e_5$};

\node[draw=green,shape=ellipse] (12) at (1*\a,1*\b) {$\e_1-\e_3$};
\node[draw=green,shape=ellipse] (23) at (3*\a,1*\b) {$\e_2-\e_4$};
\node[draw=green,shape=ellipse] (34) at (5*\a,1*\b) {$\e_3-\e_5$};

\node[draw=green,shape=ellipse] (13) at (2*\a,2*\b) {$\e_1-\e_4$};
\node [draw=green,shape=ellipse](24) at (4*\a,2*\b) {$\e_2-\e_5$};

\node [draw=green,shape=ellipse](14) at (3*\a,3*\b) {$\e_1-\e_5$};

\draw[thick,red] (1)--(12);
\draw[thick,red] (2)--(12);
\draw[thick,red] (2)--(23);
\draw[thick,red] (3)--(23);
\draw[thick,red] (3)--(34);
\draw[thick,red] (4)--(34);

\draw[thick,red] (12)--(13);
\draw[thick,red] (23)--(13);
\draw[thick,red] (23)--(24);
\draw[thick,red] (34)--(24);

\draw[thick,red] (13)--(14);
\draw[thick,red] (24)--(14);
\end{scope}

  \end{scope}
\begin{scope}[shift={(10.5,0)}]

\node[font=\scriptsize] at (0,0) {\begin{ytableau}
*(yellow)\scriptstyle \e_1-\e_5 &*(yellow)\scriptstyle \e_1-\e_4&*(yellow)\scriptstyle \e_1-\e_3&*(white)\scriptstyle \e_1-\e_2\\
*(yellow)\scriptstyle \e_2-\e_5&*(yellow)\scriptstyle \e_2-\e_4&*(white)\scriptstyle \e_2-\e_3\\
*(yellow)\scriptstyle \e_3-\e_5&*(white)\scriptstyle \e_3-\e_4\\
*(yellow)\scriptstyle \e_4-\e_5
\end{ytableau}};

\end{scope}
\begin{scope}[shift={(12.3,-1.7)},scale=.65]

\draw [help lines] (0,0) grid (5,5);
\draw [very thick, red] (0,0)--(0,1);
\draw [very thick, red] (0,1)--(1,1);
\draw [very thick, red] (1,1)--(1,3);
\draw [very thick, red] (1,3)--(2,3);
\draw [very thick, red] (2,3)--(2,4);
\draw [very thick, red] (2,4)--(3,4);
\draw [very thick, red] (3,4)--(3,5);
\draw [very thick, red] (3,5)--(5,5);

\end{scope}

\end{tikzpicture}

\caption{On the left is the nonnesting set partition
  $\pi=\{\{1,3\},\{2,4,5\}\}$. Next to it is the root poset of type
  $A_4$ with the filter corresponding to $\pi$ circled. The third
  figure is the partition inside the staircase which corresponds to
  $\pi$. On the right is the corresponding Dyck path.}
\label{fig:catObjs}
\end{figure}
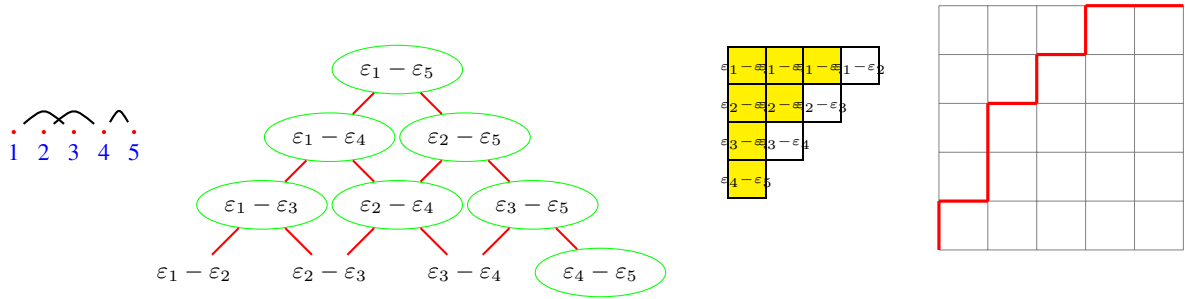

We denote the set $\{1,2,\ldots,n\}$ by $[n]$. The nonempty subsets $B_1,\ldots,B_k$ of $[n]$ are a {\em set partition} of $[n]$ if they are pairwise disjoint and their union is $[n]$. We denote the set partition $\{B_1,\ldots,B_k\}$ by $B_1|\cdots |B_k$. For example, $13|256|4$ is a partition of $[6]$. The {\em arc diagram} of a set partition $\pi$ is defined as follows: place the numbers $1,2,\ldots,n$ in order on a line and draw an arc between each pair $i<j$ such that
\begin{itemize}
\item $i$ and $j$ are in the same block of $\pi$, and
\item there is no $k$ such that $i<k<j$ and $i$, $k$, and $j$ are in the same block. 
\end{itemize}
See Figure~\ref{fig:catObjs}.

The partition $\pi$ has $k$ blocks if and only if it has $n-k$ arcs. This is easy to see if the partition has no arcs. Consider a partition with $k$ blocks and $n-k$ arcs, where $i$ and $j$ are in different blocks. Suppose we add an arc from $i$ to $j$. We have joined $i$'s and $j$'s blocks, and we now have $k-1$ blocks and $n-k+1$ arcs.

A set partition is called {\em noncrossing} if there does not exist
$i<j<k<l$ such that there is an arc from $i$ to $k$ and an arc from
$j$ to $l$. There are $C_n$ noncrossing set partitions,
where $$C_n=\frac{1}{n+1}\binom{2n}{n}$$ is the Catalan number (type
$A$). It is called {\em nonnesting} if there does not exist $i<j<k<l$
such that there is an arc from $i$ to $l$ and an arc from $j$ to $k$.
\ytableausetup{mathmode,boxsize=2em}

\subsubsection{Dyck paths}
A {\em Dyck path} of length $n$ is a lattice path which starts at
(0,0), takes only north or east steps of length 1, never goes below
the line $y=x$, and ends at $(n,n)$. A north step followed by an east
step is called a {\em valley} of the path.

\subsubsection{Root poset}
\label{subsec:rootposet}
An {\em ideal} of a poset $P$ is a subset $I$ of the elements of $P$
such that if $x\in I$ and $y\leq x$ then $y\in I$. A {\em filter} is
like an ideal, except that the condition becomes that $x\in I$ and
$x\leq y$ implies $y\in I$. A subset $X$ of the elements of $P$ is an
{\em antichain} if no two elements in $X$ are comparable. An ideal is
determined by its maximal elements, which form an antichain, just as a
filter is by its minimal elements. These antichains are also called
nonnesting partitions; in type $A$ the bijection to the nonnesting
partition described above is simply sending the root $\e_i-\e_j$ in an
antichain to the arc from $i$ to $j$. There are Catalan number of
ideals, filters, and antichains in the root poset defined in
Section~\ref{sec:cox}. We can map a filter $F$ in the root poset for
type $A_{n-1}$ to a partition $\lambda$ whose diagram fits inside the
staircase partition $(n-1,n-2,\ldots,1)$ by the rule that $\e_i-\e_j\in
F$ if and only if the box in row $i$ and column $n+1-j$ is in the
diagram of $\lambda$. The minimal elements of $F$ correspond to the
removable boxes of $\lambda$.

There are well-known bijections among all these objects; please see Stanley's book on the subject \cite{stanley2015}.

\subsection{Deformation of Coxeter arrangements}
\label{sec:hyperplanes}
A {\em (real) hyperplane arrangement} $\mathcal{H}$ is a set of
hyperplanes, possibly affine hyperplanes, in a real vector space.  For
us, the vector space will be $V$, the span of some root system
$\roots$, with a fixed inner product $\brac{}{}$ which is $W_{\roots}$
invariant.  We'll be looking at connected components of a hyperplane
arrangement's complement $V \setminus \bigcup_{H \in \mathcal{H}} H$.
We will refer to these as the {\em regions} of the arrangement. The
closure $\bar{R}$ of the region $R$ is a convex polyhedron. A {\em
  face} of $\mc{H}$ is a nonempty set of the form $\bar{R}\cap x$,
where $x$ is an intersection of hyperplanes in $\mc{H}$. The dimension
of a face is the dimension of its affine span. See Stanley
\cite{stanley1996} for more details.  A {\em wall} $H$ of $R$ is a
hyperplane $H\in\mc{H}$ such that $\dim(H\cap R)=\dim(H)$. The word
``bounded'' applied to a region has its usual meaning: a region is {\em
  bounded} if there is a real number $M$ such that all points in the
region are within distance $M$ of the origin. Let $\regNum(\mc{A})$
and $\bddRegNum(\mc{A})$ be the number of regions and number of
bounded regions, respectively, of the arrangement $\mc{A}$.

Let $\roots$ be a root system. 
The roots (plus the integers) define a system of affine hyperplanes
$$\Hak{\alpha}{k} = \{ v \in V \mid \brac{v }{\alpha} = k \}.$$  Note
$\Hak{-\alpha}{-k} =\Hak{\alpha}{k}$. In type $A$, we will sometimes
write $x_i-x_j=k$ instead of $\Hak{\alpha_i+\cdots+\alpha_{j-1}}{k}$.

The Coxeter arrangement, also called the braid arrangement, is defined 
$$\Cox{\roots}=\{\Hak{\alpha}{0}:\alpha\in\pos\}.$$ We give the regions of this arrangement the special name {\em chambers}. Each chamber
corresponds to an element of $W=W(\roots)$.   The {\em
  dominant} chamber of $V$ is $\bigcap_{i=1}^{n-1} \Hp{\alpha_i}{0}$, where $\Hp{\alpha_i}{k}$ is the half-space $\{ v \in V \mid \brac{v }{\alpha} \geq k \}$. the dominant chamber corresponds to the identity of $W$. It is also referred to as
the fundamental chamber in the literature.

The affine Coxeter arrangement is all integer translates of the
hyperplanes in $\Cox{\roots}$; that is, it is the whole system of hyperplanes
$\{\Hak{\alpha}{k}\}_{\alpha\in\pos,k\in\Z}$. In this arrangement,
each region is called an {\em alcove} and the {\em fundamental alcove}
is $\Afund$, the interior of $ \Hn{\theta}{1} \cap \bigcap_{\alpha\in\simples}\Hp{\alpha}{0}$, where $\theta$ is the \highestroot root.  A {\em dominant alcove} is one contained in the
dominant chamber.

We also have the $m$-Catalan arrangement:
$$\mCat{\roots}{m}=\{\Hak{\alpha}{r}:\alpha\in\roots,0\leq r\leq m\}.$$



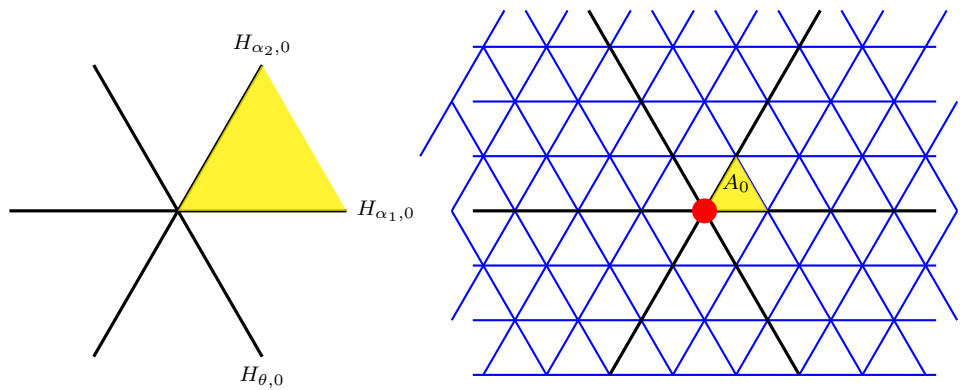
\begin{figure}[ht]
\begin{tikzpicture}[scale=.28]
\begin{scope}
\tikzstyle{every node}=[font=\footnotesize]

\path [draw = black, very thick, draw opacity = 1] (240:8) -- (60:8) node[black,above] {$H_{\alpha_2,0}$};

\path [draw = black,very thick, draw opacity = 1] (-8,0) -- (8,0) node[black,right] {$H_{\alpha_1,0}$};

\path [draw = red,very thick,draw opacity = 1] (0,0) -- (0,0);

\path [draw = black,very thick,draw opacity = 1] (120:8) -- (300:8) node[black,below] {$H_{\theta,0}$};
\filldraw[fill=yellow,draw opacity = 0, fill opacity = .8] (60:8)--(0,0)--(8,0)--cycle;

\end{scope}

  \begin{scope}[shift={(25,0)}]

\tikzstyle{every node}=[font=\footnotesize]

\def\a{2.59807}

\path [draw = blue,  thick, draw opacity = 1,shift = {(3,0)}] (240:9) -- (60:11);

\path [draw = blue,  thick, draw opacity = 1,shift = {(6,0)}] (240:9) -- (60:11);

\path [draw = blue,  thick, draw opacity = 1,shift = {(-3,0)}] (240:9) -- (60:11);

\path [draw = blue,  thick, draw opacity = 1,shift = {(9,0)}]
(240:9) -- (60:6);

\path [draw = blue,  thick, draw opacity = 1,shift = {(12,0)}] (240:9) -- (60:0);



\path [draw = blue,  thick, draw opacity = 1,shift = {(-6,0)}] (240:9) -- (60:11);

\path [draw = blue,  thick, draw opacity = 1,shift = {(-9,0)}]
(240:6) -- (60:11);

\path [draw = blue,  thick, draw opacity = 1,shift = {(-12,0)}] (240:0) -- (60:11);

\path [draw = blue,  thick, draw opacity = 1,shift = {(-13.5,\a)}]
(240:0) -- (60:8);


\path [draw = blue, thick, draw opacity = 1, shift = {(0,3*\a)}] (-11,0) -- (11,0);

\path [draw = blue, thick, draw opacity = 1, shift = {(0,\a)}] (-11,0) -- (11,0);

\path [draw = blue, thick, draw opacity = 1,shift = {(0,2*\a)}] (-11,0) -- (11,0);

\path [draw = blue, thick, draw opacity = 1,shift = {(0,-\a)}]
(-11,0) -- (11,0);

\path [draw = blue, thick, draw opacity = 1,shift = {(0,-2*\a)}] (-11,0) -- (11,0);

\path [draw = blue, thick, draw opacity = 1,shift = {(0,-3*\a)}] (-11,0) -- (11,0);


\path [draw = blue, thick,draw opacity = 1,shift={(3,0)}] (120:11) -- (300:9);

\path [draw = blue, thick,draw opacity = 1,shift={(6,0)}] (120:11)
-- (300:9);
\path [draw = blue, thick,draw opacity = 1,shift={(9,0)}] (120:11)
-- (300:6);

\path [draw = blue, thick,draw opacity = 1,shift={(12,0)}] (120:11) -- (300:0);


\path [draw = blue, thick,draw opacity = 1,shift={(-3,0)}]
(120:11) -- (300:9);

\path [draw = blue, thick,draw opacity = 1,shift={(-6,0)}]
(120:9) -- (300:9);

\path [draw = blue, thick,draw opacity = 1,shift={(-9,0)}] (120:6) -- (300:9);

\path [draw = blue, thick,draw opacity = 1,shift={(-12,0)}] (120:0) -- (300:9);

\path [draw = black, very thick, draw opacity = 1] (240:9) -- (60:11);

\path [draw = black,very thick,draw opacity = 1] (120:11) -- (300:9);

\path [draw = black,very thick, draw opacity = 1] (-11,0) -- (11,0);

\filldraw[fill=yellow,draw opacity = 0, fill opacity = .8] (60:3)--(0,0)--(3,0)--cycle;
\node at (0,0) [circle, fill = red]{};
\node at (1.5,.5*\a) {$\Afund$};


  \end{scope}
\end{tikzpicture}
\caption{The Coxeter arrangement for $A_2$ is on the left, with shaded fundamental chamber. On the right is all its translates, with fundamental alcove shaded..} 
\end{figure}

Let $W$ be an affine Weyl group and $V$ the vector space spanned by
its roots.  $W$ acts on $V$ via affine linear transformations, and
acts freely and transitively on the set of alcoves.  In affine type
$A_{n-1}$, $s_i$ reflects over $\Hak{\alpha_i}{0}$ for $1\leq i\leq
n-1$ and $s_0$ reflects over $\Hak{\theta}{1}$, where $\theta$ is the
highest root. We identify each alcove $\A$ with the unique $w \in W$
such that $\A = \Afund w$.  \omitt{Each simple generator $s\in S$,
  acts by reflection with respect to the simple root $\alpha_i$. In
  other words, it acts by reflection over the hyperplane
  $\Hak{\alpha_i}{0}$.  The element $s_0$ acts as reflection with
  respect to the affine hyperplane $\Hak{\theta}{1}$, where $\theta$
  is the highest root.} For example, if $w$ is the element of affine
$A_2$ whose reduce decomposition is $s_0s_1$, then $\Afund w$ is the
image of $\Afund$ after reflecting first across $\Hak{\theta}{1}$ and
then across $\Hak{\alpha_1}{0}$. See Figure~\ref{fig:inverses}.

We can be even more specific for type $A_{n-1}$.
The action on $V$ is given by
\begin{gather*}
s_i (a_1, \ldots, a_i, a_{i+1}, \ldots, a_\n) = (a_1, \ldots, a_{i+1},
a_i,\ldots, a_\n) \quad \text{ for $i \neq 0$, and} \\ s_0 (a_1,
\ldots, a_\n) = (a_\n +1, a_2 , \ldots, a_{\n-1}, a_1 - 1).
\end{gather*}

Note $\Sn$ preserves $\brac{\;}{\;}$, but $\affSn$ does not.

An alcove $\A$ can be described by the hyperplanes it is between. For
example, in Figure~\ref{fig:inverses}, the alcove labeled by $20$ is
between $\Hak{\alpha_1}{1}$ and $\Hak{\alpha_1}{2}$, between
$\Hak{\alpha_2}{0}$ and $\Hak{\alpha_2}{1}$ and between
$\Hak{\theta}{1}$ and $\Hak{\theta}{2}$. Given a positive root
$\alpha$, there is a unique integer $k=k_{\alpha}(\A)$ such that
$k<\brac{\alpha}{x}<k+1$ for all $x\in\A$. Let
$K(\A)=\{k_{\alpha}(\A)\}_{\alpha\in\pos}$ denote the set of
coordinates for $\A$, indexed by the positive roots. Returning to Figure~\ref{fig:inverses}, the alcove labeled by $20$ has coordinates $k_{\alpha_1}=k_{\theta}=1$ and $k_{\alpha_2}=0$.

Shi characterized the sets of integers which can arise as $K(\A)$
for some alcove $\A$; for type $A$ in \cite[Chapter 6]{shiBook}
and for general affine Weyl groups in \cite{shi1987b}. The situation
in general is rather messy, but if we assume our root system is an
irreducible crystallographic one, then Shi found (see also
\cite{athanasiadis2004}) that a collection of integers indexed by the
positive roots $\pos$ corresponds to an alcove if and only if

\begin{equation}
\label{eq:alcoveCoords}
k_{\alpha}+k_{\beta}\leq k_{\alpha+\beta}\leq k_{\alpha}+k_{\beta}+1
\end{equation}
 for all $\alpha,\beta\in\pos$ such that $\alpha+\beta\in\pos$. We call the set $K(\A)$ the coordinates of $\A$.

\section{Origin}
\label{sec:origin}
\subsection{Kazhdan-Lusztig Cells}\label{subsec:KL}
We provide a bare-bones introduction to Kazhdan-Lusztig theory. For
more information, see \cite{shiBook}, \cite{kazhdan-lusztig1979},
\cite{bjorner-brenti2005}. If you are willing to believe that Kazhdan
and Lusztig defined an equivalence relation on the elements of a
Coxeter group, then skip this section. We include this collection of
definitions for completeness.  We'll need the Hecke algebra $\hecke$
and the Kazhdan-Lusztig polynomials in order to define the $W$-graph
and then the cells. We will prove none of our claims.

Let $W$ be a Coxeter group and let $S$ be the corresponding set of
simple reflections. We first follow \cite{kazhdan-lusztig1979}, who follow
\cite{bourbaki1968}, for the definition of the Hecke algebra. Let
$\KLRing$ be the ring of Laurent polynomials in the indeterminate $q^{1/2}$ with integral coefficients. The Hecke algebra
$\hecke=\hecke(W,S)$ is a free module over $\KLRing$ with basis
$T_w$, one for each $w\in W$. The multiplication is defined by the
rules

\begin{enumerate}
\item $T_wT_{w'}=T_{ww'}$ if $\len{ww'}=\len{w}+\len{w'}$
\item $(T_s+1)(T_s-q)=0$ if $s\in S$;
\end{enumerate}

here $\len{w}$ is the length of $w$. 

Now for the polynomials. The involution on $\KLRing$ $a\mapsto \bar{a}$ defined by
$\overline{q^{1/2}}=q^{-1/2}$ extends to an involution of the ring
$\hecke$: $$\overline{\sum a_wT_w}=\sum \overline{a_w}T^{-1}_{w^{-1}}.$$
Kazhdan's and Lusztig's theorem, Theorem~\ref{thm:KL} in this survey, asserts the existence of elements
$C_w\in\hecke$, one for each $w\in W$, and simultaneously defines the
Kazhdan-Lusztig polynomials $P_{y,w}$, where $y,x\in W$. The order is
the Bruhat order on $W$.

\begin{theorem}[\cite{kazhdan-lusztig1979}] \label{thm:KL}For any $w\in W$, there is a unique element $C_w\in\hecke$ such that
\begin{itemize}
\item $\overline{C}_w=C_w$
\item $C_w=\sum\limits_{y\leq w}(-1)^{\len{y}+\len{w}}q^{\len{w}/2-\len{y}}\overline{P_{y,w}}T_y$
\end{itemize}
where $P_{y,x}\in\KLRing$ is a polynomial in $q$ of degree at most
$\frac{1}{2}(\len{w}-\len{y}-1)$ for $y<w$, and $P_{w,w}=1$.
\end{theorem}

 Kazhdan and Lusztig used the polynomials to define a graph and from
 there the cells.  Now we follow the exposition given in
 \cite{bjorner-brenti2005}, simplified just a bit because we will not
 prove anything. For $u,w\in W$, define $\mu(u,w)$ to be the
 coefficient of $q^{\frac{1}{2}(\len{w}-\len{u}-1)}$ in $P_{u,w}(q)$ if
 $u<w$ and $\frac{1}{2}(\len{w}-\len{y}-1)$ is an integer; otherwise,
 $\mu(u,w)$ is 0.  Let $\ldes{w}$ be the set of left descents of $w$:
 $\ldes{w}=\{s\in S|sw<w\}$. \omitt{Then Bj\"orner and Brenti define} The
 directed, labeled graph $\tilde{\Gamma}_{(W,S)}^L$ is the graph
 with vertices $x\in W$ and edges $x\xrightarrow{(s,\mu)} y\in E$
 . There are two types of edges in $E$:
\begin{enumerate}
\item \label{nonloop}$x,y\in W$, $x\neq y$, either $\mu(x,y)\neq 0$ or
  $\mu(y,x)\neq 0$, and $s\in\ldes{x}\setminus\ldes{y}$. Let $\mu$ be either
  $\mu(x,y)$ or $\mu(y,x)$, whichever is not 0.
\item Loops at $x$: labeled by $s\in S$ and $$\mu=\begin{cases}1&\text{if $s\not\in\ldes{x}$},\\-1&\text{if $s\in\ldes{x}$}.\end{cases}$$
\end{enumerate}
The graph $\tilde{\Gamma}_{(W,S)}^R$ has an analogous definition,
using right descents of $w$: $\rdes{w}=\{s\in S|ws<w\}$.  The graph
$\tilde{\Gamma}_{(W,S)}^{LR}$ is the superposition of the
$\tilde{\Gamma}_{(W,S)}^L$ and $\tilde{\Gamma}_{(W,S)}^R$. We describe
the cells in graph theoretic terms. A directed graph is {\em strongly
  connected} if there is a directed path between all pairs of
vertices. A {\em strongly connected component} of a directed graph is
a maximal strongly connected subgraph. Finally, the {\em left cells}
are the strongly connected components of $\tilde{\Gamma}_{(W,S)}^L$,
the {\em right cells} the strongly connected components of
$\tilde{\Gamma}_{(W,S)}^R$, and the {\em two-sided cells} the strongly
connected components of $\tilde{\Gamma}_{(W,S)}^{LR}$.


The list of areas in math where cells appear is mind-boggling. Please
see the short survey by Gunnells \cite{gunnells2006}, for example, for
references, as well as for insight into the geometry of the cells. The book by Bj\"orner and Brenti \cite{bjorner-brenti2005} explains much of the combinatorial connection.

\subsection{Shi regions and Kahzdan-Lusztig cells}
\label{subsec:cells}
Shi was studying cells in \cite{shiBook}. He concentrated on the
affine Weyl groups of type $A$, because of the following conjectures
of Lusztig.  In \cite{lusztig1983}, Lusztig defined a map $\sigma$
from $\affSn$ to partitions of $n$. He conjectured that for any
partition $\lambda$ of $n$, $\sigma^{-1}(\lambda)$, a set of affine
permutations, is in fact a two-sided cell. \omitt{He showed in
\cite{lusztig1985} that the number of two-sided cells is the number of
partitions.} What's more, Lusztig also conjectured a formula for the
number of left (or right) cells which make up the two-sided cell
$\sigma^{-1}(\lambda)$.

\Lusztig{We'd like to say a few words on the history of the proof of
  these conjecture. In 1979, Vogan \cite{vogan} defined related cells,
  which we call V-cells following \cite{lusztig1985}. Shi called the
  V-cells RL-cells in \cite{shiBook}. In his 1983 Ph.D. thesis
  \cite{shiThesis}, later included in \cite{shiBook}, Shi determined
  and enumerated the left and right cells, both V-cells and cells as defined in \cite{kazhdan-lusztig1979} and described above. He confirmed
  Lusztig's conjectures on the image of $\sigma^{-1}$ and on the number
  of 2-sided V-cells, but not for 2-sided cells in the sense of
  \cite{kazhdan-lusztig1979}. In \cite{lusztig1983b}, Lusztig defines
  a function $a$, which he used in \cite{lusztig1985} to show that
  V-cells and cells coincide. This step completed the proof of the
  conjectures.}

The description of the map $\sigma$ is simple enough
and we define it here. Let $w\in\affSn$ and define $d_k=d_k(w)$ to be the
maximum size of a subset of $\Z$ whose elements are noncongruent to
each other modulo $n$ and which is a disjoint union of $k$ subsets
each of which has its natural order reversed by $w$. The partition
$\lambda$ is given by $(d_1,d_2-d_1,\ldots,d_n-d_{n-1})$

\setcounter{MaxMatrixCols}{20}
\begin{example}
Let $n=3$. The permutation $s_0$
is $$\begin{pmatrix}\cdots&-3&-2&-1&0&1&2&3&4&5&6&\cdots\\ \cdots&-2&-3&-1&1&0&2&4&3&5&3&\cdots\end{pmatrix}.$$
    The set $\{3,4\}$ has its order reversed by $s_0$ and there is no
    larger set, so $d_1=2$. The sets $\{3,4\}$ and $\{2\}$ show that
    $d_2=3$. Therefore $\sigma(s_0)=(2,1)$. \omitt{Shi writes his permutations acting on on the right, but we let them act on the left.} In our notation, $s_0s_2$ is

    $$\begin{pmatrix}\cdots&-3&-2&-1&0&1&2&3&4&5&6&\cdots\\\cdots&-2&-4&0&1&-1&3&4&2&6&7&\cdots\end{pmatrix}$$ and $\sigma(s_0s_2)=(2,1)$ also. The identity maps to $(1,1,1)$ under $\sigma$ and $s_1s_2s_1$, for example, maps to $(3)$.
  \end{example}

\begin{figure}[h]
\begin{tikzpicture}[scale=.28]

\tikzstyle{every node}=[font=\footnotesize]

\def\a{2.59807}

\filldraw[fill=gray,draw opacity=0,fill opacity=.8](6,-4*\a)--(3+4*1.5,-4*\a)--(3,0)--(0,0)--(-6,-4*\a)--(3-4*1.5,-4*\a)--(1.5,-\a)--cycle;

\filldraw[fill=red!10,draw opacity=0,fill opacity=.8](3+4*1.5,-4*\a)--(3,0)--(15,0)--cycle;

\filldraw[fill=red!10,draw opacity=0,fill opacity=.8](6,-4*\a)--(3-4*1.5,-4*\a)--(1.5,-\a)--cycle;

\filldraw[fill=gray, draw opacity=0,fill opacity=.8](3,0)--(1.5,\a)--(7.5,5*\a)--(9,4*\a)--(4.5,\a)--(13.5,\a)--(15,0)--cycle;

\filldraw[fill=red!10, draw opacity=0,fill opacity=.8](1.5,\a)--(7.5,5*\a)--(-4.5,5*\a)--cycle;

\filldraw[fill=red!10, draw opacity=0,fill opacity=.8](9,4*\a)--(4.5,\a)--(13.5,\a)--cycle;

\filldraw[fill=gray, draw opacity=0,fill opacity=.8](0,0)--(1.5,\a)--(-4.5,5*\a)--(-6,4*\a)--(-1.5,\a)--(-12+1.5,\a)--(-12,0)--cycle;

\filldraw[fill=red!10, draw opacity=0,fill opacity=.8](0,0)--(-12,0)--(-6,-4*\a)--cycle;

\filldraw[fill=red!10, draw opacity=0,fill opacity=.8](-6,4*\a)--(-1.5,\a)--(-12+1.5,\a)--cycle;


\draw [draw = blue,  thick, draw opacity = 1,shift = {(3,0)}] (240:12) -- (60:15-3);

\draw [draw = blue,  thick, draw opacity = 1,shift = {(6,0)}] (240:12) -- (60:15-6);

\draw [draw = blue,  thick, draw opacity = 1,shift = {(9,0)}](240:12) -- (60:15-9);

\draw [draw = blue,  thick, draw opacity = 1,shift = {(12,0)}] (240:12) -- (60:15-12);

\draw [draw = blue,  thick, draw opacity = 1,shift = {(15,0)}] (240:12) -- (60:0);


\draw [draw = blue,  thick, draw opacity = 1,shift = {(-3,0)}] (240:12-3) -- (60:15);

\draw [draw = blue,  thick, draw opacity = 1,shift = {(-6,0)}] (240:12-6) -- (60:15);

\draw [draw = blue,  thick, draw opacity = 1,shift = {(-9,0)}](240:12-9) -- (60:15);

\draw [draw = blue,  thick, draw opacity = 1,shift = {(-12,0)}] (240:0) -- (60:15);




\draw [draw = blue, thick, draw opacity = 1, shift = {(0,5*\a)}] (-12+7.5,0) -- (15-7.5,0);

\draw [draw = blue, thick, draw opacity = 1, shift = {(0,4*\a)}] (-12+6,0) -- (15-6,0);

\draw [draw = blue, thick, draw opacity = 1, shift = {(0,3*\a)}] (-12+4.5,0) -- (15-4.5,0);

\draw [draw = blue, thick, draw opacity = 1, shift = {(0,\a)}] (-12+1.5,0) -- (15-1.5,0);

\draw [draw = blue, thick, draw opacity = 1,shift = {(0,2*\a)}] (-12+3,0) -- (15-3,0);

\draw [draw = blue, thick, draw opacity = 1,shift = {(0,-\a)}](-12+1.5,0) -- (15-1.5,0);

\draw [draw = blue, thick, draw opacity = 1,shift = {(0,-2*\a)}] (-12+3,0) -- (15-3,0);

\draw [draw = blue, thick, draw opacity = 1,shift = {(0,-3*\a)}] (-12+4.5,0) -- (15-4.5,0);

\draw [draw = blue, thick, draw opacity = 1, shift = {(0,-4*\a)}] (-12+6,0) -- (15-6,0);

\draw [draw = blue, thick,draw opacity = 1,shift={(3,0)}] (120:15) -- (300:12);

\draw [draw = blue, thick,draw opacity = 1,shift={(6,0)}] (120:15)-- (300:12-3);

\draw [draw = blue, thick,draw opacity = 1,shift={(9,0)}] (120:15)-- (300:12-6);

\draw [draw = blue, thick,draw opacity = 1,shift={(12,0)}] (120:15) -- (300:12-9);

\draw [draw = blue, thick,draw opacity = 1,shift={(15,0)}] (120:15) -- (300:0);

\draw [draw = blue, thick,draw opacity = 1,shift={(-3,0)}](120:9) -- (300:12);

\draw [draw = blue, thick,draw opacity = 1,shift={(-6,0)}](120:6) -- (300:12);

\draw [draw = blue, thick,draw opacity = 1,shift={(-9,0)}] (120:3) -- (300:12);

\draw [draw = blue, thick,draw opacity = 1,shift={(-12,0)}] (120:0) -- (300:12);

\draw [draw = black, very thick, draw opacity = 1] (240:12) -- (60:15);

\draw [draw = black,very thick,draw opacity = 1] (120:12) -- (300:12);

\draw [draw = black,very thick, draw opacity = 1] (-12,0) -- (15,0);

\filldraw[fill=yellow,draw opacity = 0, fill opacity = .8] (60:3)--(0,0)--(3,0)--cycle;
\node at (0,0) [circle, fill = red]{};
\node at (1.5,.5*\a) {$\Afund$};

\omitt{
\node (x) at (240:12) {$x$};
\node (y) at (3-4*1.5,-4*\a) {$y$};
\node (z) at (300:3) {$z$};
\node  (t) at (300:12){$T$};
\node (u) at  (3+4*1.5,-4*\a){$u$};
\node (v) at (3,0) {$v$};
\node (w) at (0,0) {$w$};
 \filldraw[fill=red](t)--(u)--(v)--(w)--(x)--(y)--(z)--cycle;




}

\end{tikzpicture}
\caption{The cells for affine $A_2$. The yellow region is the two-sided cell $\sigma^{-1}((1,1,1))$, which is also a single left-cell.
  The six pink regions are left-cells, whose union is the two-sided cell $\sigma^{-1}((3))$. The three gray regions are also left cells, and their union is the two-sided cell $\sigma^{-1}((2,1))$. See \cite[Page 98]{shiBook}.}
\label{fig:shiP98}
\end{figure}
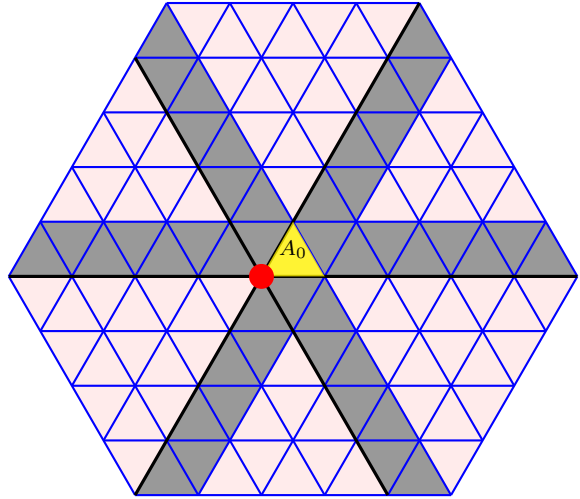

Shi proved \Lusztig{parts of }both of Lusztig's conjectures, and more. 
Shi used the identification of $\affSn$ with
alcoves to describe the cells of affine type $A$.
He showed that the two-sided cells correspond to connected sets of
alcoves, one set of alcoves for each partition $\lambda$ of $n$. A
two-sided cell is a disjoint union of left-cells. Inside the two-sided
cell corresponding to the partition $\lambda$, there is one left-cell
for each tabloid of shape $\lambda$. See Figure~\ref{fig:shiP98}.

What came to be known as the Shi arrangement was not initially defined
in terms of hyperplanes. Shi began by defining {\em rank $n$ sign types} as
triangular arrays $X=(x_{ij})_{1\leq i<j\leq n}$ with entries from
$\{+,-,\bigcirc\}$. The {\em admissible sign types} correspond to the regions of his arrangement. He defined them as the sign types which satisfy the following condition: for all $1\leq i<t<j\leq n$, the triple

\begin{center}
\begin{tikzpicture}[scale=.3, shift={(0,-1)}]\Triangle{$x_{ij}$}{$x_{it}$}{$x_{tj}$}\end{tikzpicture}
\end{center}

is a member of the set $G_A$ of admissible sign types of rank $3$ ($A_2$). $G_A$ is the set
\begin{equation}
\label{GA}
\begin{tikzpicture}[scale=.2]

 \def\a{6.5}
  \def\b{4.3}
\def\co{3.2}
  \def\c{1.4}
  \node at (-1.5,.5){${\LARGE \{}$};

  \begin{scope}[shift={(0,0)},scale=\c]
    \Triangle{$+$}{$+$}{$+$}
   \node at (\co,0){,};   
  \end{scope}

  \begin{scope}[shift={(1*\a,0)},scale=\c]
    \Triangle{$+$}{$+$}{$\bigcirc$}
      \node at (\co,0){,};
  \end{scope}
  
    \begin{scope}[shift={(2*\a,0)},scale=\c]
      \Triangle{$+$}{$\bigcirc$}{$+$}
        \node at (\co,0){,};
    \end{scope}

  \begin{scope}[shift={(3*\a,0)},scale=\c]
    \Triangle{$+$}{$+$}{$-$}
      \node at (\co,0){,};
  \end{scope}
  
    \begin{scope}[shift={(4*\a,0)},scale=\c]
      \Triangle{$+$}{$-$}{$+$}
        \node at (\co,0){,};
  \end{scope}

  \begin{scope}[shift={(5*\a,0)},scale=\c]
    \Triangle{$+$}{$\bigcirc$}{$\bigcirc$}
   \node at (\co,0){,};   
  \end{scope}

  \begin{scope}[shift={(6*\a,0)},scale=\c]
    \Triangle{$\bigcirc$}{$\bigcirc$}{$\bigcirc$}
   \node at (\co,0){,};   
  \end{scope}

    \begin{scope}[shift={(7*\a,0)},scale=\c]
    \Triangle{$\bigcirc$}{$+$}{$-$}
    \node at (\co,0){,};
  \end{scope}
    
      \begin{scope}[shift={(0*\a,-\b)},scale=\c]
    \Triangle{$\bigcirc$}{$-$}{$+$}
   \node at (\co,0){,};   
  \end{scope}

        \begin{scope}[shift={(1*\a,-\b)},scale=\c]
    \Triangle{$\bigcirc$}{$\bigcirc$}{$-$}
   \node at (\co,0){,};   
  \end{scope}

        \begin{scope}[shift={(2*\a,-\b)},scale=\c]
    \Triangle{$\bigcirc$}{$-$}{$\bigcirc$}
   \node at (\co,0){,};   
  \end{scope}

  \begin{scope}[shift={(3*\a,-\b)},scale=\c]
    \Triangle{$-$}{$-$}{$-$}
   \node at (\co,0){,};   
  \end{scope}

    \begin{scope}[shift={(4*\a,-\b)},scale=\c]
    \Triangle{$-$}{$+$}{$-$}
   \node at (\co,0){,};   
  \end{scope}

  \begin{scope}[shift={(5*\a,-\b)},scale=\c]
    \Triangle{$-$}{$-$}{$+$}
   \node at (\co,0){,};   
  \end{scope}

  \begin{scope}[shift={(6*\a,-\b)},scale=\c]
    \Triangle{$-$}{$-$}{$\bigcirc$}
   \node at (\co,0){,};   
  \end{scope}
  
  \begin{scope}[shift={(7*\a,-\b)},scale=\c]
    \Triangle{$-$}{$\bigcirc$}{$-$}
  \end{scope}

  \node at (7.3*\a+3.4,-\b+.5){${\LARGE \}}$};

\end{tikzpicture}
\end{equation}

Two comments on $G_A$. If we order the symbols $\{\bigcirc,+,-\}$ as $-<0<+$, then $G_A$ can be seen as the rank 3  sign types where either $x_{12}\leq x_{13}\leq x_{23}$ or $x_{23}\leq x_{13}\leq x_{12}$, together with $x_{13}=+,x_{12}=x_{23}=0$. The set $G_A$ has cardinality 16, which is $(n+1)^{n-1}$ for $n=3$.


Shi connected the admissible sign types to geometry using
\eqref{eq:alcoveCoords} and the map $\zeta$. If $K$ is the set of
coordinates for an alcove $\A$, then define the sign type
$X=\zeta(\A)$ by
$$x_{ij}=\begin{cases}+&\text{if $k_{ij}>0$}\\\bigcirc&\text{if
  $k_{ij}=0$}\\-&\text{if $k_{ij}<0$.}\end{cases}$$ He then calculated
the hyperplanes so that the regions defined by them were made up of
alcoves with the same image under the map $\zeta$. We use admissible
sign type, region in the Shi arrangement, and Shi region interchangeably.

Shi showed in \cite{shiBook} that the left-cells for affine type $A$
are themselves disjoint unions of admissible sign
types. \omitt{Admissible sign types are the regions of the hyperplane
  arrangement that came to be known as the Shi arrangement.}
Admissible sign types were not used directly in the proofs of the
Lusztig conjectures in Shi's monograph, but describe the structure of the cells.
They have taken on a life of
their own.

Later, in \cite{shi1987}, Shi extended the definition of admissible sign types, thereby generalizing the Shi arrangement. This is the definition we give below. 

We give the definition for any irreducible, crystallographic root system $\roots$. When the root system is type $A_{n-1}$, we will sometimes write $\Shi{n}$ instead of $\Shi{\roots}$.
\begin{definition}
\label{def:shi}
The Shi arrangement $\Shi{\roots}$ is the collection of
hyperplanes $$\{ \Hak{\alpha}{k} \mid \alpha \in \pos, 0 \leq
k \leq 1 \}.$$
\end{definition}

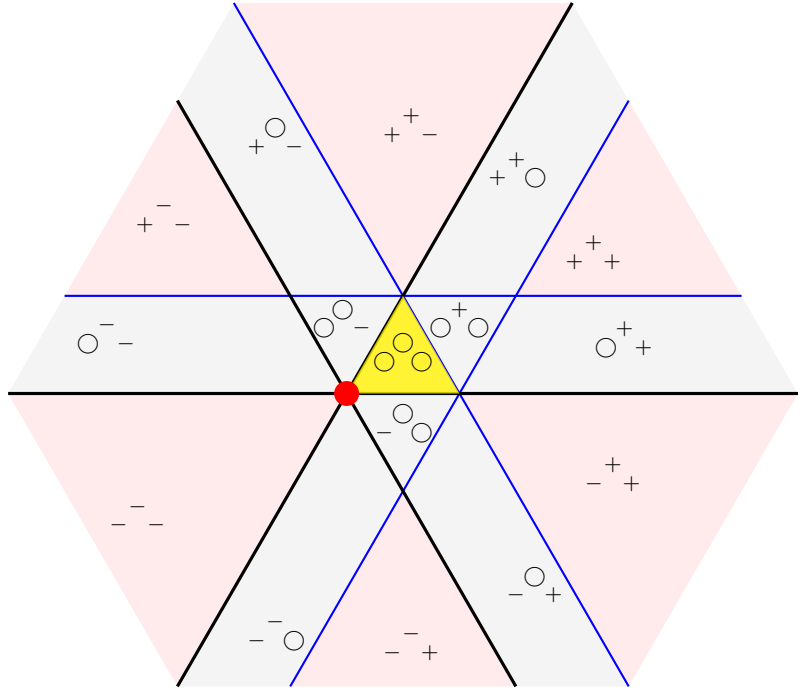
\begin{figure}[h]
\begin{tikzpicture}[scale=.5]

\tikzstyle{every node}=[font=\footnotesize]

\def\a{2.59807}

\filldraw[fill=gray!10,draw opacity=0,fill opacity=.8](6-1.5,-3*\a)--(3+3*1.5,-3*\a)--(3,0)--(0,0)--(-6+1.5,-3*\a)--(3-3*1.5,-3*\a)--(1.5,-\a)--cycle;

\filldraw[fill=red!10,draw opacity=0,fill opacity=.8](3+3*1.5,-3*\a)--(3,0)--(15-3,0)--cycle;

\filldraw[fill=red!10,draw opacity=0,fill opacity=.8](6-1.5,-3*\a)--(3-3*1.5,-3*\a)--(1.5,-\a)--cycle;

\filldraw[fill=gray!10, draw opacity=0,fill opacity=.8](3,0)--(1.5,\a)--(6,4*\a)--(7.5,3*\a)--(4.5,\a)--(12-1.5,\a)--(12,0)--cycle;

\filldraw[fill=red!10, draw opacity=0,fill opacity=.8](1.5,\a)--(6,4*\a)--(-3,4*\a)--cycle;

\filldraw[fill=red!10, draw opacity=0,fill opacity=.8](7.5,3*\a)--(4.5,\a)--(12-1.5,\a)--cycle;

\filldraw[fill=gray!10, draw opacity=0,fill opacity=.8](0,0)--(1.5,\a)--(-3,4*\a)--(-4.5,3*\a)--(-1.5,\a)--(-12+4.5,\a)--(-12+3,0)--cycle;

\filldraw[fill=red!10, draw opacity=0,fill opacity=.8](0,0)--(-12+3,0)--(-6+1.5,-3*\a)--cycle;

\filldraw[fill=red!10, draw opacity=0,fill opacity=.8](-4.5,3*\a)--(-1.5,\a)--(-12+4.5,\a)--cycle;


\draw [draw = blue,  thick, draw opacity = 1,shift = {(3,0)}] (240:12-3) -- (60:15-6);


\draw [draw = blue, thick, draw opacity = 1, shift = {(0,\a)}] (-12+4.5,0) -- (15-4.5,0);


\draw [draw = blue, thick,draw opacity = 1,shift={(3,0)}] (120:15-3) -- (300:12-3);

\draw [draw = black, very thick, draw opacity = 1] (240:12-3) -- (60:15-3);

\draw [draw = black,very thick,draw opacity = 1] (120:12-3) -- (300:12-3);

\draw [draw = black,very thick, draw opacity = 1] (-12+3,0) -- (15-3,0);

\filldraw[fill=yellow,draw opacity = 0, fill opacity = .8] (60:3)--(0,0)--(3,0)--cycle;

\def\r{7}
\begin{scope}[shift={(1,.33*\a)},scale=.5]
\Triangle{$\bigcirc$}{$\bigcirc$}{$\bigcirc$};
\end{scope}

\begin{scope}[shift={(2.5,.67*\a)},scale=.5]
\Triangle{$+$}{$\bigcirc$}{$\bigcirc$};
\end{scope}

\begin{scope}[shift={(-.6,.67*\a)},scale=.5]
\Triangle{$\bigcirc$}{$\bigcirc$}{$-$};
\end{scope}

\begin{scope}[shift={(1,-.4*\a)},scale=.5]
\Triangle{$\bigcirc$}{$-$}{$\bigcirc$};
\end{scope}

\begin{scope}[shift={(169:\r)},scale=.5]
\Triangle{$-$}{$\bigcirc$}{$-$};
\end{scope}

\begin{scope}[shift={(140:\r)},scale=.5]
\Triangle{$-$}{$+$}{$-$};
\end{scope}

\begin{scope}[shift={(110:\r)},scale=.5]
\Triangle{$\bigcirc$}{$+$}{$-$};
\end{scope}

\begin{scope}[shift={(80:\r)},scale=.5]
\Triangle{$+$}{$+$}{$-$};
\end{scope}

\begin{scope}[shift={(55:\r)},scale=.5]
\Triangle{$+$}{$+$}{$\bigcirc$};
\end{scope}

\begin{scope}[shift={(30:\r)},scale=.5]
\Triangle{$+$}{$+$}{$+$};
\end{scope}

\begin{scope}[shift={(10:\r)},scale=.5]
\Triangle{$+$}{$\bigcirc$}{$+$};
\end{scope}

\begin{scope}[shift={(340:\r)},scale=.5]
\Triangle{$+$}{$-$}{$+$};
\end{scope}

\begin{scope}[shift={(310:\r)},scale=.5]
\Triangle{$\bigcirc$}{$-$}{$+$};
\end{scope}

\begin{scope}[shift={(280:\r)},scale=.5]
\Triangle{$-$}{$-$}{$+$};
\end{scope}

\begin{scope}[shift={(250:\r)},scale=.5]
\Triangle{$-$}{$-$}{$\bigcirc$};
\end{scope}

\begin{scope}[shift={(210:\r)},scale=.5]
\Triangle{$-$}{$-$}{$-$};
\end{scope}

\node at (0,0) [circle, fill = red]{};

\end{tikzpicture}
\caption{ The Shi arrangement for type $A_2$. See \cite[Page
    102]{shiBook}. Each region has been labeled with its sign
  type. See Example~\ref{ex:shiArrWithST}}
\label{fig:shiArrWithST}
\end{figure}

\begin{example}
  \label{ex:shiArrWithST}
All alcoves in the region labeled by
\begin{center}
  \begin{tikzpicture}[scale=.3,
      shift={(0,0)}]\Triangle{$\bigcirc$}{$+$}{$-$}\end{tikzpicture}
  \end{center}
  in Figure~\ref{fig:shiArrWithST} have
  positive coordinate $k_{12}=k_{\alpha_1}$, have negative coordinate
  $k_{\alpha_2}$, and have the coordinate $k_{\theta}$ equal to
  zero.
  Likewise, all three coordinates of all alcoves in the region labeled by
  \begin{center}
  \begin{tikzpicture}[scale=.3,
      shift={(0,0)}]\Triangle{$+$}{$+$}{$+$}\end{tikzpicture}
  \end{center}
are positive.
  See Section~\ref{sec:hyperplanes} for the definition of
  coordinates of an alcove.
 \end{example}

We mention here that in the case where the
Coxeter graph of the system contains an edge with a label greater than
$3$, it is not true that all the left-cells of the affine Weyl
group are unions of admissible sign types. It may be conjectured that it holds for any affine Weyl group of simply-laced
type. The cells in affine $D_4$ have been explicitly described by Shi
in \cite{shi1994}, so the conjecture may not be difficult to verify. It is known
that any left-cell in the lowest or highest two-sided cell of any
irreducible affine Weyl group forms a single admissible sign type; see
Shi \cite{shi1987c,shi1988}. We thank Jian-Yi Shi for this information. 

By a {\em dominant region} of the Shi arrangement, we mean a connected
component of the hyperplane arrangement complement $V \setminus
\bigcup_{H \in \Shi{\roots}} H$ that is contained in the dominant
chamber. Both the formula for the number of regions in the whole
arrangement and for the number of dominant regions are intriguing and
will be discussed in the next section.

\section{Enumeration}
\label{sec:enumeration}

The Shi regions have been counted multiple times. We discuss four different approaches to enumerating them. 
\subsection{The number of Shi regions, part 1}\label{sec:shiEnum}
  Shi concentrated on the admissible sign types in Chapter 7 of his
  book \cite{shiBook}, where he introduced the arrangement for type
  $A$.

He enumerates them for type $A$ by considering the alcove closest to
the origin in each region. We'll call this the {\em minimal alcove} of
the region and denote it $\A_R$ if the region is $R$. Shi called such an alcove
the shortest alcove \cite[Section 7.3]{shiBook}. He characterized
$\A_R$ using left descents. Left descents are key to the definition of
cells, so it is not surprising that they appear in the description of
minimal alcoves. Basically, an alcove is minimal if any reflection
which brings it closer to the origin flips it out of the region.

\begin{figure}[ht]

\begin{tikzpicture}[fill opacity=.5]
\tikzstyle{every node}=[font=\footnotesize]

[-1,-2,1]
} 

\filldraw [draw=white, thick,fill = white, shift={(.0*0.375,-2*0.649519)}](0,0) -- (60:0.75cm) -- (120:0.75cm) --(0,0) node[green,fill opacity=1] at (0cm,0.6cm) {} node[black,fill opacity=1, draw opacity = 1] at (0cm,0.45cm) {$021$} node[green,fill opacity=1] at (0cm,0.3cm) {};

\filldraw [draw=black,thick,fill opacity = .5, fill = Gainsboro,shift={(0.000000,-1.29904)}](0,0) -- (60:0.75cm) -- (0:0.75cm) --(0,0) node[blue] at (0.3cm,0.15cm) {} node[black,fill opacity=1] at (0.375cm,0.3cm) {$2021$} node[blue] at (0.3cm,0.45cm) {};

\filldraw [draw=white, thick,fill = white, shift={(2*0.375,-2*0.649519)}](0,0) -- (60:0.75cm) -- (120:0.75cm) --(0,0) node[green,fill opacity=1] at (0cm,0.6cm) {} node[black,fill opacity=1, draw opacity = 1] at (0cm,0.45cm) {$201$} node[green,fill opacity=1] at (0cm,0.3cm) {};

[-1,-1,-1]
\filldraw [draw=black, thick,fill = Gainsboro, shift={(-0.375,-0.649519)}](0,0) -- (60:0.75cm) -- (120:0.75cm) --(0,0) node[green,fill opacity=1] at (0cm,0.6cm) {} node[black,fill opacity=1, draw opacity = 1] at (0cm,0.45cm) {$121$} node[green,fill opacity=1] at (0cm,0.3cm) {};
[-1,-1,0]
\filldraw [draw=black,thick,fill opacity = .5, fill = Gainsboro,shift={(-0.375000,-0.649519)}](0,0) -- (60:0.75cm) -- (0:0.75cm) --(0,0) node[blue] at (0.3cm,0.15cm) {} node[black,fill opacity=1] at (0.3cm,0.3cm) {$21$} node[blue] at (0.3cm,0.45cm) {};
[0,-1,0]
\filldraw [draw=black, thick,fill = Gainsboro, shift={(0.375,-0.649519)}](0,0) -- (60:0.75cm) -- (120:0.75cm) --(0,0) node[green,fill opacity=1] at (0cm,0.6cm) {} node[black,fill opacity=1, draw opacity = 1] at (0cm,0.45cm) {$1$} node[green,fill opacity=1] at (0cm,0.3cm) {};
[0,-1,1]
\filldraw [draw=black,thick,fill opacity = .5, fill = Gainsboro,shift={(0.375000,-0.649519)}](0,0) -- (60:0.75cm) -- (0:0.75cm) --(0,0) node[blue] at (0.3cm,0.15cm) {} node[black,fill opacity=1] at (0.3cm,0.3cm) {$01$} node[blue] at (0.3cm,0.45cm) {};
[1,-1,1]
\filldraw [draw=black, thick,fill = Gainsboro, shift={(1.125,-0.649519)}](0,0) -- (60:0.75cm) -- (120:0.75cm) --(0,0) node[green,fill opacity=1] at (0cm,0.6cm) {} node[black,fill opacity=1, draw opacity = 1] at (0cm,0.45cm) {$101$} node[green,fill opacity=1] at (0cm,0.3cm) {};

\filldraw [draw=white, thick,fill = white, shift={(-2.*0.375,0*0.649519)}](0,0) -- (60:0.75cm) -- (120:0.75cm) --(0,0) node[green,fill opacity=1] at (0cm,0.6cm) {} node[black,fill opacity=1, draw opacity = 1] at (0cm,0.45cm) {$012$} node[green,fill opacity=1] at (0cm,0.3cm) {};

[-1,0,-1]
\filldraw [draw=black,thick,fill opacity = .5, fill = Gainsboro,shift={(-0.750000,0)}](0,0) -- (60:0.75cm) -- (0:0.75cm) --(0,0) node[blue] at (0.3cm,0.15cm) {} node[black,fill opacity=1] at (0.3cm,0.3cm) {$12$} node[blue] at (0.3cm,0.45cm) {};
[0,0,-1]
\filldraw [draw=black, thick,fill = Gainsboro, shift={(0,0)}](0,0) -- (60:0.75cm) -- (120:0.75cm) --(0,0) node[green,fill opacity=1] at (0cm,0.6cm) {} node[black,fill opacity=1, draw opacity = 1] at (0cm,0.45cm) {$2$} node[green,fill opacity=1] at (0cm,0.3cm) {};
[0,0,0]
\filldraw [draw=black,thick,fill opacity = .5, fill = Goldenrod,shift={(0.000000,0)}](0,0) -- (60:0.75cm) -- (0:0.75cm) --(0,0) node[blue] at (0.3cm,0.15cm) {} node[black,fill opacity=1] at (0.3cm,0.3cm) {$e$} node[blue] at (0.3cm,0.45cm) {};
[1,0,0]
\filldraw [draw=black, thick,fill = Goldenrod, shift={(0.75,0)}](0,0) -- (60:0.75cm) -- (120:0.75cm) --(0,0) node[green,fill opacity=1] at (0cm,0.6cm) {} node[black,fill opacity=1, draw opacity = 1] at (0cm,0.45cm) {$0$} node[green,fill opacity=1] at (0cm,0.3cm) {};
[1,0,1]
\filldraw [draw=black,thick,fill opacity = .5, fill = Goldenrod,shift={(0.750000,0)}](0,0) -- (60:0.75cm) -- (0:0.75cm) --(0,0) node[blue] at (0.3cm,0.15cm) {} node[black,fill opacity=1] at (0.3cm,0.3cm) {$10$} node[blue] at (0.3cm,0.45cm) {};

\filldraw [draw=white, thick,fill = white, shift={(4*0.375,0*0.649519)}](0,0) -- (60:0.75cm) -- (120:0.75cm) --(0,0) node[green,fill opacity=1] at (0cm,0.6cm) {} node[black,fill opacity=1, draw opacity = 1] at (0cm,0.45cm) {$210$} node[green,fill opacity=1] at (0cm,0.3cm) {};

[-1,1,-2]
\filldraw [draw=black,thick,fill opacity = .5, fill = Gainsboro,shift={(-1.125000,0.649519)}](0,0) -- (60:0.75cm) -- (0:0.75cm) --(0,0) node[blue] at (0.3cm,0.15cm) {} node[black,fill opacity=1] at (0.375cm,0.3cm) {$1012$} node[blue] at (0.3cm,0.45cm) {};

\filldraw [draw=white, thick,fill = white, shift={(-1*0.375,0.649519)}](0,0) -- (60:0.75cm) -- (120:0.75cm) --(0,0) node[green,fill opacity=1] at (0cm,0.6cm) {} node[black,fill opacity=1, draw opacity = 1] at (0cm,0.45cm) {$102$} node[green,fill opacity=1] at (0cm,0.3cm) {};

[0,1,-1]
\filldraw [draw=black,thick,fill opacity = .5, fill = Gainsboro,shift={(-1*0.375000,0.649519)}](0,0) -- (60:0.75cm) -- (0:0.75cm) --(0,0) node[blue] at (0.3cm,0.15cm) {} node[black,fill opacity=1] at (0.3cm,0.3cm) {$02$} node[blue] at (0.3cm,0.45cm) {};
[1,1,-1]
\filldraw [draw=black, thick,fill = Gainsboro, shift={(0.375,0.649519)}](0,0) -- (60:0.75cm) -- (120:0.75cm) --(0,0) node[green,fill opacity=1] at (0cm,0.6cm) {} node[black,fill opacity=1, draw opacity = 1] at (0cm,0.45cm) {$202$} node[green,fill opacity=1] at (0cm,0.3cm) {};
[1,1,0]
\filldraw [draw=black,thick,fill opacity = .5, fill = Goldenrod,shift={(0.375000,0.649519)}](0,0) -- (60:0.75cm) -- (0:0.75cm) --(0,0) node[blue] at (0.3cm,0.15cm) {} node[black,fill opacity=1] at (0.4cm,0.3cm) {$20$} node[blue] at (0.3cm,0.45cm) {};

\filldraw [draw=white, thick,fill = white, shift={(3*0.375,0.649519)}](0,0) -- (60:0.75cm) -- (120:0.75cm) --(0,0) node[green,fill opacity=1] at (0cm,0.6cm) {} node[black,fill opacity=1, draw opacity = 1] at (0cm,0.45cm) {$120$} node[green,fill opacity=1] at (0cm,0.3cm) {};

[2,1,1]
\filldraw [draw=black,thick,fill opacity = .5, fill = Goldenrod,shift={(1.125000,0.649519)}](0,0) -- (60:0.75cm) -- (0:0.75cm) --(0,0) node[blue] at (0.3cm,0.15cm) {} node[black,fill opacity=1] at (0.375cm,0.3cm) {$1210$} node[blue] at (0.3cm,0.45cm) {};

\def\b{1}
\def\b{1}
\def\ma{1.94856/1.125}
\def\xa{-2.25+\b+1.94856}
\def\ya{\xa*\ma}
\def\ya{1.2099}
\def\dy{.5}
\def\dx{dy/\ma}
\def\dx{.2886747137}

\path [draw = red, very thick,  shift = {(0,0)}](-1.125+\dx,-1.94856+\dy) -- (1.125,1.94856);
\path [draw = red,very thick, shift = {(0,0)}](-2.25+\b,0) -- (2.25,0);
\path [draw = red, very thick,  shift = {(0,0)}](-1.125,1.94856) -- (1.125-\dx,-1.94856+\dy);
\path [draw = red, very thick,  shift = {(0.75,0)}](-1.125+\dx,-1.94856+\dy) -- (1.125,1.94856);
\path [draw = red,very thick, shift = {(0,0.649519)}](-2.25+\b,0) -- (2.25,0);
\path [draw = red, very thick,  shift = {(0.75,0)}](-1.125,1.94856) -- (1.125-\dx,-1.94856+\dy);
\fill[fill=black, fill opacity=1.](0.,0.) circle (0.0675cm);[-2,-4,2]

\end{tikzpicture}
  \caption{\small On the left, the minimal alcove of each Shi region
    is shaded. The dominant ones are shaded yellow, the rest
    gray. We've labeled the alcoves by the corresponding affine
    permutation, where we use $i$ for $s_i$ for space reasons. The
    alcoves corresponding to the inverses are drawn on the right. }
  \label{fig:inverses}
\end{figure}
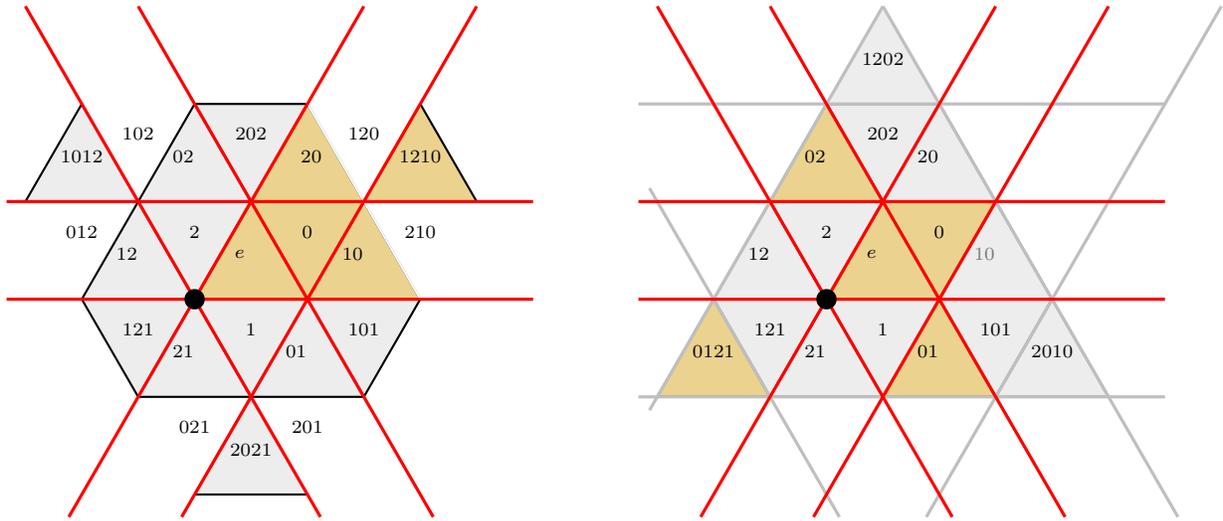

\begin{example}
  See Figure~\ref{fig:inverses}. Let $w=s_1s_2s_1s_0$. This $w$ has
  two left descents, $s_1$ and $s_2$, and both $s_1w=s_2s_1s_0$ and
  $s_2w=s_1s_2s_0$ are in different regions the minimal permutation
  $w$. On the other hand, $w=s_1s_2s_0$ is not minimal and indeed
  $\len{s_1w}<\len{w}$ and $s_1w=s_2s_0$ is in the same region as $w$.
  \end{example}

Every alcove corresponds to an affine
permutation. We'll call the affine
permutations whose alcoves are minimal {\em minimal permutations}.
Shi showed that the collection of alcoves corresponding to the
inverses of minimal permutations is exactly a scaled version of the
fundamental alcove. See Figure~\ref{fig:inverses}. Thus to calculate
the number of regions in his newfound arrangement, he calculated the
number of alcoves in this scaled fundamental alcove. He calculated
something a bit more general: if the fundamental alcove is expanded by
the positive integer $m$, then it is made up of $m^{\dim(V)}$
alcoves. The alcoves corresponding to the inverses of minimal
permutation land in the fundamental alcove scaled by $n+1$, which
showed that there are $(n+1)^{n-1}$ regions in the Shi arrangement of
type $A$. The expression $(n+1)^{n-1}$ pops up frequently in
combinatorics and algebra; see \cite{haglund2008} for their connection
to $q,t$-Catalan numbers, for example. 


Shi's enumeration of the regions is perhaps more complicated than
the others described here. However, his discovery that the inverses of
the minimal permutations correspond to a simplex is worth the price of
admission. The minimal alcoves have been useful in other enumeration;
see \cite{athanasiadis2005,fishel-vazirani2010}, for example. For
another example, Hohlweg, Nadeau, and Williams, in
\cite{hohlweg-nadeau-williams2016}, generalize the Shi arrangement to
any Coxeter group (and beyond!) and conjecture that the inverses of
the analogues of minimal permutations form a convex body. See also
Sommers \cite{sommers2005} where the simplex was generalized to what
is now called the Sommers region. Thomas and Williams
\cite{thomas-williams2014} show that the set of alcoves in this
region, and by extension the Shi regions, exhibit the cyclic sieving
phenomenon.

\begin{figure}[h]
\begin{tikzpicture}[fill opacity=.5, scale=2, font=\scriptsize]


\begin{scope}[shift={(4.2,0)}] 
\filldraw [draw=black,thick,fill opacity = .5, fill = Goldenrod,shift={(-1.125000,-0.649519)}](0,0) -- (60:0.75cm) -- (0:0.75cm) --(0,0) node[blue] at (0.3cm,0.15cm) {} node[black,fill opacity = 1] at (0.375cm,0.3cm) {$0121$} node at (0.3cm,0.45cm) {};

\filldraw [draw=black, thick,fill = Gainsboro, shift={(-0.375,-0.649519)}](0,0) -- (60:0.75cm) -- (120:0.75cm) --(0,0) node[green,fill opacity=1] at (0cm,0.6cm) {} node[black,fill opacity=1, draw opacity = 1] at (0cm,0.45cm) {$121$} node[green,fill opacity=1] at (0cm,0.3cm) {};
[-1,-1,0]
\filldraw [draw=black,thick,fill opacity = .5, fill = Gainsboro,shift={(-0.375000,-0.649519)}](0,0) -- (60:0.75cm) -- (0:0.75cm) --(0,0) node[blue] at (0.3cm,0.15cm) {} node[black,fill opacity = 1] at (0.3cm,0.3cm) {$21$} node[blue] at (0.3cm,0.45cm) {};
[0,-1,0]
\filldraw [draw=black, thick,fill = Gainsboro, shift={(0.375,-0.649519)}](0,0) -- (60:0.75cm) -- (120:0.75cm) --(0,0) node[green,fill opacity=1] at (0cm,0.6cm) {} node[black,fill opacity=1, draw opacity = 1] at (0cm,0.45cm) {$1$} node[green,fill opacity=1] at (0cm,0.3cm) {};
[0,-1,1]
\filldraw [draw=black,thick,fill opacity = .5, fill = Goldenrod,shift={(0.375000,-0.649519)}](0,0) -- (60:0.75cm) -- (0:0.75cm) --(0,0) node[blue] at (0.3cm,0.15cm) {} node[black,fill opacity = 1] at (0.3cm,0.3cm) {$01$} node [blue] at (0.3cm,0.45cm) {};
\begin{pgfonlayer}{foreground layer}
\end{pgfonlayer}{foreground layer}
[1,-1,1]
\filldraw [draw=black, thick,fill = Gainsboro, shift={(1.125,-0.649519)}](0,0) -- (60:0.75cm) -- (120:0.75cm) --(0,0) node[green,fill opacity=1] at (0cm,0.6cm) {} node[black,fill opacity=1, draw opacity = 1] at (0cm,0.45cm) {$101$} node[green,fill opacity=1] at (0cm,0.3cm) {};
[1,-1,2]
\filldraw [draw=black,thick,fill opacity = .5, fill = Gainsboro,shift={(1.125000,-0.649519)}](0,0) -- (60:0.75cm) -- (0:0.75cm) --(0,0) node[blue] at (0.3cm,0.15cm) {} node[black,fill opacity = 1] at (0.375cm,0.3cm) {$2010$} node[blue] at (0.3cm,0.45cm) {};
[-1,0,-1]
\filldraw [draw=black,thick,fill opacity = .5, fill = Gainsboro,shift={(-0.750000,0)}](0,0) -- (60:0.75cm) -- (0:0.75cm) --(0,0) node[blue] at (0.3cm,0.15cm) {} node[black,fill opacity = 1] at (0.3cm,0.3cm) {$12$} node[blue] at (0.3cm,0.45cm) {};
[0,0,-1]
\filldraw [draw=black, thick,fill = Gainsboro, shift={(0,0)}](0,0) -- (60:0.75cm) -- (120:0.75cm) --(0,0) node[green,fill opacity=1] at (0cm,0.6cm) {} node[black,fill opacity=1, draw opacity = 1] at (0cm,0.45cm) {$2$} node[green,fill opacity=1] at (0cm,0.3cm) {};
[0,0,0]
\filldraw [draw=black,thick,fill opacity = .5, fill = Goldenrod,shift={(0.000000,0)}](0,0) -- (60:0.75cm) -- (0:0.75cm) --(0,0) node[blue] at (0.3cm,0.15cm) {} node[black,fill opacity=1,draw opacity = 1] at (0.3cm,0.3cm) {$e$} node[blue] at (0.3cm,0.45cm) {};
\begin{pgfonlayer}{foreground layer}
\end{pgfonlayer}{foreground layer}
[1,0,0]
\filldraw [draw=black, thick,fill = Goldenrod, shift={(0.75,0)}](0,0) -- (60:0.75cm) -- (120:0.75cm) --(0,0) node[green,fill opacity=1] at (0cm,0.6cm) {} node[black,fill opacity=1, draw opacity = 1] at (0cm,0.45cm) {$0$} node[green,fill opacity=1] at (0cm,0.3cm) {};
\begin{pgfonlayer}{foreground layer}
\end{pgfonlayer}{foreground layer}
[1,0,1]
\filldraw [draw=black,thick,fill opacity = .5, fill = Gainsboro,shift={(0.750000,0)}](0,0) -- (60:0.75cm) -- (0:0.75cm) --(0,0) node[blue] at (0.3cm,0.15cm) {} node[black,draw opacity = 1] at (0.3cm,0.3cm) {$10$} node[blue] at (0.3cm,0.45cm) {};
[0,1,-1]
\filldraw [draw=black,thick,fill opacity = .5, fill = Goldenrod,shift={(-0.375000,0.649519)}](0,0) -- (60:0.75cm) -- (0:0.75cm) --(0,0) node[blue] at (0.3cm,0.15cm) {} node[black,fill opacity = 1] at (0.3cm,0.3cm) {$02$} node[blue] at (0.3cm,0.45cm) {};
\begin{pgfonlayer}{foreground layer}
\end{pgfonlayer}{foreground layer}
[1,1,-1]
\filldraw [draw=black, thick,fill = Gainsboro, shift={(0.375,0.649519)}](0,0) -- (60:0.75cm) -- (120:0.75cm) --(0,0) node[green,fill opacity=1] at (0cm,0.6cm) {} node[black,fill opacity=1, draw opacity = 1] at (0cm,0.45cm) {$202$} node[green,fill opacity=1] at (0cm,0.3cm) {};
[1,1,0]
\filldraw [draw=black,thick,fill opacity = .5, fill = Gainsboro,shift={(0.375000,0.649519)}](0,0) -- (60:0.75cm) -- (0:0.75cm) --(0,0) node[blue] at (0.3cm,0.15cm) {} node[black,fill opacity = 1] at (0.3cm,0.3cm) {$20$} node[blue] at (0.3cm,0.45cm) {};

[1,2,-1]

\fill [fill opacity = .5, fill = Gainsboro,shift={(0.000000,1.29904)}](0,0) -- (60:0.75cm) -- (0:0.75cm) --(0,0) node[blue] at (0.3cm,0.15cm) {} node[black,fill opacity = 1] at (0.375cm,0.3cm) {$1202$} node[blue] at (0.3cm,0.45cm) {};

\def\b{1}
\def\ma{1.94856/1.125}
\def\xa{-2.25+\b+1.94856}
\def\ya{\xa*\ma}
\def\ya{1.2099}
\def\dy{.5}
\def\dx{dy/\ma}
\def\dx{.2886747137}

\path [draw = gray!50, very thick,  shift = {(2*0.75,0)}](-1.125+\dx,-1.94856+\dy) -- (1.125,1.94856);
\path [draw = gray!50,very thick, shift = {(0,2*0.649519)}](-2.25+\b,0) -- (2.25,0);
\path [draw = gray!50, very thick,  shift = {(2*0.75,0)}](-1.125,1.94856) -- (1.125-\dx,-1.94856+\dy);

\path [draw = gray!50, very thick,  shift = {(-1*0.75,0)}](-1.125+\xa,-1.94856+\ya) -- (1.125,1.94856);
\path [draw = gray!50,very thick, shift = {(0,-1*0.649519)}](-2.25+\b,0) -- (2.25,0);
\path [draw = gray!50, very thick,  shift = {(-1*0.75,0)}](-1.125+\xa,1.94856-\ya) -- (1.125-\dx,-1.94856+\dy);

\path [draw = red, very thick,  shift = {(0,0)}](-1.125+\dx,-1.94856+\dy) -- (1.125,1.94856);
\path [draw = red,very thick, shift = {(0,0)}](-2.25+\b,0) -- (2.25,0);
\path [draw = red, very thick,  shift = {(0,0)}](-1.125,1.94856) -- (1.125-\dx,-1.94856+\dy);

\path [draw = red, very thick,  shift = {(0.75,0)}](-1.125+\dx,-1.94856+\dy) -- (1.125,1.94856);
\path [draw = red,very thick, shift = {(0,0.649519)}](-2.25+\b,0) -- (2.25,0);
\path [draw = red, very thick,  shift = {(0.75,0)}](-1.125,1.94856) -- (1.125-\dx,-1.94856+\dy);
\def\g{20}

\fill[fill=black, fill opacity=1.](0.,0.) circle (0.0675cm);[-4,-2,-2]
\end{scope} 

\omitt{
\end{tikzpicture}
  \caption{\small 
$w \Afund$ for the $\m$-minimal alcoves $w^{-1} \Afund$ in Figure
\ref{fig;onetwoshi}
below, $m=1,2$.  Note $w \Afund \subseteq \Region$. 
Each $\gamma \in Q$ is in precisely one yellow/blue alcove,
so this  illustrates the second statement of Theorem \ref{thm;haiman}. }
  \label{fig;dilate}
\end{figure}

In \cite{shi1987}, Shi generalized sign types to other affine Weyl
groups. He defined sets analogous to \eqref{GA} for other types. The
hyperplane arrangements were still given by Definition~\ref{def:shi},
but now he considered root systems other than type $A$. He used the
map $\zeta$ on alcoves and described the sign types which arose and again
characterized the element in each region with the minimal number of
hyperplanes separating it from the origin. As above, we identify
elements $w\in W$ with $\Afund w$ and refer to $w$ as minimal if its
alcove is minimal.  The fact that
$$\bigcup_{w\text{ minimal}}\Afund w^{-1}$$ is a simplex is not just a
type $A$ phenomenon. Shi proved it for other affine Weyl groups and
used it to prove that there are $(h+1)^n$ regions, where $h$ is the Coxeter
number of the system.

 Shi counted the number of regions in the dominant chamber for affine
 Weyl groups in \cite{shi1997}.\omitt{$A_{\ell}$, $B_m$, $C_m$, $D_n$,
   $E_6$, $E_7$, $E_8$, $F_4$, $G_2$} He calls the admissible sign
 types corresponding to regions in the dominant chamber $\oplus$-sign
 types. For types $A$, $B$, $C$, and $D$, he finds a bijection from
 $\oplus$-sign types to filters in the root poset for $\pos$. Here is
 a technical detail: Shi finds the bijection to the positive coroots
 $(\pos)^{\vee}$, which we won't define, then mentions that it has the
 same type as $\pos$ except when $\pos$ has type $B$ and $C$. He deals
 with types $B$ and $C$ separately. We will continue using $\pos$. He
 further maps the filters to subdiagrams of certain Young diagrams;
 see Section~\ref{subsec:rootposet}. For example, for type $A$, the
 subdiagrams are those of partitions whose diagrams fit inside the
 staircase shape. In the exceptional cases, he enumerates increasing
 subsets directly. He shows the $\oplus$-sign types for affice Weyl groups are enumerated by
 the Catalan numbers, although Shi does not mention them. 

A few more words are in order on this important bijection to filters
in the root poset. The key proposition, from Section 1.2 of
\cite{shi1997}, follows (using roots instead of coroots).

\begin{proposition}
  \label{prop:keyoplus}
  Assume that $X=(X_{\alpha})_{\alpha\in\pos}$ is a $\pos$-tuple with
  $X_{\alpha}\in\{+,\bigcirc\}$. Then $X$ is an $\oplus$-sign type if
  and only if the following condition on $X$ holds: if
  $\alpha,\beta\in\pos$ satisfy $\beta>\alpha$ and $X_{\alpha}=+$, then
    $X_{\beta}=+$.
  \end{proposition}

We'll use type $A$ as an example. In the set $G_A$ in \eqref{GA},
there are five triples which contain only $\bigcirc$ and $+$, but
there are eight which are possible. The condition in
Proposition~\ref{prop:keyoplus} rules the other three out, proving
sufficiency. For necessity, Shi uses induction to reduce to the rank
two case and shows that the five $\oplus$-sign types in $G$ satisfy
the condition.

\subsection{Interlude}
We'll need these standard definitions for
Sections~\ref{subsec:headleyEnum} and \ref{subsec:athanEnum}. See
\cite{EC1} for the definitions of the rank function $\rho$ and M\"obius function $\mu$ of a
poset.
\begin{definition}\cite{EC1}
Let $P$ be a finite graded poset with $\hat{0}$. Let $\rho$ be its
rank function and $n$ the rank of $P$. Define the characteristic
polynomial $\chi_P(x)$ of $P$ by $$\chi_P(x)=\sum_{t\in
  P}\mu(\hat{0},t)x^{n-\rho(t)}.$$
\end{definition}

\begin{definition}\cite{EC1}
  Let $\mc{A}$ be a hyperplane arrangement in a vector space $V$ and
  let $L(\mc{A})$ be the set of all nonempty intersections of
  hyperplanes in $\mc{A}$. Include $V$ itself, by considering it as
  the intersection over the empty set. Order $L(\mc{A})$ by reverse
  inclusion, so that $\hat{0}$ is $V$.
  \end{definition}

\begin{figure}[h]
\begin{tikzpicture}

\def\a{3}
\def\b{1.5}
\def\h{.5}
\begin{scope}[font=\scriptsize,every node/.append style={text width=,  minimum size=0cm, inner sep=1mm,fill=blue!20,rectangle, rounded corners}]
\node (i1) at (-2.5*\a,2.7*\b)  {$x_1=x_2=x_3$};
\node (i2) at (-1.5*\a,2.7*\b) [] {$x_1=x_2=x_3+1$};
\node (i3) at (-.5*\a,2.7*\b) [] {$x_1=x_2+1=x_3$};
\node (i4) at (.5*\a,2.7*\b) [] {$x_1-1=x_2=x_3$};
\node (i5) at (1.5*\a,2.7*\b) [] {$x_1-1=x_2=x_3+1$};
\node (i6) at (2.5*\a,2.7*\b) [] {$x_1=x_2-1=x_3$};

\node (p1) at (-2.5*\a,1*\b) {$x_1=x_2$};
\node (p2) at (-1.5*\a,1*\b) {$x_1=x_2+1$};
\node (p3) at (-.5*\a,1*\b) {$x_1=x_3$};
\node (p4) at (.5*\a,1*\b) {$x_1=x_3+1$};
\node (p5) at (1.5*\a,1*\b) {$x_2=x_3$};
\node (p6) at (2.5*\a,1*\b) {$x_2=x_3+1$};

\node (v) at (0,0) {$V$};
\end{scope}

\begin{pgfonlayer}{foreground layer}
\begin{scope}[font=\scriptsize,every node/.append style={text width=,  minimum size=0cm, inner sep=1mm,fill=orange!30,circle}]
\node (i1m) at (-2.5*\a,2.7*\b+\h)  {$+2$};
\node (i2m) at (-1.5*\a,2.7*\b+\h) [] {$+2$};
\node (i3m) at (-.5*\a,2.7*\b+\h) [] {$+1$};
\node (i4m) at (.5*\a,2.7*\b+\h) [] {$+1$};
\node (i5m) at (1.5*\a,2.7*\b+\h) [] {$+2$};
\node (i6m) at (2.5*\a,2.7*\b+\h) [] {$+1$};

\node (p1m) at (-2.5*\a,1*\b-\h) {$-1$};
\node (p2m) at (-1.5*\a,1*\b-\h) {$-1$};
\node (p3m) at (-.5*\a,1*\b-\h) {$-1$};
\node (p4m) at (.5*\a,1*\b-\h) {$-1$};
\node (p5m) at (1.5*\a,1*\b-\h) {$-1$};
\node (p6m) at (2.5*\a,1*\b-\h) {$-1$};

\node (vm) at (0,0-\h) {$+1$};
\end{scope}
\end{pgfonlayer}{foreground layer}

\draw[thick, red] (v)--(p1);
\draw[thick, red] (v)--(p2);
\draw[thick, red] (v)--(p3);
\draw[thick, red] (v)--(p4);
\draw[thick, red] (v)--(p5);
\draw[thick, red] (v)--(p6);

\draw[thick, red] (p1)--(i1);
\draw[thick, red] (p1)--(i2);
\draw[thick, red] (p2)--(i3);
\draw[thick, red] (p2)--(i4);
\draw[thick, red] (p2)--(i5);
\draw[thick, red] (p3)--(i1);
\draw[thick, red] (p3)--(i3);
\draw[thick, red] (p3)--(i6);
\draw[thick, red] (p4)--(i2);
\draw[thick, red] (p4)--(i5);
\draw[thick, red] (p5)--(i1);
\draw[thick, red] (p5)--(i4);
\draw[thick, red] (p6)--(i2);
\draw[thick, red] (p6)--(i5);
\draw[thick, red] (p6)--(i6);
\end{tikzpicture}

  \caption{The poset of intersections for type $A_2$.  For example, the
  node labeled $x_1=x_2=x_3+1$ represents
  $\Hak{\alpha_1}{0}\cap\Hak{\alpha_2}{0}\cap\Hak{\alpha_3}{1}$. See
  \cite[Page 36]{headley1994}. The numbers in the circle above or
  below each element is the M\"obius function of the interval from the
  element to $V=\hat{0}$. The Poincar\'e polynomial is
  $\poin{\roots_{A_2}}(x)=(1+3x)^2$.}
\label{fig:posetOfInts}
\end{figure}
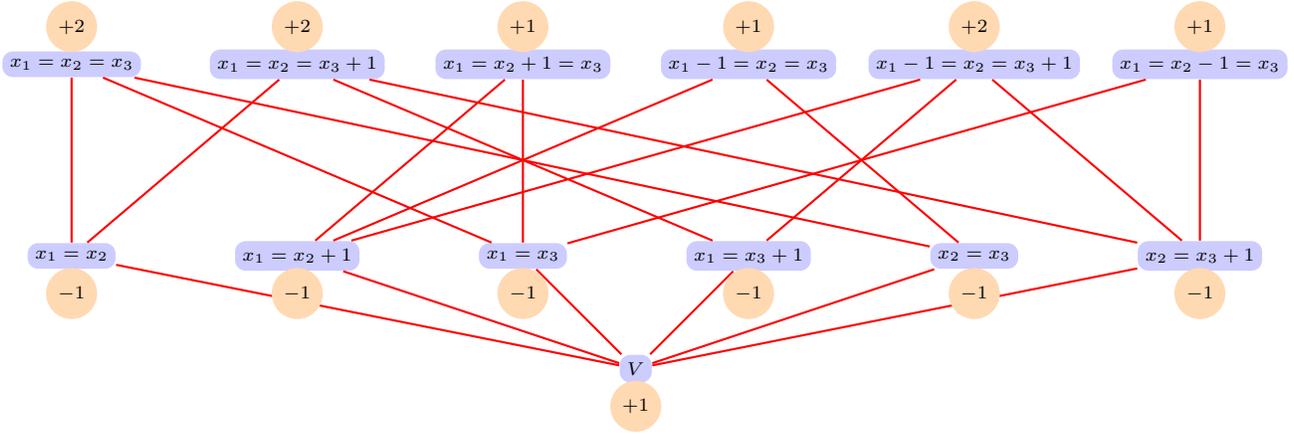

See Figure~\ref{fig:posetOfInts} for an example of the poset of intersections.

If the intersection of all the hyperplanes in $\mc{A}$ is nonempty,
then $L(\mc{A})$ is a lattice. The intersection of all hyperplanes in
$\Shi{\roots}$ is empty and $L(\Shi{\roots})$ will only be a meet
semi-lattice. It is finite and graded by $\rho(t)=n-\dim(t)$, where
$n=\dim(V)$.  The characteristic polynomial of an arrangement is
$$\chi_{\mc{A}}(x)=\sum_{t\in L(\mc{A})}\mu(\hat{0},t)x^{n-\rho(t)}$$
and the Poincar\'e polynomial is $$\poin{\mc{A}}(x)=\sum_{t\in
  L(\mc{A})}\mu(\hat{0},t)(-x)^{\rho(t)}.$$ The Poincar\'e polynomial
is a rescaled version of the characteristic polynomial.


The characteristic polynomial is invaluable for studying hyperplane
arrangements, thanks to a theorem of Zaslavsky \cite{zaslavsky1975,
  zaslavsky1975a}. See also Stanley's notes on hyperplanes
\cite{stanley2007}. In Section~\ref{sec:hyperplanes}, we defined
$\regNum(\mc{A})$ and $\bddRegNum(\mc{A})$ be the number of regions
and number of bounded regions of the arrangement $\mc{A}$.

\begin{theorem}[\cite{zaslavsky1975}]
  Let $\mc{A}$ be an arrangement in an $n$-dimensional real vector space. Then
\begin{eqnarray*}
\regNum(\mc{A})&=&(-1)^n\chi_{\mc{A}}(-1)=\poin{\mc{A}}(1)\\
\bddRegNum(\mc{A})&=&(-1)^{\rank{\mc{A}}}\chi_{\mc{A}}(1)=\poin{\mc{A}}(-1).\\
\end{eqnarray*}
\end{theorem}

\subsection{The number of Shi regions, part 2}\label{subsec:headleyEnum}
Headley \cite{headley1997,
headley1994} calculated the Poincar\'e polynomial $\poin{\Shi{\roots}}$
of the Shi arrangement for an irreducible root system.  He found a
recursion for its coefficients, which we will now present.

Let $\overline{\Shi{\roots}}$ be the subarrangement of $\Shi{\roots}$
consisting of all hyperplanes which contain the origin. For $Y\in
L(\overline{\Shi{\roots}})$, let $W_Y$ be the group generated by the
reflection through all the hyperplanes containing $Y$. For a polynomial $p(t)$, let $[t^k]p(t)$ be the coefficient of
$t^k$ in $p(t)$.

\begin{lemma}[\cite{headley1997}]\label{lem:H}
For $Y\in L(\overline{\Shi\roots})$, let $W_{Y,1}\times\cdots W_{Y,m}$ be the decomposition of $W_{Y}$ into irreducible Coxeter groups. Let $S_i=\Shi{W_{Y,i}}$ be the Shi arrangement associated to the Coxeter group $W_{Y,i}$. Then
$$[t^k]\poin{\Shi{\roots}}(t)=[t^k]\sum_{Y\in L(\overline{\Shi{\roots}}):\rank(Y)=k}\poin{S_1}(t)\cdots\poin{S_m}(t).$$
\end{lemma}

Headley uses induction on the number of generators to determine every coefficient except the leading one. For this, he relies on Shi's enumeration and the relationship between the Poincar\'e polynomial and the number of regions. His analysis is done case by case for each Coxeter type. His theorem is

\begin{theorem}
Let $\roots$ be an irreducible crystallographic root system, with Coxeter number $h$ and rank $n$. Then $$\poin{\Shi{\roots}}(t)=(1+ht)^n.$$
\end{theorem}

In type $A$, the argument for calculating the coefficients is simple enough to repeat: he matches an element $Y$ of the intersection poset $L(\hat{\Shi{\roots}})$ with the set partition $B=(B_1,\ldots,B_m)$ of $[n+1]$ by
$$Y=\cap\{x_i-x_j=0:i,j\text{ are in the same block of }B\}.$$ In this case, $W_Y$ is isomorphic to $A_{|B_1|-1}\times\cdots\times A_{|B_m|-1}$ and $\rank(Y)=n+1-m$. Therefore, by induction and Lemma~\ref{lem:H},

\begin{equation}\label{eq:headleyInd}
[t^k]\poin{\Shi{\roots}}(t)=\sum_{\text{Partitions of $[n+1]$ into
    $n+1-k$
    blocks}}|B_1|^{|B_1|-1}\cdots|B_{n+1-k}|^{|B_{n+1-k}|-1}.\end{equation}
In his thesis\cite{headley1994}, he used Lagrange inversion to
calculate $$[t^k]\poin{\Shi{\roots}}=(n+1)^k\binom{n}{n-k}.$$ Later, in
\cite{headley1997}, he recognized the sum in \eqref{eq:headleyInd} to
be the number of labeled forests on $n+1$ vertices of $n+1-k$ trees
and used \cite{moon1970}. In both his thesis and later paper, he
showed that the coefficient of $t^k$ in $\poin{\Shi{\roots}}(t)$ and
$(1+(n+1))^n$ are the same for $1\leq k\leq n-1$. Then since the
degree of $\poin{\Shi{\roots}}(t)$ is $n$ and since
$\poin{\Shi{\roots}}(1)=(n+2)^n$ by Shi's result, he showed
\begin{equation}
\label{eq:poinA}
\poin{\Shi{\roots}}(t)=(1+(n+1)t)^n=(1+ht)^n
\end{equation}
 in $A_n$.

\subsection{The number of Shi regions, part 3} 
\label{subsec:athanEnum}
Crapo and Rota \cite[Chapter 16]{crapo-rota1970} described the {\em critical
  problem}: let $S$ be a set of points in an $n$-dimensional vector
space $V_n$ over the field $\mathbb{F}_q$ with $q$ elements. The set $S$
must not contain the origin. Find the minimum number $c$ of projective
hyperplanes $H_1,\ldots,H_c$ with the property that the intersection
$H_1\cap\ldots H_c\cap S$ is null. They were able to solve the problem
using the poset of intersections and characteristic
  polynomial. \omitt{We'll see both of these several times, so we might as
well define them now.}

Athanasiadis \cite{athanasiadis1996} turned Crapo and Rota's theorem
around and used it to calculate the characteristic polynomial of
subspace arrangements. Blass and Sagan \cite{blass-sagan1998} had
previously used a similar idea, but not for all subspaces and not for
the Shi arrangement.  We present first the the Crapo and Rota theorem,
then describe how Athanasiadis used it to get his hands on the
characteristic polynomial for the Shi arrangements for irreducible
crystallographic root systems.

\begin{theorem}[\cite{crapo-rota1970}]  The number of linearly
  ordered sequences $(L_1,L_2,\ldots,L_k)$ of $k$ linear functionals
  in $V_n$ which distinguish the set $S$ is given by $p(q^k)$ where
  $p(v)$ is the characteristic polynomial of the geometric lattice
  spanned by the set $S$.
  \end{theorem}

Athanasiadis needed to count the number of $n$-tuples
$(x_1,\ldots,x_n)\in\F_q^n$ which satisfy $x_i\neq x_j$ and $x_i\neq
x_j+1$ for $i<j$. The argument is simple (and lovely) enough in type 
$A_{n-1}$ for the full Shi arrangement that we reproduce it here. See also \cite{EC1}.

We first solve a related problem. Find the number of ways there are to place $n$
labeled balls in $q$ unlabeled boxes, where

\begin{enumerate}
\item the boxes are in a circle,\item \label{cond:ballsInj} there is
  never more than one ball in a box, and \item \label{cond:ballsAff}
  if $i<j$, then ball $i$ is not placed in the box immediately
  following, in the clockwise direction, the box holding ball
  $j$. \end{enumerate}

There will be $q-n$ empty boxes, so first place them in a
circle. There is one way to do that. There are now $q-n$ spaces
between the empty boxes, where the boxes holding the balls will go. By
cyclic symmetry, there is one way to place the box holding the
$1$-ball. Then there are $(q-n)^{n-1}$ ways to place the rest of the
boxes holding balls in the empty spaces. It is enough to pick the
space between empty boxes: to avoid violating condition
\eqref{cond:ballsAff}, the boxes between a consecutive pair of empty
boxes must placed in increasing order of the labels on the balls
inside. That is, our final answer to the related problem is
$(q-n)^{n-1}$.


Now back to counting $n$-tuples. We are essentially done, if we think
of each $n$-tuple representing a distribution of $n$ labeled balls
into a circle of $q$ labeled boxes, where the distribution satisfies
conditions \eqref{cond:ballsInj} and \eqref{cond:ballsAff}. We place
the ball labeled $i$ in box $x_i$. We need only label the boxes, and
there are $q$ ways to do this. Thus there are $q(q-n)^{n-1}$
$n$-tuples which satisfy $x_i\neq x_j$ and $x_i\neq x_j+1$ for $i<j$
and we have that the characteristic polynomial is $\chi_{L(\Shi{n})}=q(q-n)^{n-1}$.

Crapo and Rota's finite field method has since been used to calculate
other characteristic polynomials. See Armstrong \cite{armstrong2013},
Armstrong and Rhoades \cite{armstrong-rhoades2012}, and Ardila
\cite{ardila2007}, for example. See Athanasiadis
\cite{athanasiadis2010} for reciprocity results for the characteristic
polynomial for the Shi arrangement.

\Christos{In Section 5 of \cite{yoshinaga18}, Yoshinaga gave a
  uniform, lattice point counting proof of the formula for the
  characteristic polynomial of the Shi arrangement.}

\Zaslavsky{Berthom\'e, Cordovil, Forge, Ventos, and Zaslavsky use gain graphs in \cite{BCFVZ} to calculate the chromatic polynomial and thereby the characteristic polynomial (and more).}

\Levear{Athanasiadis further used the finite field method to count the
  faces in the Shi arrangement in
  \cite{athanasiadisThesis,athanasiadis1996}. If we view the regions
  of an arrangement as polytopes, the faces of the arrangement are the
  open faces of the these polytopes. The dimension of a face is the
  dimension of its affine span. Athanasiadis proved
  $$f_k=\binom{n}{k}\sum_{i=0}^k(-1)^{i}\binom{k}{i}(n-i+1)^{n-1},$$
  where $f_k$ is the number of $k$-codimensional faces of the Shi
  arrangement for $0\leq k\leq n.$ He then uses inclusion-exclusion to
  write
\begin{equation}\label{eqn:faceNum} f_k=\binom{n}{k}\left|\{f:[n-1]\to[n+1]|[k]\subseteq\im f\}\right|, \end{equation}
  for $0\leq k\leq n-1$ and $f_n=0$.}

\Levear{Levear \cite{levearFPSAC19,levearThesis,levear2020bijections}
  gave a bijective explanation of \eqref{eqn:faceNum}, and a
  generalization for the $m$-Shi arrangement. He also has a
  similar result for the $m$-Catalan arrangement. }

\Levear{Ehrenborg  generalized Athanasiasis's face counting result in \cite{ehrenborg19}.}

\subsection{The number of Shi regions, part 4}
\label{sec:stanEnum}
Pak and Stanley, in \cite{stanley1996}, give a bijection from Shi regions (type $A$) to parking functions, which
refines  \eqref{eq:poinA}.  It is proved to be a bijection in \cite[Lecture
  6]{stanley2007}. A {\em parking function of length $n$} is a tuple
of nonnegative integers $(p_1,\ldots,p_n)$ such that when rearranged in nondecreasing order and relabeled as $b_1\leq b_2\leq\cdots\leq b_n$, then $b_i\leq i-1$. Parking functions generalize inversion vectors of
permutations.  Pak and Stanley recursively defined the label
$\lambda(R)$ of a region $R$. We use the description given in
\cite[Lecture 6]{stanley2007}.

Let $R_0$ be the fundamental alcove $\Afund$.
Set $\lambda(R_0)=(0,\cdots,0)$. Suppose we have labeled the region
$R$ and its label $\lambda(R)$ is $(a_1,\ldots,a_n)$.
\begin{itemize}
\item If the regions $R$ and $R'$ are separated by the single
  hyperplane $H$ with the equation $x_i-x_j=0$, $i<j$, and if $R$ and
  $R_0$ lie on the same side of $H$, then
  $\lambda(R')=(a_1,\ldots,a_{i-1},a_i+1,a_{i+1}\ldots,a_{j-1}a_j,a_{j+1},\dots,a_n)$.
\item If the regions $R$ and $R'$ are separated by the single
  hyperplane $H$ with the equation $x_i-x_j=1$, $i<j$, and if $R$ and
  $R_0$ lie on the same side of $H$, then
  $\lambda(R')=(a_1,\ldots a_{i-1},a_i,a_{i+1}\ldots,a_{j-1},a_j+1,a_{j+1},\dots,a_n)$.
\end{itemize}

The bijection generalizes the well-known bijection from permutations
to inversion vectors \cite[Chapter 1]{EC1}.  Although the map
$\lambda$ is simply stated, the proof that it is a bijection is not
simple. To show that the labeling is a bijection, Stanley encodes each
region as a permutation and antichain pair. He builds the inverse map
step-by-step from the parking functions to the pairs. The summary by
Armstrong \cite[Theorem 3]{armstrong2013} of the proof that the
Pak-Stanley map is bijective is particularly good. Recall that the
filters/antichains/ideals in the root poset for type $A$ correspond to
partitions in a staircase, and define the non-inversions of a
permutation $w$ to be the pairs $(i,j)$ such that $i<j$ and $w(i)<w(j)$.
Then the proof can be summarized as showing that the Shi regions are
in bijection with pairs $(w,\mc{I})$ where $w\in\Sn$ and $\mc{I}$ is
an ideal in the root poset $\pos$ such that the minimal elements of
$\mc{I}$, which are labels in the valleys of the Dyck path
corresponding to $\mc{I}$, are non-inversions of $w$.


The Pak and Stanley bijection from regions to
  parking functions ($m=1$) case can be composed with a bijection from trees to parking functions. \omitt{Additionally} The number of regions
  $R$ for which $i$ hyperplanes separate $R$ from the region $R_0$ is
  equal to the number of trees on the vertices $0,\ldots,n$ with
  $\binom{n}{2}-i$ inversions. The pair $(i,j)$, where $1\leq i<j$,
  is an inversion for $T$ if the vertex $j$ lies on the unique path in
  $T$ from $0$ to $i$. See \cite[Theorem 5.1]{stanley1996}.

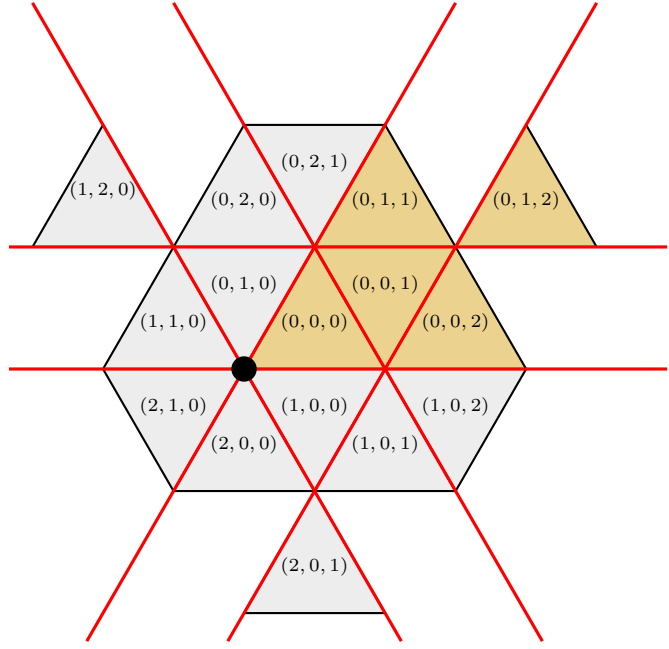
\begin{figure}[ht]

\begin{tikzpicture}[fill opacity=.5, scale=2.5, font=\scriptsize]

\filldraw [draw=black,thick,fill opacity = .5, fill = Gainsboro,shift={(0.000000,-1.29904)}](0,0) -- (60:0.75cm) -- (0:0.75cm) --(0,0) node[blue] at (0.3cm,0.15cm) {} node[black,fill opacity=1] at (0.375cm,0.25cm) {$(2,0,1)$} node[blue] at (0.3cm,0.45cm) {};
[-1,-1,-1]
\filldraw [draw=black, thick,fill = Gainsboro, shift={(-0.375,-0.649519)}](0,0) -- (60:0.75cm) -- (120:0.75cm) --(0,0) node[green,fill opacity=1] at (0cm,0.6cm) {} node[black,fill opacity=1, draw opacity = 1] at (0cm,0.45cm) {$(2,1,0)$} node[green,fill opacity=1] at (0cm,0.3cm) {};
[-1,-1,0]
\filldraw [draw=black,thick,fill opacity = .5, fill = Gainsboro,shift={(-0.375000,-0.649519)}](0,0) -- (60:0.75cm) -- (0:0.75cm) --(0,0) node[blue] at (0.3cm,0.15cm) {} node[black,fill opacity=1] at (0.375cm,0.25cm) {$(2,0,0)$} node[blue] at (0.3cm,0.45cm) {};
[0,-1,0]
\filldraw [draw=black, thick,fill = Gainsboro, shift={(0.375,-0.649519)}](0,0) -- (60:0.75cm) -- (120:0.75cm) --(0,0) node[green,fill opacity=1] at (0cm,0.6cm) {} node[black,fill opacity=1, draw opacity = 1] at (0cm,0.45cm) {$(1,0,0)$} node[green,fill opacity=1] at (0cm,0.3cm) {};
[0,-1,1]
\filldraw [draw=black,thick,fill opacity = .5, fill = Gainsboro,shift={(0.375000,-0.649519)}](0,0) -- (60:0.75cm) -- (0:0.75cm) --(0,0) node[blue] at (0.3cm,0.15cm) {} node[black,fill opacity=1] at (0.375cm,0.25cm) {$(1,0,1)$} node[blue] at (0.3cm,0.45cm) {};
[1,-1,1]
\filldraw [draw=black, thick,fill = Gainsboro, shift={(1.125,-0.649519)}](0,0) -- (60:0.75cm) -- (120:0.75cm) --(0,0) node[green,fill opacity=1] at (0cm,0.6cm) {} node[black,fill opacity=1, draw opacity = 1] at (0cm,0.45cm) {$(1,0,2)$} node[green,fill opacity=1] at (0cm,0.3cm) {};
[-1,0,-1]
\filldraw [draw=black,thick,fill opacity = .5, fill = Gainsboro,shift={(-0.750000,0)}](0,0) -- (60:0.75cm) -- (0:0.75cm) --(0,0) node[blue] at (0.3cm,0.15cm) {} node[black,fill opacity=1] at (0.375cm,0.25cm) {$(1,1,0)$} node[blue] at (0.3cm,0.45cm) {};
[0,0,-1]
\filldraw [draw=black, thick,fill = Gainsboro, shift={(0,0)}](0,0) -- (60:0.75cm) -- (120:0.75cm) --(0,0) node[green,fill opacity=1] at (0cm,0.6cm) {} node[black,fill opacity=1, draw opacity = 1] at (0cm,0.45cm) {$(0,1,0)$} node[green,fill opacity=1] at (0cm,0.3cm) {};
[0,0,0]
\filldraw [draw=black,thick,fill opacity = .5, fill = Goldenrod,shift={(0.000000,0)}](0,0) -- (60:0.75cm) -- (0:0.75cm) --(0,0) node[blue] at (0.3cm,0.15cm) {} node[black,fill opacity=1] at (0.375cm,0.25cm) {$(0,0,0)$} node[blue] at (0.3cm,0.45cm) {};
[1,0,0]
\filldraw [draw=black, thick,fill = Goldenrod, shift={(0.75,0)}](0,0) -- (60:0.75cm) -- (120:0.75cm) --(0,0) node[green,fill opacity=1] at (0cm,0.6cm) {} node[black,fill opacity=1, draw opacity = 1] at (0cm,0.45cm) {$(0,0,1)$} node[green,fill opacity=1] at (0cm,0.3cm) {};
[1,0,1]
\filldraw [draw=black,thick,fill opacity = .5, fill = Goldenrod,shift={(0.750000,0)}](0,0) -- (60:0.75cm) -- (0:0.75cm) --(0,0) node[blue] at (0.3cm,0.15cm) {} node[black,fill opacity=1] at (0.375cm,0.25cm) {$(0,0,2)$} node[blue] at (0.3cm,0.45cm) {};
[-1,1,-2]
\filldraw [draw=black,thick,fill opacity = .5, fill = Gainsboro,shift={(-1.125000,0.649519)}](0,0) -- (60:0.75cm) -- (0:0.75cm) --(0,0) node[blue] at (0.3cm,0.15cm) {} node[black,fill opacity=1] at (0.375cm,0.3cm) {$(1,2,0)$} node[blue] at (0.3cm,0.45cm) {};
[0,1,-1]
\filldraw [draw=black,thick,fill opacity = .5, fill = Gainsboro,shift={(-0.375000,0.649519)}](0,0) -- (60:0.75cm) -- (0:0.75cm) --(0,0) node[blue] at (0.3cm,0.15cm) {} node[black,fill opacity=1] at (0.375cm,0.25cm) {$(0,2,0)$} node[blue] at (0.3cm,0.45cm) {};
[1,1,-1]
\filldraw [draw=black, thick,fill = Gainsboro, shift={(0.375,0.649519)}](0,0) -- (60:0.75cm) -- (120:0.75cm) --(0,0) node[green,fill opacity=1] at (0cm,0.6cm) {} node[black,fill opacity=1, draw opacity = 1] at (0cm,0.45cm) {$(0,2,1)$} node[green,fill opacity=1] at (0cm,0.3cm) {};
[1,1,0]
\filldraw [draw=black,thick,fill opacity = .5, fill = Goldenrod,shift={(0.375000,0.649519)}](0,0) -- (60:0.75cm) -- (0:0.75cm) --(0,0) node[blue] at (0.3cm,0.15cm) {} node[black,fill opacity=1] at (0.375cm,0.25cm) {$(0,1,1)$} node[blue] at (0.3cm,0.45cm) {};
[2,1,1]
\filldraw [draw=black,thick,fill opacity = .5, fill = Goldenrod,shift={(1.125000,0.649519)}](0,0) -- (60:0.75cm) -- (0:0.75cm) --(0,0) node[blue] at (0.3cm,0.15cm) {} node[black,fill opacity=1] at (0.375cm,0.25cm) {$(0,1,2)$} node[blue] at (0.3cm,0.45cm) {};

\def\b{1}
\def\b{1}
\def\ma{1.94856/1.125}
\def\xa{-2.25+\b+1.94856}
\def\ya{\xa*\ma}
\def\ya{1.2099}
\def\dy{.5}
\def\dx{dy/\ma}
\def\dx{.2886747137}

\path [draw = red, very thick,  shift = {(0,0)}](-1.125+\dx,-1.94856+\dy) -- (1.125,1.94856);
\path [draw = red,very thick, shift = {(0,0)}](-2.25+\b,0) -- (2.25,0);
\path [draw = red, very thick,  shift = {(0,0)}](-1.125,1.94856) -- (1.125-\dx,-1.94856+\dy);
\path [draw = red, very thick,  shift = {(0.75,0)}](-1.125+\dx,-1.94856+\dy) -- (1.125,1.94856);
\path [draw = red,very thick, shift = {(0,0.649519)}](-2.25+\b,0) -- (2.25,0);
\path [draw = red, very thick,  shift = {(0.75,0)}](-1.125,1.94856) -- (1.125-\dx,-1.94856+\dy);
\fill[fill=black, fill opacity=1.](0.,0.) circle (0.0675cm);[-2,-4,2]

\end{tikzpicture}
  \caption{\small  The Shi regions labeled by parking functions, using the Pak-Stanley labeling. A label is nondecreasing if and only if the region is in the dominant region. }
  \label{fig:stanPak}
\end{figure}

We mention a few papers which build on the Pak-Stanley bijection. Duarte and Guedes de Oliveira, for example, further analyzed this bijection in \cite{duarte-guedes2015}. \DGdO{In \cite{duarte-guedes2019arxiv} they labelled the regions of the related $m$-Catalan arrangement.} Rinc\'on \cite{rincon2007} extended the Pak-Stanley labeling to the poset of faces of the Shi arrangement. See also Section~\ref{sec:further}. 

\subsection{More}

Believe it or not, there are still other wonderful proofs concerning the
number of regions.

For example, Athanasiadis and Linusson
\cite{athanasiadis-linusson1999} defined a bijection, different from
Pak and Stanley's, from parking functions to the Shi regions (type
$A$). Theirs gives a simple proof of the number of regions.
Their bijection was generalized to type
$C$ by M\'esz\'aros in \cite{meszaros2013}.  In \cite[Section
  5.2]{armstrong-rhoades2012} there is another proof of the formula
for the number of regions, using Armstrong's and Rhoades' ceiling
diagrams, which we define in Section~\ref{sec:ish}. The ceiling
diagrams are related to the diagrams Athanasiadis and Linusson used.
Armstrong, Reiner, and Rhoades \cite{armstrong-reiner-rhoades2015}
define {\em nonnesting parking functions} using the root poset and
permutations from the finite Weyl group and label the Shi regions with
these.  Their definition can also be used for types which are not
crystallographic.\omitt{Fomin and Mikhalkin mention that their floor
  diagrams of genus 0 and degree $d$ are in bijection with the
  regions\cite{fomin-mikhalkin2010}.} Other, more recent and more
general bijections include \cite{hopkins-perkinson2016,
  beck-berrizbeitia-dairyko-rodriguez-ruiz-veeneman2015} for example.

\subsection{The Ish and the Shi}\label{sec:ish}
We'll start off by writing the $q,t$-Catalan polynomial
combinatorially:
\begin{equation}\label{eq:cat} C_n(q,t)=\sum_{\pi}q^{\area(\pi)}t^{\bounce(\pi)},\end{equation}

where the sum is over all Dyck paths of length $n$. The
$(q,t)$-Catalan polynomials are remarkable generating functions coming
from representation theory. They have been intensely studied since
their introduction by Garsia and Haiman in
\cite{garsia-haiman1996}. See Haglund's monograph \cite[Chapter
  3]{haglund2008} for more information.  There are the same number of
dominant Shi regions as there are Dyck paths, and one of Armstrong's
results in \cite{armstrong2013} (and the one we'll describe) was to
transfer the statistics $\area$ and $\bounce$ to dominant regions. His
statistics are actually for all regions. The statistic $\shi$ will
correspond to the statistic $\area$ and $\shi(R)$ is defined as the
number of hyperplanes which must be crossed on a trip to the region
$R$ from $R_0=\Afund$. The statistic $\ish$ is defined using a second
hyperplane arrangement, the Ish arrangement $\Ish{n}$. It is defined
for type $A$ and is a deformation of the Coxeter arrangement. Let
$\roots$ be the set of roots for type $A$, so $\alpha_i = \e_i -
\e_{i+1}$ as in Section~\ref{subsec:affSn}. Denote
$\alpha_i+\alpha_{i+1}+\cdots+\alpha_{n-1}$ by
$\tilde{\alpha}_i$. Then the definition of the Ish arrangement is

$$\Ish{n}=\Cox{n}\cup\{\Hak{\tilde{\alpha}_j}{k}| 1\leq j\leq n-1, k\in\{1,2,\ldots, n-j\}\}=\Cox{n}\cup\{x_j-x_n=k|k\in\{1,\ldots,n-j\},1\leq j\leq n-1\}.$$

The $\ish$ statistic is defined on Shi regions using the hyperplanes
in $\Ish{n}$.  \omitt{We now let the affine symmetric group act on the left, following \cite{armstrong13}.} Let $R$ be region of the arrangement $\Shi{n}$ with
minimal alcove $\A_R$. There is a unique $w\in\affSn$ such that
$\A_R=\Afund w.$ The affine permutation $w$ has a unique factorization
$w^I\cdot w_I$ \cite{bjorner-brenti2005}, where $w_I\in\Sn$ and $w^I$
is a {\em minimal length coset representative}, which we won't
define. What is important for us is that $\Afund w^I$ is an alcove in
the dominant chamber since $w^I$ is a minimal length coset
representative. Then $\ish(R)=\ish(\A_R)$ is the number of hyperplanes
in $\Ish{n}$ which must be crossed in traveling from $\Afund$ to $\A w^I$.

\begin{figure}[ht]

\begin{tikzpicture}[fill opacity=.5, scale=2., font=\scriptsize,shift={(4.2,0)}]

\def\b{1}
\def\ma{1.94856/1.125}
\def\xa{-2.25+\b+1.94856}
\def\ya{\xa*\ma}
\def\ya{1.2099}
\def\dy{.5}
\def\dx{dy/\ma}
\def\dx{.2886747137}

\path [draw = red, very thick,  shift = {(0,0)}](-1.125+\dx,-1.94856+\dy) -- (1.125,1.94856);
\path [draw = red,very thick, shift = {(0,0)}](-2.25+\b,0) -- (2.25,0);
\path [draw = red, very thick,  shift = {(0,0)}](-1.125,1.94856) -- (1.125-\dx,-1.94856+\dy);
\path [draw = red, very thick,  shift = {(0.75,0)}](-1.125+\dx,-1.94856+\dy) -- (1.125,1.94856);
\path [draw = red, very thick,  shift = {(0.75,0)}](-1.125,1.94856) -- (1.125-\dx,-1.94856+\dy);
\path [draw = red, very thick,  shift = {(1.5,0)}](-1.125,1.94856) -- (1.125-\dx,-1.94856+\dy);

\fill[fill=black, fill opacity=1.](0.,0.) circle (0.0675cm);[-2,-4,2]

  \begin{scope}[shift={(4.2,0)}]
\filldraw [draw=black,thick,fill opacity = .5, fill = Gainsboro,shift={(0.000000,-1.29904)}](0,0) -- (60:0.75cm) -- (0:0.75cm) --(0,0) node[blue] at (0.3cm,0.15cm) {} node[black,fill opacity=1] at (0.375cm,0.25cm) {$(3,1)$} node[blue] at (0.3cm,0.45cm) {};
[-1,-1,-1]
\filldraw [draw=black, thick,fill = Gainsboro, shift={(-0.375,-0.649519)}](0,0) -- (60:0.75cm) -- (120:0.75cm) --(0,0) node[green,fill opacity=1] at (0cm,0.6cm) {} node[black,fill opacity=1, draw opacity = 1] at (0cm,0.45cm) {$(3,0)$} node[green,fill opacity=1] at (0cm,0.3cm) {};
[-1,-1,0]
\filldraw [draw=black,thick,fill opacity = .5, fill = Gainsboro,shift={(-0.375000,-0.649519)}](0,0) -- (60:0.75cm) -- (0:0.75cm) --(0,0) node[blue] at (0.3cm,0.15cm) {} node[black,fill opacity=1] at (0.375cm,0.25cm) {$(2,0)$} node[blue] at (0.3cm,0.45cm) {};
[0,-1,0]
\filldraw [draw=black, thick,fill = Gainsboro, shift={(0.375,-0.649519)}](0,0) -- (60:0.75cm) -- (120:0.75cm) --(0,0) node[green,fill opacity=1] at (0cm,0.6cm) {} node[black,fill opacity=1, draw opacity = 1] at (0cm,0.45cm) {$(1,0)$} node[green,fill opacity=1] at (0cm,0.3cm) {};
[0,-1,1]
\filldraw [draw=black,thick,fill opacity = .5, fill = Gainsboro,shift={(0.375000,-0.649519)}](0,0) -- (60:0.75cm) -- (0:0.75cm) --(0,0) node[blue] at (0.3cm,0.15cm) {} node[black,fill opacity=1] at (0.375cm,0.25cm) {$(2,1)$} node[blue] at (0.3cm,0.45cm) {};
[1,-1,1]
\filldraw [draw=black, thick,fill = Gainsboro, shift={(1.125,-0.649519)}](0,0) -- (60:0.75cm) -- (120:0.75cm) --(0,0) node[green,fill opacity=1] at (0cm,0.6cm) {} node[black,fill opacity=1, draw opacity = 1] at (0cm,0.45cm) {$(3,2)$} node[green,fill opacity=1] at (0cm,0.3cm) {};
[-1,0,-1]
\filldraw [draw=black,thick,fill opacity = .5, fill = Gainsboro,shift={(-0.750000,0)}](0,0) -- (60:0.75cm) -- (0:0.75cm) --(0,0) node[blue] at (0.3cm,0.15cm) {} node[black,fill opacity=1] at (0.375cm,0.25cm) {$(2,0)$} node[blue] at (0.3cm,0.45cm) {};
[0,0,-1]
\filldraw [draw=black, thick,fill = Gainsboro, shift={(0,0)}](0,0) -- (60:0.75cm) -- (120:0.75cm) --(0,0) node[green,fill opacity=1] at (0cm,0.6cm) {} node[black,fill opacity=1, draw opacity = 1] at (0cm,0.45cm) {$(1,0)$} node[green,fill opacity=1] at (0cm,0.3cm) {};
[0,0,0]
\filldraw [draw=black,thick,fill opacity = .5, fill = Goldenrod,shift={(0.000000,0)}](0,0) -- (60:0.75cm) -- (0:0.75cm) --(0,0) node[blue] at (0.35cm,0.15cm) {$q^3$} node[black,fill opacity=1] at (0.375cm,0.35cm) {$(0,0)$} node[blue] at (0.3cm,0.5cm) {};
[1,0,0]
\filldraw [draw=black, thick,fill = Goldenrod, shift={(0.75,0)}](0,0) -- (60:0.75cm) -- (120:0.75cm) --(0,0) node[blue] at (0cm,0.25cm) {$q^2t$} node[black,fill opacity=1, draw opacity = 1] at (0cm,0.45cm) {$(1,1)$} node[green,fill opacity=1] at (0cm,0.3cm) {};
          [1,0,1]
\filldraw [draw=black,thick,fill opacity = .5, fill = Goldenrod,shift={(0.750000,0)}](0,0) -- (60:0.75cm) -- (0:0.75cm) --(0,0) node[blue] at (0.35cm,0.15cm) {$qt^2$} node[black,fill opacity=1] at (0.375cm,0.35cm) {$(2,2)$} node[blue] at (0.3cm,0.45cm) {};
[-1,1,-2]
\filldraw [draw=black,thick,fill opacity = .5, fill = Gainsboro,shift={(-1.125000,0.649519)}](0,0) -- (60:0.75cm) -- (0:0.75cm) --(0,0) node[blue] at (0.35cm,0.15cm) {} node[black,fill opacity=1] at (0.375cm,0.3cm) {$(3,2)$} node[blue] at (0.3cm,0.45cm) {};
[0,1,-1]
\filldraw [draw=black,thick,fill opacity = .5, fill = Gainsboro,shift={(-0.375000,0.649519)}](0,0) -- (60:0.75cm) -- (0:0.75cm) --(0,0) node[blue] at (0.3cm,0.15cm) {} node[black,fill opacity=1] at (0.375cm,0.35cm) {$(2,1)$} node[blue] at (0.3cm,0.45cm) {};
[1,1,-1]
\filldraw [draw=black, thick,fill = Gainsboro, shift={(0.375,0.649519)}](0,0) -- (60:0.75cm) -- (120:0.75cm) --(0,0) node[green,fill opacity=1] at (0cm,0.6cm) {} node[black,fill opacity=1, draw opacity = 1] at (0cm,0.45cm) {$(3,1)$} node[green,fill opacity=1] at (0cm,0.3cm) {};
          [1,1,0]
\filldraw [draw=black,thick,fill opacity = .5, fill = Goldenrod,shift={(0.375000,0.649519)}](0,0) -- (60:0.75cm) -- (0:0.75cm) --(0,0) node[blue] at (0.35cm,0.15cm) {$qt$} node[black,fill opacity=1] at (0.375cm,0.35cm) {$(2,1)$} node[blue] at (0.3cm,0.45cm) {};
[2,1,1]
\filldraw [draw=black,thick,fill opacity = .5, fill = Goldenrod,shift={(1.125000,0.649519)}](0,0) -- (60:0.75cm) -- (0:0.75cm) --(0,0) node[blue] at (0.35cm,0.15cm) {$t^3$} node[black,fill opacity=1] at (0.375cm,0.25cm) {$(3,3)$} node[blue] at (0.3cm,0.45cm) {};

\def\b{1}
\def\b{1}
\def\ma{1.94856/1.125}
\def\xa{-2.25+\b+1.94856}
\def\ya{\xa*\ma}
\def\ya{1.2099}
\def\dy{.5}
\def\dx{dy/\ma}
\def\dx{.2886747137}

\path [draw = red, very thick,  shift = {(0,0)}](-1.125+\dx,-1.94856+\dy) -- (1.125,1.94856);
\path [draw = red,very thick, shift = {(0,0)}](-2.25+\b,0) -- (2.25,0);
\path [draw = red, very thick,  shift = {(0,0)}](-1.125,1.94856) -- (1.125-\dx,-1.94856+\dy);

\path [draw = red, very thick,  shift = {(0.75,0)}](-1.125+\dx,-1.94856+\dy) -- (1.125,1.94856);

\path [draw = red,very thick, shift = {(0,0.649519)}](-2.25+\b,0) -- (2.25,0);
\path [draw = red, very thick,  shift = {(0.75,0)}](-1.125,1.94856) -- (1.125-\dx,-1.94856+\dy);


\fill[fill=black, fill opacity=1.](0.,0.) circle (0.0675cm);[-2,-4,2]

\end{scope}

\end{tikzpicture}
\caption{\small  The Ish arrangement for $n=3$ is on the left. The Shi arrangement is on the right, where region $R$ is labeled by the pair $(\shi(R),\ish(R))$. The dominant regions are also labeled by $$q^{\binom{3}{2}-\shi(R)}t^{\ish(R)}.$$ The sum of the monomials is $C_3(q,t)=q^3+q^2t+qt+qt^2+t^3$.}

  \label{fig:ishShi}
\end{figure}


Each dominant Shi region $R$ corresponds to a Dyck path
$\pi_R$. Armstrong showed that $\binom{n}{2}-\shi(R)=\area(\pi_R)$ and
$\ish(R)=\bounce(\pi_R)$ Notice that the $\ish$ and $\shi$ statistics
are defined on all regions, not just the dominant ones. Armstrong was
able to show they agree with $\bounce$ and $\area$ on all {\em
  diagonally labeled} Dyck paths. See \cite[Chapter 5]{haglund2008}
and \cite[Section 3]{armstrong2013}.


Armstrong and Rhoades concentrated on properties of the Ish
arrangement, especially its uncanny similarities to the Shi
arrangement, in \cite{armstrong-rhoades2012}. Their definition of the arrangement changes just a bit: replace $x_j-x_n=i$ by $x_1-x_{n-j+1}=n-i+1$. 
Their main
theorem is for {\em deleted} versions (more general) of the arrangements (see
Section~\ref{sec:themes}), but we'll stick with the full
arrangements. That is, $G$ is the complete graph in this survey. We need to define a
few terms before we can state the main theorem. The wall $\Hy$ of
a region $R$ is called a {\em ceiling} if it does not contain the
origin and if the origin and $R$ are not separated by $\Hy$ (they
lie in the same half-space of $\Hy$). The regions of both
$\Ish{n}$ and $\Shi{n}$ are convex, so every region has a {\em recession
cone}:
$$\rec{R}=\{v\in V:v+R\subseteq R\}.$$ The cone is closed under
nonnegative linear combinations and has a dimension. The dimension of
$\rec{R}$ is called the {\em degrees of freedom} of $R$. It's worth
mentioning that the region $R$ is bounded if and only if
$\rec{R}=\{0\}$. 

A simplified version of their main theorem can now be stated: 

\begin{theorem}[\cite{armstrong-rhoades2012}] 
\label{thm:ARMain}
Let $c$ and $d$ be nonnegative integers. The $\Ish{n}$ and $\Shi{n}$ have the same 
\begin{enumerate}
\item \label{ishShiChar} characteristic polynomial,
\item \label{ishShicd} number of dominant regions with $c$ ceilings, and
\item \label{ishShicd} number of regions with $c$ ceilings and $d$ degrees of freedom. 
\end{enumerate}
\end{theorem}

We are not presenting their theorem in its full generality, and as
written here, \eqref{ishShiChar} was proved in \cite{armstrong2013}.

We want to define the ceiling diagrams for the Shi arrangement because
they show the properties of the corresponding region so clearly.  We
will also need variations on the root poset, which they use to prove
Theorem~\ref{thm:ARMain}. First the definition of the Shi ceiling
diagram of a region $R$. Suppose the region is in the chamber
$\domCham w$, where $\domCham$ is the dominant chamber and
$w\in\Sn$. Then we define the set partition $\sigma_R$: there is an arc
from $i$ to $j$, $i<j$, in the diagram of $\sigma_R$ if and only if the
hyperplane $x_{w(i)}-x_{w(j)}=1$ is a ceiling of $R$. We draw the {\em
  Shi ceiling diagram $(w,\sigma_R)$} by placing the arc diagram for
$\sigma_R$ above $w(1),w(2),\ldots,w(n)$. See
Example~\ref{ex:ceilingDiagrams}. Armstrong and Rhoades show that
$\sigma_R$ is a nonnesting set partition. The number of arcs is $c$. What
about $d$?  Let $d'$ be the number of $k$, $1\leq k\leq n-1$ where
there is no arc covering the space between $k$ and $k+1$; that is, the
number of $k$ where for which there is no $i<j$ such that $i\leq k<j$
and there is an arc from $i$ to $j$. For example, $d'=0$ for the set
partition in Figure~\ref{fig:catObjs} and $d'=1$ for the set partition
$12|35|4$. Then set $d=d'+1$. Additionally, the recession cone
$\rec{R}$ can be read from the diagram.

There is still a key point: for a fixed $w\in\Sn$, both the regions
and the ceiling diagrams are in bijection with antichains in
$\pos(w)$. The poset $\pos(w)$ is the first variation on the root
poset: $$\pos(w)=\{x_{w(i)}-x_{w(j)}=1:w(i)<w(j)\}.$$ The elements of
$\pos(w)$ are the affine hyperplanes in the Shi arrangement which
intersect $\domCham w$. The partial order on the hyperplanes is given
by $$x_{w(i')}-x_{w(j')}=1\leq x_{w(i)}-x_{w(j)}=1$$ if $w(i)\leq
w(i')<w(j')\leq w(j)$. The partial order is defined so that ceilings
of any Shi region $R$, $R\in \domCham w$, are the maximal elements of
an order ideal. The number $c$ shows up as the number of these maximal
elements.\omitt{see Thm 4.5 of AR, and right after}

We'll need the second variation on $\pos$, $\Psi^+(w)$, for discussing
the Ish arrangement: $$\Psi^+(w)=\{x_1-x_j=i:w^{-1}(i)<w^{-1}(j)\}.$$
Its elements are the Ish hyperplanes that intersect $w\domCham$ and its
partial order is chosen so that the ceilings of any Ish region $R$,
$R\in w\domCham$, are minimal elements of a filter, and $c$ for the
region is the number of these minimal elements.

\begin{figure}[h]
\begin{tikzpicture}
\begin{scope}[scale=.4]

\def\a{2.59807}

\omitt{
\filldraw[fill=gray!10,draw opacity=0,fill opacity=.8](6-1.5,-3*\a)--(3+3*1.5,-3*\a)--(3,0)--(0,0)--(-6+1.5,-3*\a)--(3-3*1.5,-3*\a)--(1.5,-\a)--cycle;

\filldraw[fill=red!10,draw opacity=0,fill opacity=.8](3+3*1.5,-3*\a)--(3,0)--(15-3,0)--cycle;

\filldraw[fill=red!10,draw opacity=0,fill opacity=.8](6-1.5,-3*\a)--(3-3*1.5,-3*\a)--(1.5,-\a)--cycle;

}

\filldraw[fill=gray!10, draw opacity=0,fill opacity=.8](1.5,\a)--(6,4*\a)--(7.5,3*\a)--(4.5,\a)--cycle;

\omitt{

\filldraw[fill=red!10, draw opacity=0,fill opacity=.8](1.5,\a)--(6,4*\a)--(-3,4*\a)--cycle;

\filldraw[fill=red!10, draw opacity=0,fill opacity=.8](7.5,3*\a)--(4.5,\a)--(12-1.5,\a)--cycle;

}


\filldraw[fill=gray!10, draw opacity=0,fill opacity=.8](1.5,\a)--(-3,4*\a)--(-4.5,3*\a)--(-1.5,\a)--cycle;

\omitt{
\filldraw[fill=red!10, draw opacity=0,fill opacity=.8](0,0)--(-12+3,0)--(-6+1.5,-3*\a)--cycle;

\filldraw[fill=red!10, draw opacity=0,fill opacity=.8](-4.5,3*\a)--(-1.5,\a)--(-12+4.5,\a)--cycle;

} 


\draw [draw = blue,  thick, draw opacity = 1,shift = {(3,0)}] (240:12-3) -- (60:15-6);


\draw [draw = blue, thick, draw opacity = 1, shift = {(0,\a)}] (-12+4.5,0) -- (15-4.5,0);


\draw [draw = blue, thick,draw opacity = 1,shift={(3,0)}] (120:15-3) -- (300:12-3);

\draw [draw = black, ultra thick, draw opacity = 1] (240:12-3) -- (60:15-3);

\draw [draw = black,ultra thick,draw opacity = 1] (120:12-3) -- (300:12-3);

\draw [draw = black,ultra thick, draw opacity = 1] (-12+3,0) -- (15-3,0);


\def\r{7}
\omitt{
\begin{scope}[shift={(1,.33*\a)},scale=.5]
\Triangle{$\bigcirc$}{$\bigcirc$}{$\bigcirc$};
\end{scope}

\begin{scope}[shift={(2.5,.67*\a)},scale=.5]
\Triangle{$+$}{$\bigcirc$}{$\bigcirc$};
\end{scope}

\begin{scope}[shift={(-.6,.67*\a)},scale=.5]
\Triangle{$\bigcirc$}{$\bigcirc$}{$-$};
\end{scope}

\begin{scope}[shift={(1,-.4*\a)},scale=.5]
\Triangle{$\bigcirc$}{$-$}{$\bigcirc$};
\end{scope}

\begin{scope}[shift={(169:\r)},scale=.5]
\Triangle{$-$}{$\bigcirc$}{$-$};
\end{scope}
} 

\begin{scope}[shift={(140:\r)},scale=.5,font=\LARGE]
\node[ellipse,fill=yellow!20,font=\LARGE] at (0,0) {$s_1s_2$};
\end{scope}

\begin{scope}[shift={(110:\r)},scale=.5]
\node at (0,0) {$R_2$};
\end{scope}

\begin{scope}[thick, shift={(110:\r-2.5)}]
\tikzstyle{every node}=[font=\scriptsize]
\def\a{3}
\def\b{-2}
\def\c{.4}
\begin{scope}[shift={(0*\a,0*\b)},scale=.8]
  \node[below,blue] (1) at (0,0) {};
  \node[below,blue] at (0,-\c) {$1$};
\fill[red] ($(1)$) circle (2pt);

\node[below,blue] (2) at (1,0) {};
\node[below,blue]  at (1,-\c) {$3$};
\fill[red] ($(2)$) circle (2pt);

\node[below,blue] (3) at (2,0) {};
\node[below,blue]  at (2,-\c) {$2$};
\fill[red] ($(3)$) circle (2pt);

\omitt{
\node[below,blue] (4) at (3,0) {};
\node[below,blue]  at (3,-\c) {4};
\fill[red] ($(4)$) circle (2pt);
\node[below,blue] (5) at (4,0) {};
\node[below,blue]  at (4,-\c) {5};
\fill[red] ($(5)$) circle (2pt);
}

\draw (1) .. controls (.5,.5) .. (2);
\omitt{
\draw (4) .. controls (3.5,.5) .. (5);
\draw (2) .. controls (2,.5) .. (4);
}

\end{scope}
\end{scope}

\begin{scope}[shift={(77:\r)},scale=.5,font=\LARGE]
\node[ellipse,fill=yellow!20] at (0,0) {$s_2$};
\end{scope}

\begin{scope}[shift={(50:\r)},scale=.5]
\node at (0,0) {$R_1$};
\end{scope}

\begin{scope}[thick, shift={(55:\r-2.0)}]
\tikzstyle{every node}=[font=\scriptsize]
\def\a{3}
\def\b{-2}
\def\c{.4}
\begin{scope}[shift={(0*\a,0*\b)},scale=.8]
  \node[below,blue] (1) at (0,0) {};
  \node[below,blue] at (0,-\c) {$1$};
\fill[red] ($(1)$) circle (2pt);

\node[below,blue] (2) at (1,0) {};
\node[below,blue]  at (1,-\c) {$2$};
\fill[red] ($(2)$) circle (2pt);

\node[below,blue] (3) at (2,0) {};
\node[below,blue]  at (2,-\c) {$3$};
\fill[red] ($(3)$) circle (2pt);

\draw (2) .. controls (1.5,.5) .. (3);
\omitt{
\draw (4) .. controls (3.5,.5) .. (5);
\draw (2) .. controls (2,.5) .. (4);
}

\end{scope}
\end{scope}

\begin{scope}[shift={(30:\r)},scale=.5,font=\LARGE]
\node[ellipse,fill=yellow!20] at (0,0) {$e$};
\end{scope}

\omitt{
\begin{scope}[shift={(10:\r)},scale=.5]
\Triangle{$+$}{$\bigcirc$}{$+$};
\end{scope}
} 

\begin{scope}[shift={(340:\r)},scale=.5,font=\LARGE]
\node[ellipse,fill=yellow!20] at (0,0) {$s_1$};
\end{scope}

\omitt{
\begin{scope}[shift={(310:\r)},scale=.5]
\Triangle{$\bigcirc$}{$-$}{$+$};
\end{scope}
} 

\begin{scope}[shift={(280:\r)},scale=.5,font=\LARGE]
\node[ellipse,fill=yellow!20] at (0,0) {$s_2s_1$};
\end{scope}

\omitt{
\begin{scope}[shift={(250:\r)},scale=.5]
\Triangle{$-$}{$-$}{$\bigcirc$};
\end{scope}
}

\begin{scope}[shift={(210:\r)},scale=2,font=\LARGE]
\node[ellipse,fill=yellow!20] at (0,0) {$s_1s_2s_1$};
\end{scope}

\node at (0,0) [circle, fill = red]{};

\end{scope}

\end{tikzpicture}
\caption{The Shi arrangement for $A_2$. The chamber $(\domCham)w$ is labeled by $w$.}
\label{fig:ceilingDiagrams}
\end{figure}
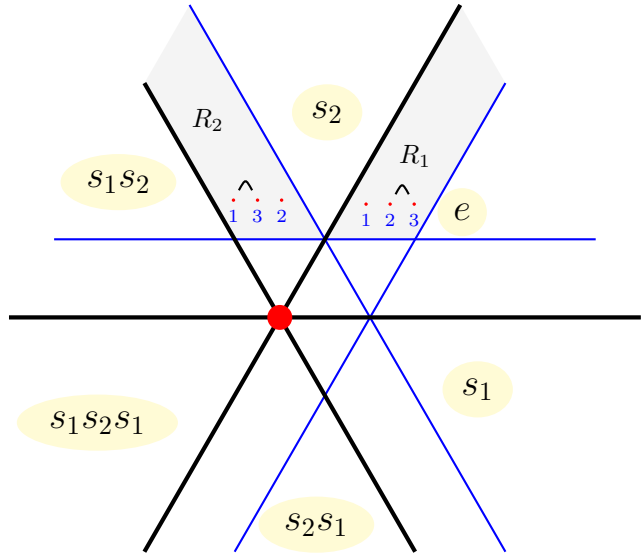

\begin{figure}[h]
\begin{tikzpicture}

\def\a{2}
\def\b{2}
\tikzstyle{every node}=[font=\footnotesize]
\node[ellipse,fill=green!20,font=\small] (13) at (0,0) {$x_1-x_3=1$};
\node[ellipse,fill=green!20,font=\small] (12) at (-\a,-\b) {$x_1-x_2=1$};
\node[ellipse,fill=green!20,font=\small] (23) at (\a,-\b) {$x_2-x_3=1$};
\draw[thick](13)--(12);
\draw[thick](13)--(23);

\begin{scope}[shift={(3*\b,0)}]
\node[ellipse,fill=green!20,font=\small] (13) at (0,0) {$x_1-x_3=1$};
\node[ellipse,fill=green!20,font=\small] (12) at (0,-\b) {$x_1-x_2=1$};
\draw[thick](13)--(12);

\end{scope}

\end{tikzpicture}
\caption{On the left, the poset $\pos(e)$.  The poset $\pos(s_2)$ is on the right.}
\label{fig:moreRootPosets}
\end{figure}

\begin{example}
This example refers to Figures~\ref{fig:ceilingDiagrams} and
\ref{fig:moreRootPosets}. In \ref{fig:ceilingDiagrams}, the chamber
$\domCham w$ is labeled by $w\in\Sym{3}.$ We've picked two Shi regions
to consider in this example-the ones we have labeled $R_1$ and $R_2$.

$R_1$ is in the chamber $\domCham w$ for $w=e$, the identity. There
are three hyperplanes of the form $x_i-x_j=1$ which intersect
$\domCham w$ and the poset $\pos(w)$ is on the left in
Figure~\ref{fig:moreRootPosets}. The hyperplane $x_2-x_3=1$ is a
ceiling for $R_1$ and the region is also labeled with its ceiling
diagram in Figure~\ref{fig:ceilingDiagrams}.

The region $R_2$ is in the chamber $\domCham w$ for $w=s_2$ and poset
$\pos(s_2)$ is on the right in Figure~\ref{fig:moreRootPosets}. The
poset $\pos(s_2)$ has three ideals, corresponding to the three Shi
regions in $\domCham s_2$. Our region $R_2$ has ceiling
$x_1-x_3=x_{w(1)}-x_{w(2)}=1$ and we have placed the arc diagram for
the partition $12|3$ above $w(1)w(2)w(3)$ to build the ceiling
diagram.
\label{ex:ceilingDiagrams}
\end{example}

The relationship between the filters and the regions is bijective in
the Ish case, just as between ideals and regions in the Shi case. The
posets $\pos(w)$ and $\Psi^+(w)$ are dual to each other when $w$ is
the identity permutation $e$. The final step in the proof of
Theorem~\ref{thm:ARMain}, part~\eqref{ishShicd}, is simply to send an
order ideal in $\pos(e)$ to the corresponding filter in
$\Psi^+(e)$. Since the maximal elements in the ideal become the
minimal elements in the filter, $c$ is preserved. There are also
ceiling diagrams for the Ish arrangement, but we won't define them.


To prove Theorem~\ref{thm:ARMain}, part~\eqref{ishShicd},
Armstrong and Rhoades used {\em ceiling partitions}, which are set partitions of
$[n]$.  We now define a simplified version of them. First suppose $R$
is a Shi region. The ceiling partition $\pi_R$ has an arc from $i$ to
$j$, $i<j$, if and only if the hyperplane $x_i-x_j=1$ is a ceiling of
$R$.  Next suppose we have an Ish region $R$. Its ceiling partition
has an arc from $i$ to $j$, $i<j$, if and only if $x_1-x_j=i$ is a
ceiling of $R$. The definition of the ceiling partition does not
depend on the chamber of $R$ for either arrangement.  Surprisingly,
the distribution of the set partitions is the same for the Ish and Shi
arrangements.

\begin{theorem}[\cite{armstrong-rhoades2012}]
  \label{thm:AR5.1}
  Let $\mc{A}$ be either the Ish or the Shi arrangement. Let $\pi$ be a partition of $[n]$ with $k$ blocks and let $1\leq d\leq k$.
  \begin{enumerate}
  \item The number of regions of $\mc{A}$ with ceiling partition $\pi$ is
    $$\frac{n!}{(n-k+1)!}.$$
  \item \label{AR5.1part2} The number of regions of $\mc{A}$ with ceiling partition $\pi$ and $d$ degrees of freedom is
    $$\frac{d(n-d-1)!(k-1)!}{(n-k-1)!(k-d)!}.$$
    \end{enumerate}

  \end{theorem}

To obtain the number of regions with $c$ ceilings and $d$ degrees of freedom, thereby proving Theorem~\ref{thm:ARMain}, part \eqref{ishShicd}, sum the expression in Theorem~\ref{thm:AR5.1}, part \eqref{AR5.1part2}, over all partitions $\pi$ with $k=n-c$ blocks.


For space reasons, we cannot include the arguments here for Theorem~\ref{thm:ARMain} and Theorem~\ref{thm:AR5.1}. This is a shame, because we thereby don't present evidence for their observation \cite{armstrong-rhoades2012}:
\begin{quote}
  The Ish arrangement is something of a ``toy model'' for the Shi arrangement (and other Catalan objects). That is, for any property $P$ that $\Shi{n}$ and $\Ish{n}$ share, the proof that $\Ish{n}$ satisfies $P$ is easier than the proof that $\Shi{n}$ satisfies $P$.
  \end{quote}

Many of the theorems of \cite{armstrong-rhoades2012} are proved bijectively by Leven, Rhoades, and Wilson in \cite{leven-rhoades-wilson2014}. 

\DGdO{Duarte and Guedes de Oliveira introduced a new family of hyperplane arrangements in dimension $n\geq3$ that includes both the Shi arrangement and the Ish arrangement in \cite{duarte-guedes2018} and characterized its Pak-Stanley labeling in \cite{duarte-guedes2019}.}

\subsection{Extended Shi arrangement} 
\label{sec:extShi}
In \cite{postnikov-stanley2000}, Postnikov and Stanley introduced the extended Shi arrangement of type $A_{n-1}$:

 $$\mShi{n}{m}=\{ \Hak{\alpha}{k} \mid \alpha \in \pos, -m+1 \leq k \le m \}.$$

This kind of extension is sometimes denoted by Fuss, as in {\em
  Fuss-Catalan} \cite{armstrong2009}.  Up until now, we have been
discussing $m=1$. Postnikov and Stanley show that $\mShi{n}{m}$ has
$(mn+1)^{n-1}$ regions. They fix $m$, set $f_n=\regNum(\mShi{n}{m})$
to be the number of regions, and show that the exponential generating
function $$f=\sum_{n\geq 0}f_n\frac{x^n}{n!}$$
satisfies $$f=e^{xf^m}.$$ The extended Shi arrangement is a special
case ($a=m,b=m+1$) of what they named {\em truncated affine
  arrangements}; see \cite[Section 9]{postnikov-stanley2000} for more
details. The dominant regions of the $m$-Shi are the same as the
dominant regions of the $m$-Catalan.


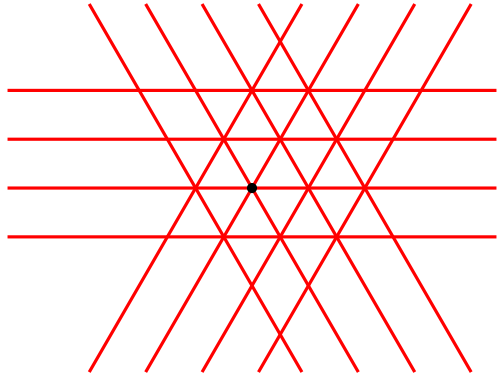
\begin{figure}[ht]

\begin{tikzpicture}[fill opacity=.5, scale=1., font=\scriptsize]

\def\b{0}
\def\ma{0}
\def\xa{0}
\def\ya{0}
\def\ya{0}
\def\dy{0}
\def\dx{0}
\def\dx{0}

\def\b{1}
\def\ma{1.94856/1.125}
\def\xa{-2.25+\b+1.94856}
\def\ya{\xa*\ma}
\def\ya{1.2099}
\def\dy{.5}
\def\dx{dy/\ma}
\def\dx{.2886747137}

\path [draw = red, very thick,  shift = {(2*0.75,0)}](-1.125-\dx,-1.94856-\dy) -- (1.125+\dx,1.94856+\dy);
\path [draw = red,very thick, shift = {(0,2*0.649519)}](-2.25-\b,0) -- (2.25+\b,0);
\path [draw = red, very thick,  shift = {(2*0.75,0)}](-1.125-\dx,1.94856+\dy) -- (1.125+\dx,-1.94856-\dy);

\path [draw = red, very thick,  shift = {(-1*0.75,0)}](-1.125-\dx,-1.94856-\dy) -- (1.125+\dx,1.94856+\dy);
\path [draw = red,very thick, shift = {(0,-1*0.649519)}](-2.25-\b,0) -- (2.25+\b,0);
\path [draw = red, very thick,  shift = {(-1*0.75,0)}](-1.125-\dx,1.94856+\dy) -- (1.125+\dx,-1.94856-\dy);

\path [draw = red, very thick,  shift = {(0,0)}](-1.125-\dx,-1.94856-\dy) -- (1.125+\dx,1.94856+\dy);
\path [draw = red,very thick, shift = {(0,0)}](-2.25-\b,0) -- (2.25+\b,0);
\path [draw = red, very thick,  shift = {(0,0)}](-1.125-\dx,1.94856+\dy) -- (1.125+\dx,-1.94856-\dy);

\path [draw = red, very thick,  shift = {(0.75,0)}](-1.125-\dx,-1.94856-\dy) -- (1.125+\dx,1.94856+\dy);
\path [draw = red,very thick, shift = {(0,0.649519)}](-2.25-\b,0) -- (2.25+\b,0);
\path [draw = red, very thick,  shift = {(0.75,0)}](-1.125-\dx,1.94856+\dy) -- (1.125+\dx,-1.94856-\dy);
\def\g{20}

\fill[fill=black, fill opacity=1.](0.,0.) circle (0.0675cm);[-4,-2,-2]

\end{tikzpicture}
  \caption{\small The $\mShi{3}{2}$ arrangement. There are $49=(2\cdot 3+1)^2$ regions. }
  \label{fig:shi2}
\end{figure}

Many of the other enumerative treats of the Shi arrangement generalize
well to the extended Shi arrangement. For example, Stanley
\cite{stanley1998} labeled the $m$-Shi regions with $m$-parking
functions of length $n$ using an extended version of the bijection
described here in Section~\ref{sec:stanEnum}. Stanley defined {\em
  $m$-parking functions of length $n$}. He replaced the condition that
``$b_i\leq i-1$'' in the definition of parking function (see
Section~\ref{sec:stanEnum}) by ``$b_i\leq m(i-1)$.'' If we set
$d(R)=\area(R)$ (see Section~\ref{sec:ish}) , then Stanley's bijection showed
\cite[Corollary 2.2]{stanley1998} that
\begin{equation}
\label{eq:distEnum}\sum_{R}q^{d(R)}=\sum_{(p_1,\ldots,p_n)}q^{p_1+\cdots+p_n},
\end{equation} where
the sum on the left is over all regions in $\mShi{n}{m}$ and the sum
on the right is over all $m$-parking functions of length $n$.


 In 2004, Athanasiadis wrote two papers on the extended Catalan
 arrangement for crystallographic $\roots$, concentrating on the
 dominant regions. The Catalan arrangement has more hyperplanes than
 the Shi arrangement, but it has the same dominant regions, so we
 record his results here in terms of the $m$-Shi arrangement.  We'll
 need a definition.  The {\em Narayana numbers} (type $A$) are given
 by \cite{petersen2015}
$$N_{n,k}=\frac{1}{k+1}\binom{n}{k}\binom{n-1}{k}.$$  They refine the Catalan numbers by counting the Dyck paths of length $n$ with $k$ peaks. In other words, $C_n=\sum_{k=0}^{n-1}N_{n,k}$. Athanasiadis
\begin{enumerate}
\item generalized and extended the Narayana numbers, finding what they enumerate in terms of dominant $m$-Shi regions;
\item counted the number of $m$-Shi regions in the dominant chamber, generalizing Shi's result described in Section~\ref{sec:shiEnum}; and
\item  \label{athan2004Ideals} used co-filtered chains of ideals in the root poset to describe the dominant $m$-Shi regions.
\end{enumerate}
We will describe \eqref{athan2004Ideals} in a bit more depth. Let $\pos=I_0\supseteq I_1\supseteq I_2\supseteq\cdots\supseteq I_m$  be a a decreasing chain $\mc{I}$ of ideals in $\pos$, set $I_i=I_m$ for all $i>m$, and set $J_i=\pos\setminus I_i$. The chain $\mc{I}$ is a {\em co-filtered chain of ideals of length $m$} if
\begin{enumerate}
\item $(I_i+I_j)\cap\pos\subseteq I_{i+j}$ and 
\item $(J_i+J_j)\cap\pos\subseteq J_{i+j}$
\end{enumerate}  
is true for all indices $i,j\geq 1$ with $i+j\leq m$. 
The coordinates of the chain are $$k_{\alpha}(\mc{I})=\max\{k_1+k_2+\cdots +k-r:\alpha=\beta_1+\cdots+\beta_r\text{ with $\beta_{i}\in I_{k_i}$ for all $i$}\}.$$
Athanasiadis showed that 
\begin{equation}
  \label{eq:mcoords}
  k_{\alpha}(\mc{I})+k_{\beta}(\mc{I})\leq k_{\alpha+\beta}(\mc{I})\leq k_{\alpha}(\mc{I})+k_{\beta}(\mc{I})+1
\end{equation}
whenever $\alpha,\beta,\alpha+\beta\in\pos$. Equation
\eqref{eq:mcoords} generalizes Shi's bijection between filters in the root poset and dominant Shi regions (see Section~\ref{sec:shiEnum}). Finally, to
  define the fundamental $m$-Shi region associated to $\mc{I}$, he
  sets $R_{\mc{I}}$ to be the set of points $x\in V$ which satisfy
\begin{enumerate}
\item $\brac{\alpha}{x}>k$, if $\alpha\in I_k$ and
\item $0<\brac{\alpha}{x}<k$, if $\alpha\in J_k$, for $0\leq k\leq m$. The coordinates of the ideal are then the coordinates of a region.
\end{enumerate}
Certain elements in an ideal are called indecomposable. These elements
correspond to the walls of $R_{\mc{I}}$ which separate $R_{\mc{I}}$
from $R_0$, and take the place of peaks in Dyck paths when defining
the Narayana numbers in terms of Shi regions.

Here we mention a few other enumerative results concerning the regions of
the extended Shi arrangement.  Any fixed
hyperplane in the $m$-Shi arrangement is dissected into regions by the
other hyperplanes in the arrangements. Fishel, Tzanaki, and Vazirani
enumerate the number of regions for certain fixed hyperplanes in type
$A$ in \cite{fishel-tzanaki-vazirani2013}. Fishel, Kallipoliti, and
Tzanaki \cite{fishel-kallipoliti-tzanaki2013} defined a bijection
between dominant regions of the $m$-Shi arrangement in type $A_n$ and
dissections of an $m(n+1)+2$-gon. These dissections represent facets
of the $m$-generalized cluster complex. In 2008, Sivasubramanian
\cite{sivasubramanian2008} gave combinatorial interpretations for the
coeffiencients of a two-variable version of Stanley's distance
enumerator \eqref{eq:distEnum} in type $A$ for $m=1$. Forge and
Zaslavsky study the integral points in $[m]^m$ that do not lie in
any hyperplane of the arrangement \cite{forge-zaslavsky2007}. Thiel
resolves a conjecture of Armstrong \cite[Conjecture
  5.1.24]{armstrong2009} on the distribution of floors and ceilings in
the dominant regions of the $m$-Shi arrangements for all types. See
Section~\ref{sec:ish} for the definition of ceiling.

\section{Connections}\label{sec:connections}
\subsection{Decompositions numbers and the Shi arrangement}
\label{subsec:Richards}
The dominant regions make an appearance in the study of decomposition
numbers for certain Hecke algebras. To describe this appearance, we'll
first need a host of combinatorial definitions, then we'll indicate
briefly how these arose from algebra, and finally we'll relate this
back to the Shi arrangement. We thank Matthew Fayers for not only
pointing out this connection, but carefully explaining it.

\subsubsection{Combinatorics}
Here we define $n$-cores, review some well-known facts about them, and review
the abacus construction, which will be useful for us.  Details can be
found in \cite{james-kerber1981}.

The {\em $(k,l)$-hook} of an integer partition $\lambda$ consists of the box in
row $k$ and column $l$ of $\lambda$, all the boxes to the right of it
in row $k$ together with all the nodes below it and in column $l$. The
{\it hook length\/} $h_{(k,l)}^\lambda$ of this box is the number of
boxes in the $(k,l)$-hook. Let $n$ be a positive integer. An {\em
  $n$-core\/} is a partition $\lambda$ such that $n \nmid
h_{(k,l)}^\lambda$ for all boxes $\lambda$. An {\em $n$-regular\/}
partition has no (nonzero) $n$ parts which equal each other. For example, $(7,6,6,6)$ is not $3$-regular. We'll
sometimes use $p$ or $e$ instead of $n$, depending on the context. The
definition is the same.

\begin{figure}[ht]
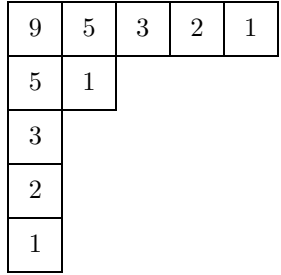


\begin{ytableau}9&5&3&2&1\\5&1\\3\\2\\1\end{ytableau}

\caption{The Young diagram of the partition $\lambda=(5,2,1,1,1)$. The
  hooklengths are the entries in the boxes of its Young diagram. The partition $\lambda$ is a $4$-core but not a $5$-core.}
\label{fig:lambda52111}  
  \end{figure}

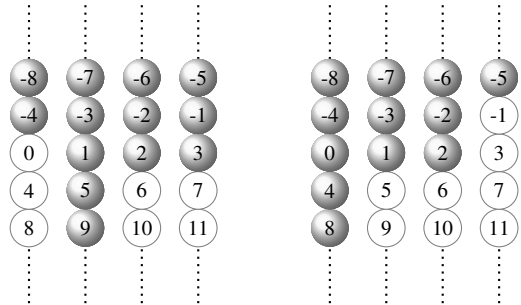
\begin{figure}[ht]
\begin{tikzpicture}
\begin{scope}[fill opacity=1,every to/.style={draw,thick, black}]
\tikzstyle{every node}=[font=\footnotesize]
\def\blankSpot{\filldraw[fill=white, draw=gray]}
\def\beadSpot{\shade[shading=ball, ball color = Gainsboro!30]}

\begin{scope}[shift={(0,0)}]
\draw[thick,black,dotted] (0,0) -- (0,4);
\draw[thick,black,dotted,shift = {(.75,0)}] (0,0) to (0,4);
\draw[thick,black,dotted,shift = {(1.5,0)}] (0,0) to (0,4);
\draw[thick,black,dotted, shift = {(2.25,0)}] (0,0) to (0,4);
\beadSpot (.75,1) circle (.25) node {9};
\beadSpot (.75,1.5) circle (.25) node{5};
\beadSpot (.75,2) circle (.25) node {1};
\beadSpot (.75,2.5) circle (.25) node {-3};
\beadSpot (.75,3) circle (.25) node {-7};
\blankSpot (1.5,1) circle (.25) node {10};
\blankSpot (1.5,1.5) circle (.25) node{6};
\beadSpot (1.5,2) circle (.25) node {2};
\beadSpot (1.5,2.5) circle (.25) node {-2};
\beadSpot (1.5,3) circle (.25) node {-6};

\blankSpot (2.25,1) circle (.25) node {11};
\blankSpot (2.25,1.5) circle (.25) node {7};
\beadSpot (2.25,2) circle (.25) node {3};
\beadSpot (2.25,2.5) circle (.25) node {-1};
\beadSpot (2.25,3) circle (.25) node {-5};

\blankSpot (0,1) circle (.25) node {8};
\blankSpot(0,1.5) circle (.25) node {4};
\blankSpot(0,2) circle (.25) node {0};
 \beadSpot (0,2.5) circle (.25) node {-4};
\beadSpot (0,3) circle (.25) node {-8};
\end{scope}
\omitt{
\begin{scope}[shift={(3,-.5)}]
\node[shift = {(-1,-1)}] {$\lambda =$};
\def\boxpath{-- +(-0.5,0) -- +(-0.5,-0.5) -- +(0,-0.5) -- cycle};
\draw[step =0.5,shift={(0,0)}, thick](0,-0.5) \boxpath node[anchor= south east ]{9};
\draw[step =0.5,shift={(0,0)}, thick](0.5,-0.5) \boxpath node[anchor= south east ]{5};
\draw[step =0.5,shift={(0,0)}, thick](1,-0.5) \boxpath node[anchor= south east ]{3};
\draw[step =0.5,shift={(0,0)}, thick](1.5,-0.5) \boxpath node[anchor= south east ]{2};
\draw[step =0.5,shift={(0,0)}, thick](2,-0.5) \boxpath node[anchor= south east ]{1};
\draw[step =0.5,shift={(0,0)}, thick](0,-1) \boxpath node[anchor= south east ]{5};
\draw[step =0.5,shift={(0,0)}, thick](0.5,-1) \boxpath node[anchor= south east ]{1};
\draw[step =0.5,shift={(0,0)}, thick](0,-1.5) \boxpath node[anchor= south east ]{3};
\draw[step =0.5,shift={(0,0)}, thick](0,-2) \boxpath node[anchor= south east ]{2};
\draw[step =0.5,shift={(0,0)}, thick](0,-2.5) \boxpath node[anchor= south east ]{1};

\end{scope}
}
\end{scope}

  \begin{scope}[shift={(4,0)}]

\begin{scope}[fill opacity=1,every to/.style={draw,thick, dashed}]
\tikzstyle{every node}=[font=\footnotesize]
\draw[thick,dotted] (0,0) to (0,4);
\draw[thick,dotted,shift = {(.75,0)}] (0,0) to (0,4);
\draw[thick,dotted,shift = {(1.5,0)}] (0,0) to (0,4);
\draw[thick,dotted,shift = {(2.25,0)}] (0,0) to (0,4);
\shade[shading=ball, ball color = Gainsboro!30] (0,1) circle (.25) node {8};
\shade[shading=ball, ball color = Gainsboro!30] (0,1.5) circle (.25) node
{4};
\shade[shading=ball, ball color = Gainsboro!30] (0,2) circle (.25) node {0};
\shade[shading=ball, ball color = Gainsboro!30] (0,2.5) circle (.25) node {-4};
\shade[shading=ball, ball color = Gainsboro!30] (0,3) circle (.25) node {-8};
\filldraw[fill=white, draw=gray] (.75,1) circle (.25) node {9};
\filldraw[fill=white, draw=gray](.75,1.5) circle (.25) node
{5};
\shade[shading=ball, ball color = Gainsboro!30] (.75,2) circle (.25) node {1};
\shade[shading=ball, ball color = Gainsboro!30] (.75,2.5) circle (.25) node {-3};
\shade[shading=ball, ball color = Gainsboro!30] (.75,3) circle (.25) node {-7};

\filldraw[fill=white, draw=gray](1.5,1) circle (.25) node {10};
\filldraw[fill=white, draw=gray](1.5,1.5) circle (.25) node {6};
\shade[shading=ball, ball color = Gainsboro!30] (1.5,2) circle (.25) node {2};
\shade[shading=ball, ball color = Gainsboro!30] (1.5,2.5) circle (.25) node {-2};
\shade[shading=ball, ball color = Gainsboro!30] (1.5,3) circle (.25) node {-6};

\filldraw[fill=white, draw=gray](2.25,1) circle (.25) node {11};
\filldraw[fill=white, draw=gray] (2.25,1.5) circle (.25) node
{7};
\filldraw[fill=white, draw=gray](2.25,2) circle (.25) node {3};
\filldraw[fill=white, draw=gray] (2.25,2.5) circle (.25) node {-1};
\shade[shading=ball, ball color = Gainsboro!30] (2.25,3) circle (.25) node {-5};
\end{scope}

  \end{scope}
\end{tikzpicture}

\caption{On the left, the beads are placed on the positions labeled by
  the first column hook-lengths of $\lambda=(5,2,1,1,1)$. On the
  right is an equivalent abacus, where all bead positions have been
  shifted by $C=-1$.}
\label{fig:aba52111}

\end{figure}

\omitt{We'll use a simplified version of $\beta$-numbers, which is
  enough for our purposes.} The {\em $\beta$-numbers} of the partition
$\lambda$ are the hook lengths from its first column. The
$\beta$-numbers can be displayed on an abacus: a {\em $p$-abacus} is a
diagram with $p$ runners, labeled $0, 1,\ldots,p-1$. Runner $i$ has
positions labeled by integers $pj+i$, for all $j\in\Z$.  We make a
$p$-abacus for $\lambda$ by placing a bead at position $\beta_k$, for
each $\beta$-number $\beta_k$ of $\lambda$ and at all negative
positions. We say two $p$-abaci are equivalent if we can change one to
the other by moving the bead at position $i$ to position $i+C$ for
some $C\in\Z$ and for all positions $i$ where there is a bead.  The
positive integer $p$ is arbitrary for now, but will be related to the
characteristic of a field when we see abaci in their algebraic
context. See the $4$-abacus in Figure~\ref{fig:aba52111}.

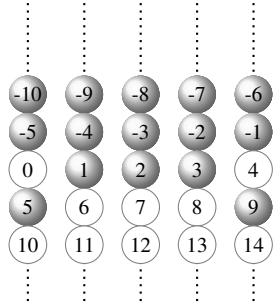
\begin{figure}[ht]
\begin{tikzpicture}
\begin{scope}[fill opacity=1,every to/.style={draw,dotted},
 information text/.style={rounded corners,fill=red!10,inner sep=1ex}]
\tikzstyle{every node}=[font=\footnotesize]
\def\blankSpot{\filldraw[fill=white, draw=gray]}
\def\beadSpot{\shade[shading=ball, ball color = Gainsboro!30]}

\draw[thick,black,dotted,shift={(0,0)}] (0,1) to (0,5.0);
\draw[thick,black,dotted,shift={(0.75,0)}] (0,1) to (0,5.0);
\draw[thick,black,dotted,shift={(1.5,0)}] (0,1) to (0,5.0);
\draw[thick,black,dotted,shift={(2.25,0)}] (0,1) to (0,5.0);
\draw[thick,black,dotted,shift={(3,0)}] (0,1) to (0,5.0);

\beadSpot (0,3.75) circle (0.25) node {-10};
\beadSpot (0.75,3.75) circle (0.25) node {-9};
\beadSpot (1.5,3.75) circle (0.25) node {-8};
\beadSpot (2.25,3.75) circle (0.25) node {-7};
\beadSpot (3,3.75) circle (0.25) node {-6};

\beadSpot (0,3.25) circle (0.25) node {-5};
\beadSpot (0.75,3.25) circle (0.25) node {-4};
\beadSpot (1.5,3.25) circle (0.25) node {-3};
\beadSpot (2.25,3.25) circle (0.25) node {-2};
\beadSpot (3,3.25) circle (0.25) node {-1};

\blankSpot (0,2.75) circle (0.25) node {0};
\beadSpot (0.75,2.75) circle (0.25) node {1};
\beadSpot (1.5,2.75) circle (0.25) node {2};
\beadSpot (2.25,2.75) circle (0.25) node {3};
\blankSpot (3,2.75) circle (0.25) node {4};

\beadSpot (0,2.25) circle (0.25) node {5};
\blankSpot (0.75,2.25) circle (0.25) node {6};
\blankSpot (1.5,2.25) circle (0.25) node {7};
\blankSpot (2.25,2.25) circle (0.25) node {8};
\beadSpot (3,2.25) circle (0.25) node {9};

\blankSpot (0,1.75) circle (0.25) node {10};
\blankSpot (0.75,1.75) circle (0.25) node {11};
\blankSpot (1.5,1.75) circle (0.25) node {12};
\blankSpot (2.25,1.75) circle (0.25) node {13};
\blankSpot (3,1.75) circle (0.25) node {14};

\omitt{
\blankSpot (0,1.25) circle (0.25) node {15};
\blankSpot (0.75,1.25) circle (0.25) node {16};
\blankSpot (1.5,1.25) circle (0.25) node {17};
\blankSpot (2.25,1.25) circle (0.25) node {18};
\blankSpot (3,1.25) circle (0.25) node {19};
}

\end{scope}

\end{tikzpicture}
\caption{The beads are placed on a $4$-abacus on the positions labeled by
  the first column hook-lengths of $\lambda=(5,2,1,1,1)$. The partition $\lambda$ is not a $5$-core, since there are gaps on the $5$-abacus. If we push the beads in positions $5$ and 9 up to fill in the gaps, we obtain the $5$-abacus of the empty partition. }
\end{figure}

We can give an equivalent description of $p$-core partitions: a
partition $\lambda$ is an $p$-core if and only if whenever there is a
bead at position $j$ of its $p$-abacus, there is also a bead at
position $j-p$ \cite{james-kerber1981}.

Suppose we have a partition which is not an $p$-core. Then there is at
least one bead at a position $j$ of its $p$-abacus which can be pushed
up into the vacant position $j-p$. This gives us the $p$-abacus of
another partition. We can repeat this until no beads can be pushed up,
at which point we have the $p$-abacus of an $p$-core. The final
partition $\gamma$ is called the $p$-core of $\lambda$ and the number of beads
we moved in called the $p$-weight of $\lambda$. With a little work, which we won't do, it is possible to show that $|\lambda|=|\gamma|+wp$.

The last combinatorial ingredient we need are residues for the boxes
of a partition. Let $n$ be a positive integer. We call the box in row
$i$, column $j$ a {\em $k$-box} if $(j-i)\bmod n$ is equal to $k$.

\begin{figure}[ht]

\begin{ytableau}0&1&2&3&0\\3&0\\2\\1\\0\end{ytableau}

\caption{The Young diagram of the partition $\lambda=(5,2,1,1,1)$. The
  entries are the residues modulo $4$.}
\label{fig:lambda52111Res}  
  \end{figure}

\subsubsection{Algebra}
\label{sec:algConn}
It is time to say a few words about how the combinatorics from the last section relate to  algebra.

Roughly speaking, the $p$-regular partitions of $n$ index the
irreducible modules $D^{\lambda}$ of $F\S_n$, where $p$ is the
characteristic of the field $F$. The {\em $p$-blocks} are the
equivalence classes of a certain equivalence relation on the
irreducible modules. By Nakayama's celebrated conjecture, and Brauer
and Robinson's theorem, two irreducibles $D^{\lambda}$ and $D^{\mu}$
belong to the same $p$-block if and only $\lambda$ and $\mu$ have the
same $p$-core. Thus each $p$-block is labeled by a $p$-core and
$p$-weight $w$. The $p$-weight keeps track of the difference between the
$p$-core labeling a block and the partitions in the block: the
$p$-core is a partition of $n-pw$. The weights are nonnegative
integers. We will assume $w$ is at least one, because the blocks where
$w=0$ are singletons consisting of the $p$-core.

\begin{figure}[ht]
\begin{tikzpicture}
\begin{scope}[fill opacity=1,every to/.style={draw,thick, black}]
\tikzstyle{every node}=[font=\footnotesize]
\def\blankSpot{\filldraw[fill=white, draw=gray]}
\def\beadSpot{\shade[shading=ball, ball color = Gainsboro!30]}

\begin{scope}[shift={(0,0)}]
\draw[thick,black,dotted] (0,0) -- (0,4);
\draw[thick,black,dotted,shift = {(.75,0)}] (0,0) to (0,4);
\draw[thick,black,dotted,shift = {(1.5,0)}] (0,0) to (0,4);
\draw[thick,black,dotted, shift = {(2.25,0)}] (0,0) to (0,4);
\blankSpot (.75,1) circle (.25) node {9};
\blankSpot (.75,1.5) circle (.25) node{5};
\beadSpot (.75,2) circle (.25) node {1};
\beadSpot (.75,2.5) circle (.25) node {-3};
\beadSpot (.75,3) circle (.25) node {-7};
\beadSpot (1.5,1) circle (.25) node {10};
\beadSpot (1.5,1.5) circle (.25) node{6};
\beadSpot (1.5,2) circle (.25) node {2};
\beadSpot (1.5,2.5) circle (.25) node {-2};
\beadSpot (1.5,3) circle (.25) node {-6};

\blankSpot (2.25,1) circle (.25) node {11};
\blankSpot (2.25,1.5) circle (.25) node {7};
\beadSpot (2.25,2) circle (.25) node {3};
\beadSpot (2.25,2.5) circle (.25) node {-1};
\beadSpot (2.25,3) circle (.25) node {-5};

\blankSpot (0,1) circle (.25) node {8};
\blankSpot(0,1.5) circle (.25) node {4};
\blankSpot(0,2) circle (.25) node {0};
 \beadSpot (0,2.5) circle (.25) node {-4};
\beadSpot (0,3) circle (.25) node {-8};
\end{scope}
\omitt{
\begin{scope}[shift={(3,-.5)}]
\node[shift = {(-1,-1)}] {$\lambda =$};
\def\boxpath{-- +(-0.5,0) -- +(-0.5,-0.5) -- +(0,-0.5) -- cycle};
\draw[step =0.5,shift={(0,0)}, thick](0,-0.5) \boxpath node[anchor= south east ]{9};
\draw[step =0.5,shift={(0,0)}, thick](0.5,-0.5) \boxpath node[anchor= south east ]{5};
\draw[step =0.5,shift={(0,0)}, thick](1,-0.5) \boxpath node[anchor= south east ]{3};
\draw[step =0.5,shift={(0,0)}, thick](1.5,-0.5) \boxpath node[anchor= south east ]{2};
\draw[step =0.5,shift={(0,0)}, thick](2,-0.5) \boxpath node[anchor= south east ]{1};
\draw[step =0.5,shift={(0,0)}, thick](0,-1) \boxpath node[anchor= south east ]{5};
\draw[step =0.5,shift={(0,0)}, thick](0.5,-1) \boxpath node[anchor= south east ]{1};
\draw[step =0.5,shift={(0,0)}, thick](0,-1.5) \boxpath node[anchor= south east ]{3};
\draw[step =0.5,shift={(0,0)}, thick](0,-2) \boxpath node[anchor= south east ]{2};
\draw[step =0.5,shift={(0,0)}, thick](0,-2.5) \boxpath node[anchor= south east ]{1};

\end{scope}
}
\end{scope}

  \begin{scope}[shift={(4,0)}]
\begin{scope}[fill opacity=1,every to/.style={draw,thick, black}]
\tikzstyle{every node}=[font=\footnotesize]
\def\blankSpot{\filldraw[fill=white, draw=gray]}
\def\beadSpot{\shade[shading=ball, ball color = Gainsboro!30]}

\begin{scope}[shift={(0,0)}]
\draw[thick,black,dotted] (0,0) -- (0,4);
\draw[thick,black,dotted,shift = {(.75,0)}] (0,0) to (0,4);
\draw[thick,black,dotted,shift = {(1.5,0)}] (0,0) to (0,4);
\draw[thick,black,dotted, shift = {(2.25,0)}] (0,0) to (0,4);
\beadSpot (.75,1) circle (.25) node {9};
\beadSpot (.75,1.5) circle (.25) node{5};
\beadSpot (.75,2) circle (.25) node {1};
\beadSpot (.75,2.5) circle (.25) node {-3};
\beadSpot (.75,3) circle (.25) node {-7};
\blankSpot (1.5,1) circle (.25) node {10};
\blankSpot (1.5,1.5) circle (.25) node{6};
\beadSpot (1.5,2) circle (.25) node {2};
\beadSpot (1.5,2.5) circle (.25) node {-2};
\beadSpot (1.5,3) circle (.25) node {-6};

\blankSpot (2.25,1) circle (.25) node {11};
\blankSpot (2.25,1.5) circle (.25) node {7};
\beadSpot (2.25,2) circle (.25) node {3};
\beadSpot (2.25,2.5) circle (.25) node {-1};
\beadSpot (2.25,3) circle (.25) node {-5};

\blankSpot (0,1) circle (.25) node {8};
\blankSpot(0,1.5) circle (.25) node {4};
\blankSpot(0,2) circle (.25) node {0};
 \beadSpot (0,2.5) circle (.25) node {-4};
\beadSpot (0,3) circle (.25) node {-8};
\end{scope}
\omitt{
\begin{scope}[shift={(3,-.5)}]
\node[shift = {(-1,-1)}] {$\lambda =$};
\def\boxpath{-- +(-0.5,0) -- +(-0.5,-0.5) -- +(0,-0.5) -- cycle};
\draw[step =0.5,shift={(0,0)}, thick](0,-0.5) \boxpath node[anchor= south east ]{9};
\draw[step =0.5,shift={(0,0)}, thick](0.5,-0.5) \boxpath node[anchor= south east ]{5};
\draw[step =0.5,shift={(0,0)}, thick](1,-0.5) \boxpath node[anchor= south east ]{3};
\draw[step =0.5,shift={(0,0)}, thick](1.5,-0.5) \boxpath node[anchor= south east ]{2};
\draw[step =0.5,shift={(0,0)}, thick](2,-0.5) \boxpath node[anchor= south east ]{1};
\draw[step =0.5,shift={(0,0)}, thick](0,-1) \boxpath node[anchor= south east ]{5};
\draw[step =0.5,shift={(0,0)}, thick](0.5,-1) \boxpath node[anchor= south east ]{1};
\draw[step =0.5,shift={(0,0)}, thick](0,-1.5) \boxpath node[anchor= south east ]{3};
\draw[step =0.5,shift={(0,0)}, thick](0,-2) \boxpath node[anchor= south east ]{2};
\draw[step =0.5,shift={(0,0)}, thick](0,-2.5) \boxpath node[anchor= south east ]{1};

\end{scope}
}
\end{scope}

  \end{scope}
\end{tikzpicture}

\caption{On the left, the abacus for $\gamma=(6,3,1,1,1)$, on the right, the abacus for $\bar{\gamma}=(5,2,1,1,1)$. We have moved $k=2$ beads from runner 2 to runner 1.  }
\label{fig:scopes}

\end{figure}
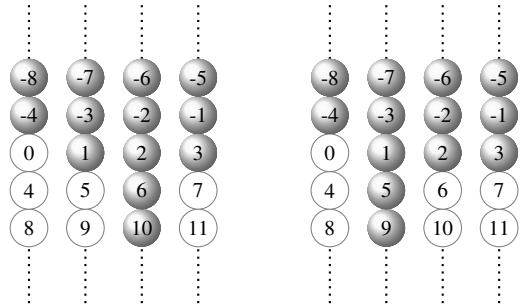

Scopes \cite{scopes1991} investigated the classes of $p$-blocks under
Morita equivalence. She characterized families of Morita equivalent
$p$-blocks using the $p$-core and $p$-weight which label
$p$-blocks. Suppose $B$ is $p$-block for $F\Sn$ labeled by weight $w$
and $p$-core $\gamma$. Let $k$ be an integer at least as large as $w$,
and suppose that in a $p$-abacus for $\gamma$, there is a runner $i$
which has $k$ more beads in positive positions than runner $(i-1)$
has. Now move $k$ beads from runner $i$ to runner $i-1$ in such a
$p$-abacus for $\gamma$. This new abacus determines another $p$-core,
say $\bar{\gamma}$. See Figure~\ref{fig:scopes}. The operation changes the size of the partition:
$|\gamma|-k=|\bar{\gamma}|$. We are glossing over details here involving the
$\beta$-numbers. Let $\bar{B}$ be the block of $F\Sym{n-k}$ labeled by
$w$ and $\bar{\gamma}$. Scopes took the transitive closure of the relation
$B\sim\bar{B}$ and showed that within an equivalence class, the
$p$-blocks have the same decomposition matrix, among other results.


Richards was studying the decomposition numbers for the Hecke algebra;
see Section~\ref{subsec:KL} for at least the definition.  He used the
classes from the Scopes equivalence on $e$-cores, where $e$ depends on
the characteristic $p$ of the field and the element $q$ used in the
definition of the Hecke algebra. He called the classes {\em families}.
Richards was interested in these families because the blocks of the
Hecke algebras for $\Sn$ and for $\Sym{n-k}$ corresponding to $\gamma$
and $\bar{\gamma}$ respectively have essentially the same
decomposition numbers \cite{richards1996}.

Richards wanted to count such families. He built the following {\em
  pyramid} $\{_ua_v\}_{0\leq u<v\leq e-1}$ for an $e$-core $\gamma$
based on $\gamma$'s $e$-abacus. Note the similarity in shape to
admissible sign types and to the arrangement of roots in a staircase
shape diagram in Figure~\ref{fig:catObjs}.  For $i=0,1,\ldots, e-1$,
let $p'_i$ be the position of the first free space on runner
$i$. Arrange these $e$ numbers in ascending order and relabel as
$p_0<p_1<\cdots< p_{e-1}$. If $0\leq u<v\leq e-1$, then $p_u-p_v$ is a
positive integer not divisible by $e$. We may use any $e$-abacus for
$\gamma$; it doesn't affect the set of differences. Richards defined
the pyramid of numbers by


$$\acoord{u}{v}=\begin{cases}
w-1&\text{if $0<p_v-p_v<e$}\\ 
w-2&\text{if $e<p_v-p_v<2e$}\\
&\vdots\\
1&\text{if $(w-2)e<p_v-p_v<(w-1)e$}\\
0&\text{if $(w-1)e<p_v-p_v$}\\
\end{cases}$$

where $w$ is the weight.
Richards proved that two $e$-cores are in the same family if and only
if they have the same pyramid and that there are
exactly $$\frac{1}{e}\binom{ew}{e-1}$$ families. What's more, he
characterized the triangles of numbers which form a pyramid.

To show the connection to the Shi arrangement, we  transform $\acoord{u}{v}$ into $\aprimecoord{u}{v}$ by $\acoord{u}{v}+\aprimecoord{u}{v}=w-1$. Then Richard's Proposition~3.4 becomes

\begin{proposition}[\cite{richards1996}] Let $e\geq 2$ and $w>0$, and for $0\leq u<v<e-1$ let $0\leq \aprimecoord{u}{v}\leq w-1$. Then the $\aprimecoord{u}{v}$ form a pyramid if and only if for all $0\leq u<t<v\leq e-1$,
$$\aprimecoord{u}{v}=\begin{cases}\aprimecoord{u}{t}+\aprimecoord{t}{v}\text{ or }\aprimecoord{u}{t}+\aprimecoord{t}{v}+1&\text{if both of these have all entries no bigger than $w-1$}\\w-1&\text{otherwise.}\end{cases}$$
\end{proposition}

Please see Example~\ref{ex:richards} for the calculation of a few pyramids from cores. 

Let's examine the case $e=3$ and $w=2$. There are five families. We
choose \omitt{the }
five $3$-cores $\{\emptyset, (1),(2),(1,1), (3,1,1)\}$ and
calculate their pyramids: \ytableausetup{smalltableaux}

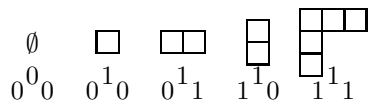
\begin{center}
\begin{tikzpicture}[scale=.2]
\def\a{5}
  \def\b{3.2}

  \begin{scope}[shift={(0,0)}]
\node at (1,\b) {$\emptyset$};
    \Triangle{$0$}{$0$}{$0$}
  \end{scope} 

  \begin{scope}[shift={(\a,0)}]
\node at (1,\b) {\ydiagram{1}};
    \Triangle{$1$}{$0$}{$0$}
  \end{scope}

  \begin{scope}[shift={(2*\a,0)}]
\node at (1,\b) {\ydiagram{2}};
    \Triangle{$1$}{$0$}{$1$}
  \end{scope}

  \begin{scope}[shift={(3*\a,0)}]
\node at (1,\b) {\ydiagram{1,1}};
    \Triangle{$1$}{$1$}{$0$}
  \end{scope}

  \begin{scope}[shift={(4*\a,0)}]
\node at (1,\b) {\ydiagram{3,1,1}};
    \Triangle{$1$}{$1$}{$1$}
  \end{scope}

\end{tikzpicture}
\end{center}


The pyramids are all different, so we have found all the families. If
we look back at the set $G$ in Section~\ref{subsec:cells} and consider
the subset where all entries are either $+$ or $\bigcirc$, we see a
similarity to the pyramids (replace $+$ with $1$). This is true in
general. Richard's proposition is the type $A$ version of
\eqref{eq:mcoords} from Section~\ref{sec:extShi}. The pyramid is also
an admissible sign type for a dominant $m$-Shi region of type $A$,
where $m=w-1$.

\subsubsection{Geometry}
\begin{figure}[h]
\begin{tikzpicture}
\begin{scope}[fill opacity=.5]
\tikzstyle{every node}=[font=\footnotesize]
[0,0,0]
\filldraw [draw=gray,thick,fill opacity = .5, fill = yellow,shift={(0.000000,0)}](0,0) -- (60:2cm) -- (0:2cm) --(0,0) node[blue] at (0.8cm,0.4cm) {} node[black,fill opacity=1, draw opacity = 1] at (1cm,0.7cm) {\tiny{$\emptyset$}} node[blue] at (0.8cm,1.2cm) {};

\def\boxpath{-- +(-0.154701,0) -- +(-0.154701,-0.154701) -- +(0,-0.154701) -- cycle};
[1,0,0]
\filldraw [draw=gray, thick,fill = yellow, shift={(2,0)}](0,0) -- (60:2cm) -- (120:2cm) --(0,0) node[green,fill opacity=1] at (0cm,1.6cm) {} node[black,fill opacity=1, draw opacity = 1] at (0cm,1.2cm) {} node[green,fill opacity=1] at (0cm,0.8cm) {};

\def\boxpath{-- +(-0.154701,0) -- +(-0.154701,-0.154701) -- +(0,-0.154701) -- cycle};
\draw[step =0.154701,shift={(1.9359,1.15855)}](0,-0.154701) \boxpath;
[1,0,1]
\filldraw [draw=gray,thick,fill opacity = .5, fill = yellow,shift={(2.000000,0)}](0,0) -- (60:2cm) -- (0:2cm) --(0,0) node[blue] at (0.8cm,0.4cm) {} node[black,fill opacity=1, draw opacity = 1] at (1cm,0.7cm) {} node[blue] at (0.8cm,1.2cm) {};

\def\boxpath{-- +(-0.154701,0) -- +(-0.154701,-0.154701) -- +(0,-0.154701) -- cycle};
\draw[step =0.154701,shift={(2.8359,1.1282)}](0,-0.154701) \boxpath;
\draw[step =0.154701,shift={(2.8359,1.1282)}](0.154701,-0.154701) \boxpath;

[2,0,1]

\filldraw [draw=gray, thick,fill = pink!30, shift={(4,0)}](0,0) -- (60:2cm) -- (120:2cm) --(0,0) node[green,fill opacity=1] at (0cm,1.6cm) {} node[black,fill opacity=1, draw opacity = 1] at (0cm,1.2cm) {} node[green,fill opacity=1] at (0cm,0.8cm) {};
\def\boxpath{-- +(-0.154701,0) -- +(-0.154701,-0.154701) -- +(0,-0.154701) -- cycle};
\draw[step =0.154701,shift={(3.9359,1.31325)}](0,-0.154701) \boxpath;
\draw[step =0.154701,shift={(3.9359,1.31325)}](0.154701,-0.154701) \boxpath;
\draw[step =0.154701,shift={(3.9359,1.31325)}](0.309401,-0.154701) \boxpath;
\draw[step =0.154701,shift={(3.9359,1.31325)}](0,-0.309401) \boxpath;
[2,0,2]
\filldraw [draw=gray,thick,fill opacity = .5, fill = pink!30,shift={(4.000000,0)}](0,0) -- (60:2cm) -- (0:2cm) --(0,0) node[blue] at (0.8cm,0.4cm) {} node[black,fill opacity=1, draw opacity = 1] at (1cm,0.7cm) {} node[blue] at (0.8cm,1.2cm) {};
\def\boxpath{-- +(-0.154701,0) -- +(-0.154701,-0.154701) -- +(0,-0.154701) -- cycle};
\draw[step =0.154701,shift={(4.8359,1.1282)}](0,-0.154701) \boxpath;
\draw[step =0.154701,shift={(4.8359,1.1282)}](0.154701,-0.154701) \boxpath;
\draw[step =0.154701,shift={(4.8359,1.1282)}](0.309401,-0.154701) \boxpath;
\draw[step =0.154701,shift={(4.8359,1.1282)}](0.464102,-0.154701) \boxpath;
\draw[step =0.154701,shift={(4.8359,1.1282)}](0,-0.309401) \boxpath;
\draw[step =0.154701,shift={(4.8359,1.1282)}](0.154701,-0.309401) \boxpath;
[3,0,2]
\filldraw [draw=gray, thick,fill = white, shift={(6,0)}](0,0) -- (60:2cm) -- (120:2cm) --(0,0) node[green,fill opacity=1] at (0cm,1.6cm) {} node[black,fill opacity=1, draw opacity = 1] at (0cm,1.2cm) {} node[green,fill opacity=1] at (0cm,0.8cm) {};
\def\boxpath{-- +(-0.1326,0) -- +(-0.1326,-0.1326) -- +(0,-0.1326) -- cycle};
\draw[step =0.1326,shift={(5.9359,1.40165)}](0,-0.1326) \boxpath;
\draw[step =0.1326,shift={(5.9359,1.40165)}](0.1326,-0.1326) \boxpath;
\draw[step =0.1326,shift={(5.9359,1.40165)}](0.265201,-0.1326) \boxpath;
\draw[step =0.1326,shift={(5.9359,1.40165)}](0.397801,-0.1326) \boxpath;
\draw[step =0.1326,shift={(5.9359,1.40165)}](0.530402,-0.1326) \boxpath;
\draw[step =0.1326,shift={(5.9359,1.40165)}](0,-0.265201) \boxpath;
\draw[step =0.1326,shift={(5.9359,1.40165)}](0.1326,-0.265201) \boxpath;
\draw[step =0.1326,shift={(5.9359,1.40165)}](0.265201,-0.265201) \boxpath;
\draw[step =0.1326,shift={(5.9359,1.40165)}](0,-0.397801) \boxpath;
[3,0,3]
\filldraw [draw=gray,thick,fill opacity = .5, fill = white,shift={(6.000000,0)}](0,0) -- (60:2cm) -- (0:2cm) --(0,0) node[blue] at (0.8cm,0.4cm) {} node[black,fill opacity=1, draw opacity = 1] at (1cm,0.7cm) {} node[blue] at (0.8cm,1.2cm) {};
\def\boxpath{-- +(-0.116025,0) -- +(-0.116025,-0.116025) -- +(0,-0.116025) -- cycle};
\draw[step =0.116025,shift={(6.8359,1.1282)}](0,-0.116025) \boxpath;
\draw[step =0.116025,shift={(6.8359,1.1282)}](0.116025,-0.116025) \boxpath;
\draw[step =0.116025,shift={(6.8359,1.1282)}](0.232051,-0.116025) \boxpath;
\draw[step =0.116025,shift={(6.8359,1.1282)}](0.348076,-0.116025) \boxpath;
\draw[step =0.116025,shift={(6.8359,1.1282)}](0.464102,-0.116025) \boxpath;
\draw[step =0.116025,shift={(6.8359,1.1282)}](0.580127,-0.116025) \boxpath;
\draw[step =0.116025,shift={(6.8359,1.1282)}](0,-0.232051) \boxpath;
\draw[step =0.116025,shift={(6.8359,1.1282)}](0.116025,-0.232051) \boxpath;
\draw[step =0.116025,shift={(6.8359,1.1282)}](0.232051,-0.232051) \boxpath;
\draw[step =0.116025,shift={(6.8359,1.1282)}](0.348076,-0.232051) \boxpath;
\draw[step =0.116025,shift={(6.8359,1.1282)}](0,-0.348076) \boxpath;
\draw[step =0.116025,shift={(6.8359,1.1282)}](0.116025,-0.348076) \boxpath;
[4,0,3]
\filldraw [draw=gray, thick,fill = white, shift={(8,0)}](0,0) -- (60:2cm) -- (120:2cm) --(0,0) node[green,fill opacity=1] at (0cm,1.6cm) {} node[black,fill opacity=1, draw opacity = 1] at (0cm,1.2cm) {} node[green,fill opacity=1] at (0cm,0.8cm) {};
\def\boxpath{-- +(-0.103134,0) -- +(-0.103134,-0.103134) -- +(0,-0.103134) -- cycle};
\draw[step =0.103134,shift={(7.9359,1.41638)}](0,-0.103134) \boxpath;
\draw[step =0.103134,shift={(7.9359,1.41638)}](0.103134,-0.103134) \boxpath;
\draw[step =0.103134,shift={(7.9359,1.41638)}](0.206267,-0.103134) \boxpath;
\draw[step =0.103134,shift={(7.9359,1.41638)}](0.309401,-0.103134) \boxpath;
\draw[step =0.103134,shift={(7.9359,1.41638)}](0.412535,-0.103134) \boxpath;
\draw[step =0.103134,shift={(7.9359,1.41638)}](0.515668,-0.103134) \boxpath;
\draw[step =0.103134,shift={(7.9359,1.41638)}](0.618802,-0.103134) \boxpath;
\draw[step =0.103134,shift={(7.9359,1.41638)}](0,-0.206267) \boxpath;
\draw[step =0.103134,shift={(7.9359,1.41638)}](0.103134,-0.206267) \boxpath;
\draw[step =0.103134,shift={(7.9359,1.41638)}](0.206267,-0.206267) \boxpath;
\draw[step =0.103134,shift={(7.9359,1.41638)}](0.309401,-0.206267) \boxpath;
\draw[step =0.103134,shift={(7.9359,1.41638)}](0.412535,-0.206267) \boxpath;
\draw[step =0.103134,shift={(7.9359,1.41638)}](0,-0.309401) \boxpath;
\draw[step =0.103134,shift={(7.9359,1.41638)}](0.103134,-0.309401) \boxpath;
\draw[step =0.103134,shift={(7.9359,1.41638)}](0.206267,-0.309401) \boxpath;
\draw[step =0.103134,shift={(7.9359,1.41638)}](0,-0.412535) \boxpath;
[4,0,4]
\filldraw [draw=gray,thick,fill opacity = .5, fill = white,shift={(8.000000,0)}](0,0) -- (60:2cm) -- (0:2cm) --(0,0) node[blue] at (0.8cm,0.4cm) {} node[black,fill opacity=1, draw opacity = 1] at (1cm,0.7cm) {} node[blue] at (0.8cm,1.2cm) {};
\def\boxpath{-- +(-0.0928203,0) -- +(-0.0928203,-0.0928203) -- +(0,-0.0928203) -- cycle};
\draw[step =0.0928203,shift={(8.8359,1.1282)}](0,-0.0928203) \boxpath;
\draw[step =0.0928203,shift={(8.8359,1.1282)}](0.0928203,-0.0928203) \boxpath;
\draw[step =0.0928203,shift={(8.8359,1.1282)}](0.185641,-0.0928203) \boxpath;
\draw[step =0.0928203,shift={(8.8359,1.1282)}](0.278461,-0.0928203) \boxpath;
\draw[step =0.0928203,shift={(8.8359,1.1282)}](0.371281,-0.0928203) \boxpath;
\draw[step =0.0928203,shift={(8.8359,1.1282)}](0.464102,-0.0928203) \boxpath;
\draw[step =0.0928203,shift={(8.8359,1.1282)}](0.556922,-0.0928203) \boxpath;
\draw[step =0.0928203,shift={(8.8359,1.1282)}](0.649742,-0.0928203) \boxpath;
\draw[step =0.0928203,shift={(8.8359,1.1282)}](0,-0.185641) \boxpath;
\draw[step =0.0928203,shift={(8.8359,1.1282)}](0.0928203,-0.185641) \boxpath;
\draw[step =0.0928203,shift={(8.8359,1.1282)}](0.185641,-0.185641) \boxpath;
\draw[step =0.0928203,shift={(8.8359,1.1282)}](0.278461,-0.185641) \boxpath;
\draw[step =0.0928203,shift={(8.8359,1.1282)}](0.371281,-0.185641) \boxpath;
\draw[step =0.0928203,shift={(8.8359,1.1282)}](0.464102,-0.185641) \boxpath;
\draw[step =0.0928203,shift={(8.8359,1.1282)}](0,-0.278461) \boxpath;
\draw[step =0.0928203,shift={(8.8359,1.1282)}](0.0928203,-0.278461) \boxpath;
\draw[step =0.0928203,shift={(8.8359,1.1282)}](0.185641,-0.278461) \boxpath;
\draw[step =0.0928203,shift={(8.8359,1.1282)}](0.278461,-0.278461) \boxpath;
\draw[step =0.0928203,shift={(8.8359,1.1282)}](0,-0.371281) \boxpath;
\draw[step =0.0928203,shift={(8.8359,1.1282)}](0.0928203,-0.371281) \boxpath;
[5,0,4]
\filldraw [draw=gray, thick,fill = white, shift={(10,0)}](0,0) -- (60:2cm) -- (120:2cm) --(0,0) node[green,fill opacity=1] at (0cm,1.6cm) {} node[black,fill opacity=1, draw opacity = 1] at (0cm,1.2cm) {} node[green,fill opacity=1] at (0cm,0.8cm) {};
\def\boxpath{-- +(-0.0843821,0) -- +(-0.0843821,-0.0843821) -- +(0,-0.0843821) -- cycle};
\draw[step =0.0843821,shift={(9.9359,1.42576)}](0,-0.0843821) \boxpath;
\draw[step =0.0843821,shift={(9.9359,1.42576)}](0.0843821,-0.0843821) \boxpath;
\draw[step =0.0843821,shift={(9.9359,1.42576)}](0.168764,-0.0843821) \boxpath;
\draw[step =0.0843821,shift={(9.9359,1.42576)}](0.253146,-0.0843821) \boxpath;
\draw[step =0.0843821,shift={(9.9359,1.42576)}](0.337528,-0.0843821) \boxpath;
\draw[step =0.0843821,shift={(9.9359,1.42576)}](0.421911,-0.0843821) \boxpath;
\draw[step =0.0843821,shift={(9.9359,1.42576)}](0.506293,-0.0843821) \boxpath;
\draw[step =0.0843821,shift={(9.9359,1.42576)}](0.590675,-0.0843821) \boxpath;
\draw[step =0.0843821,shift={(9.9359,1.42576)}](0.675057,-0.0843821) \boxpath;
\draw[step =0.0843821,shift={(9.9359,1.42576)}](0,-0.168764) \boxpath;
\draw[step =0.0843821,shift={(9.9359,1.42576)}](0.0843821,-0.168764) \boxpath;
\draw[step =0.0843821,shift={(9.9359,1.42576)}](0.168764,-0.168764) \boxpath;
\draw[step =0.0843821,shift={(9.9359,1.42576)}](0.253146,-0.168764) \boxpath;
\draw[step =0.0843821,shift={(9.9359,1.42576)}](0.337528,-0.168764) \boxpath;
\draw[step =0.0843821,shift={(9.9359,1.42576)}](0.421911,-0.168764) \boxpath;
\draw[step =0.0843821,shift={(9.9359,1.42576)}](0.506293,-0.168764) \boxpath;
\draw[step =0.0843821,shift={(9.9359,1.42576)}](0,-0.253146) \boxpath;
\draw[step =0.0843821,shift={(9.9359,1.42576)}](0.0843821,-0.253146) \boxpath;
\draw[step =0.0843821,shift={(9.9359,1.42576)}](0.168764,-0.253146) \boxpath;
\draw[step =0.0843821,shift={(9.9359,1.42576)}](0.253146,-0.253146) \boxpath;
\draw[step =0.0843821,shift={(9.9359,1.42576)}](0.337528,-0.253146) \boxpath;
\draw[step =0.0843821,shift={(9.9359,1.42576)}](0,-0.337528) \boxpath;
\draw[step =0.0843821,shift={(9.9359,1.42576)}](0.0843821,-0.337528) \boxpath;
\draw[step =0.0843821,shift={(9.9359,1.42576)}](0.168764,-0.337528) \boxpath;
\draw[step =0.0843821,shift={(9.9359,1.42576)}](0,-0.421911) \boxpath;
[1,1,0]
\filldraw [draw=gray,thick,fill opacity = .5, fill = yellow,shift={(1.000000,1.73205)}](0,0) -- (60:2cm) -- (0:2cm) --(0,0) node[blue] at (0.8cm,0.4cm) {} node[black,fill opacity=1, draw opacity = 1] at (1cm,0.7cm) {} node[blue] at (0.8cm,1.2cm) {};
\def\boxpath{-- +(-0.154701,0) -- +(-0.154701,-0.154701) -- +(0,-0.154701) -- cycle};
\draw[step =0.154701,shift={(1.8359,2.86025)}](0,-0.154701) \boxpath;
\draw[step =0.154701,shift={(1.8359,2.86025)}](0,-0.309401) \boxpath;
[2,1,0]
\filldraw [draw=gray, thick,fill = pink!30, shift={(3,1.73205)}](0,0) -- (60:2cm) -- (120:2cm) --(0,0) node[green,fill opacity=1] at (0cm,1.6cm) {} node[black,fill opacity=1, draw opacity = 1] at (0cm,1.2cm) {} node[green,fill opacity=1] at (0cm,0.8cm) {};
\def\boxpath{-- +(-0.154701,0) -- +(-0.154701,-0.154701) -- +(0,-0.154701) -- cycle};
\draw[step =0.154701,shift={(2.9359,3.2)}](0,-0.154701) \boxpath;
\draw[step =0.154701,shift={(2.9359,3.2)}](0.154701,-0.154701) \boxpath;
\draw[step =0.154701,shift={(2.9359,3.2)}](0,-0.309401) \boxpath;
\draw[step =0.154701,shift={(2.9359,3.2)}](0,-0.464102) \boxpath;
[2,1,1]
\filldraw [draw=gray,thick,fill opacity = .5, fill = yellow,shift={(3.000000,1.73205)}](0,0) -- (60:2cm) -- (0:2cm) --(0,0) node[blue] at (0.8cm,0.4cm) {} node[black,fill opacity=1, draw opacity = 1] at (1cm,0.7cm) {} node[blue] at (0.8cm,1.2cm) {};
\def\boxpath{-- +(-0.154701,0) -- +(-0.154701,-0.154701) -- +(0,-0.154701) -- cycle};
\draw[step =0.154701,shift={(3.8359,2.86025)}](0,-0.154701) \boxpath;
\draw[step =0.154701,shift={(3.8359,2.86025)}](0.154701,-0.154701) \boxpath;
\draw[step =0.154701,shift={(3.8359,2.86025)}](0.309401,-0.154701) \boxpath;
\draw[step =0.154701,shift={(3.8359,2.86025)}](0,-0.309401) \boxpath;
\draw[step =0.154701,shift={(3.8359,2.86025)}](0,-0.464102) \boxpath;
[3,1,1]
\filldraw [draw=gray, thick,fill = white, shift={(5,1.73205)}](0,0) -- (60:2cm) -- (120:2cm) --(0,0) node[green,fill opacity=1] at (0cm,1.6cm) {} node[black,fill opacity=1, draw opacity = 1] at (0cm,1.2cm) {} node[green,fill opacity=1] at (0cm,0.8cm) {};
\def\boxpath{-- +(-0.154701,0) -- +(-0.154701,-0.154701) -- +(0,-0.154701) -- cycle};
\draw[step =0.154701,shift={(4.9359,3.3547)}](0,-0.154701) \boxpath;
\draw[step =0.154701,shift={(4.9359,3.3547)}](0.154701,-0.154701) \boxpath;
\draw[step =0.154701,shift={(4.9359,3.3547)}](0.309401,-0.154701) \boxpath;
\draw[step =0.154701,shift={(4.9359,3.3547)}](0.464102,-0.154701) \boxpath;
\draw[step =0.154701,shift={(4.9359,3.3547)}](0,-0.309401) \boxpath;
\draw[step =0.154701,shift={(4.9359,3.3547)}](0.154701,-0.309401) \boxpath;
\draw[step =0.154701,shift={(4.9359,3.3547)}](0,-0.464102) \boxpath;
\draw[step =0.154701,shift={(4.9359,3.3547)}](0,-0.618802) \boxpath;
[3,1,2]
\filldraw [draw=gray,thick,fill opacity = .5, fill = pink!30,shift={(5.000000,1.73205)}](0,0) -- (60:2cm) -- (0:2cm) --(0,0) node[blue] at (0.8cm,0.4cm) {} node[black,fill opacity=1, draw opacity = 1] at (1cm,0.7cm) {} node[blue] at (0.8cm,1.2cm) {};
\def\boxpath{-- +(-0.1326,0) -- +(-0.1326,-0.1326) -- +(0,-0.1326) -- cycle};
\draw[step =0.1326,shift={(5.8359,2.86025)}](0,-0.1326) \boxpath;
\draw[step =0.1326,shift={(5.8359,2.86025)}](0.1326,-0.1326) \boxpath;
\draw[step =0.1326,shift={(5.8359,2.86025)}](0.265201,-0.1326) \boxpath;
\draw[step =0.1326,shift={(5.8359,2.86025)}](0.397801,-0.1326) \boxpath;
\draw[step =0.1326,shift={(5.8359,2.86025)}](0.530402,-0.1326) \boxpath;
\draw[step =0.1326,shift={(5.8359,2.86025)}](0,-0.265201) \boxpath;
\draw[step =0.1326,shift={(5.8359,2.86025)}](0.1326,-0.265201) \boxpath;
\draw[step =0.1326,shift={(5.8359,2.86025)}](0.265201,-0.265201) \boxpath;
\draw[step =0.1326,shift={(5.8359,2.86025)}](0,-0.397801) \boxpath;
\draw[step =0.1326,shift={(5.8359,2.86025)}](0,-0.530402) \boxpath;
[4,1,2]
\filldraw [draw=gray, thick,fill = white, shift={(7,1.73205)}](0,0) -- (60:2cm) -- (120:2cm) --(0,0) node[green,fill opacity=1] at (0cm,1.6cm) {} node[black,fill opacity=1, draw opacity = 1] at (0cm,1.2cm) {} node[green,fill opacity=1] at (0cm,0.8cm) {};
\def\boxpath{-- +(-0.116025,0) -- +(-0.116025,-0.116025) -- +(0,-0.116025) -- cycle};
\draw[step =0.116025,shift={(6.9359,3.31603)}](0,-0.116025) \boxpath;
\draw[step =0.116025,shift={(6.9359,3.31603)}](0.116025,-0.116025) \boxpath;
\draw[step =0.116025,shift={(6.9359,3.31603)}](0.232051,-0.116025) \boxpath;
\draw[step =0.116025,shift={(6.9359,3.31603)}](0.348076,-0.116025) \boxpath;
\draw[step =0.116025,shift={(6.9359,3.31603)}](0.464102,-0.116025) \boxpath;
\draw[step =0.116025,shift={(6.9359,3.31603)}](0.580127,-0.116025) \boxpath;
\draw[step =0.116025,shift={(6.9359,3.31603)}](0,-0.232051) \boxpath;
\draw[step =0.116025,shift={(6.9359,3.31603)}](0.116025,-0.232051) \boxpath;
\draw[step =0.116025,shift={(6.9359,3.31603)}](0.232051,-0.232051) \boxpath;
\draw[step =0.116025,shift={(6.9359,3.31603)}](0.348076,-0.232051) \boxpath;
\draw[step =0.116025,shift={(6.9359,3.31603)}](0,-0.348076) \boxpath;
\draw[step =0.116025,shift={(6.9359,3.31603)}](0.116025,-0.348076) \boxpath;
\draw[step =0.116025,shift={(6.9359,3.31603)}](0,-0.464102) \boxpath;
\draw[step =0.116025,shift={(6.9359,3.31603)}](0,-0.580127) \boxpath;
[4,1,3]
\filldraw [draw=gray,thick,fill opacity = .5, fill = white,shift={(7.000000,1.73205)}](0,0) -- (60:2cm) -- (0:2cm) --(0,0) node[blue] at (0.8cm,0.4cm) {} node[black,fill opacity=1, draw opacity = 1] at (1cm,0.7cm) {} node[blue] at (0.8cm,1.2cm) {};
\def\boxpath{-- +(-0.103134,0) -- +(-0.103134,-0.103134) -- +(0,-0.103134) -- cycle};
\draw[step =0.103134,shift={(7.8359,2.86025)}](0,-0.103134) \boxpath;
\draw[step =0.103134,shift={(7.8359,2.86025)}](0.103134,-0.103134) \boxpath;
\draw[step =0.103134,shift={(7.8359,2.86025)}](0.206267,-0.103134) \boxpath;
\draw[step =0.103134,shift={(7.8359,2.86025)}](0.309401,-0.103134) \boxpath;
\draw[step =0.103134,shift={(7.8359,2.86025)}](0.412535,-0.103134) \boxpath;
\draw[step =0.103134,shift={(7.8359,2.86025)}](0.515668,-0.103134) \boxpath;
\draw[step =0.103134,shift={(7.8359,2.86025)}](0.618802,-0.103134) \boxpath;
\draw[step =0.103134,shift={(7.8359,2.86025)}](0,-0.206267) \boxpath;
\draw[step =0.103134,shift={(7.8359,2.86025)}](0.103134,-0.206267) \boxpath;
\draw[step =0.103134,shift={(7.8359,2.86025)}](0.206267,-0.206267) \boxpath;
\draw[step =0.103134,shift={(7.8359,2.86025)}](0.309401,-0.206267) \boxpath;
\draw[step =0.103134,shift={(7.8359,2.86025)}](0.412535,-0.206267) \boxpath;
\draw[step =0.103134,shift={(7.8359,2.86025)}](0,-0.309401) \boxpath;
\draw[step =0.103134,shift={(7.8359,2.86025)}](0.103134,-0.309401) \boxpath;
\draw[step =0.103134,shift={(7.8359,2.86025)}](0.206267,-0.309401) \boxpath;
\draw[step =0.103134,shift={(7.8359,2.86025)}](0,-0.412535) \boxpath;
\draw[step =0.103134,shift={(7.8359,2.86025)}](0,-0.515668) \boxpath;
[5,1,3]
\filldraw [draw=gray, thick,fill = white, shift={(9,1.73205)}](0,0) -- (60:2cm) -- (120:2cm) --(0,0) node[green,fill opacity=1] at (0cm,1.6cm) {} node[black,fill opacity=1, draw opacity = 1] at (0cm,1.2cm) {} node[green,fill opacity=1] at (0cm,0.8cm) {};
\def\boxpath{-- +(-0.0928203,0) -- +(-0.0928203,-0.0928203) -- +(0,-0.0928203) -- cycle};
\draw[step =0.0928203,shift={(8.9359,3.29282)}](0,-0.0928203) \boxpath;
\draw[step =0.0928203,shift={(8.9359,3.29282)}](0.0928203,-0.0928203) \boxpath;
\draw[step =0.0928203,shift={(8.9359,3.29282)}](0.185641,-0.0928203) \boxpath;
\draw[step =0.0928203,shift={(8.9359,3.29282)}](0.278461,-0.0928203) \boxpath;
\draw[step =0.0928203,shift={(8.9359,3.29282)}](0.371281,-0.0928203) \boxpath;
\draw[step =0.0928203,shift={(8.9359,3.29282)}](0.464102,-0.0928203) \boxpath;
\draw[step =0.0928203,shift={(8.9359,3.29282)}](0.556922,-0.0928203) \boxpath;
\draw[step =0.0928203,shift={(8.9359,3.29282)}](0.649742,-0.0928203) \boxpath;
\draw[step =0.0928203,shift={(8.9359,3.29282)}](0,-0.185641) \boxpath;
\draw[step =0.0928203,shift={(8.9359,3.29282)}](0.0928203,-0.185641) \boxpath;
\draw[step =0.0928203,shift={(8.9359,3.29282)}](0.185641,-0.185641) \boxpath;
\draw[step =0.0928203,shift={(8.9359,3.29282)}](0.278461,-0.185641) \boxpath;
\draw[step =0.0928203,shift={(8.9359,3.29282)}](0.371281,-0.185641) \boxpath;
\draw[step =0.0928203,shift={(8.9359,3.29282)}](0.464102,-0.185641) \boxpath;
\draw[step =0.0928203,shift={(8.9359,3.29282)}](0,-0.278461) \boxpath;
\draw[step =0.0928203,shift={(8.9359,3.29282)}](0.0928203,-0.278461) \boxpath;
\draw[step =0.0928203,shift={(8.9359,3.29282)}](0.185641,-0.278461) \boxpath;
\draw[step =0.0928203,shift={(8.9359,3.29282)}](0.278461,-0.278461) \boxpath;
\draw[step =0.0928203,shift={(8.9359,3.29282)}](0,-0.371281) \boxpath;
\draw[step =0.0928203,shift={(8.9359,3.29282)}](0.0928203,-0.371281) \boxpath;
\draw[step =0.0928203,shift={(8.9359,3.29282)}](0,-0.464102) \boxpath;
\draw[step =0.0928203,shift={(8.9359,3.29282)}](0,-0.556922) \boxpath;
[5,1,4]
\filldraw [draw=gray,thick,fill opacity = .5, fill = white,shift={(9.000000,1.73205)}](0,0) -- (60:2cm) -- (0:2cm) --(0,0) node[blue] at (0.8cm,0.4cm) {} node[black,fill opacity=1, draw opacity = 1] at (1cm,0.7cm) {} node[blue] at (0.8cm,1.2cm) {};
\def\boxpath{-- +(-0.0843821,0) -- +(-0.0843821,-0.0843821) -- +(0,-0.0843821) -- cycle};
\draw[step =0.0843821,shift={(9.8359,2.86025)}](0,-0.0843821) \boxpath;
\draw[step =0.0843821,shift={(9.8359,2.86025)}](0.0843821,-0.0843821) \boxpath;
\draw[step =0.0843821,shift={(9.8359,2.86025)}](0.168764,-0.0843821) \boxpath;
\draw[step =0.0843821,shift={(9.8359,2.86025)}](0.253146,-0.0843821) \boxpath;
\draw[step =0.0843821,shift={(9.8359,2.86025)}](0.337528,-0.0843821) \boxpath;
\draw[step =0.0843821,shift={(9.8359,2.86025)}](0.421911,-0.0843821) \boxpath;
\draw[step =0.0843821,shift={(9.8359,2.86025)}](0.506293,-0.0843821) \boxpath;
\draw[step =0.0843821,shift={(9.8359,2.86025)}](0.590675,-0.0843821) \boxpath;
\draw[step =0.0843821,shift={(9.8359,2.86025)}](0.675057,-0.0843821) \boxpath;
\draw[step =0.0843821,shift={(9.8359,2.86025)}](0,-0.168764) \boxpath;
\draw[step =0.0843821,shift={(9.8359,2.86025)}](0.0843821,-0.168764) \boxpath;
\draw[step =0.0843821,shift={(9.8359,2.86025)}](0.168764,-0.168764) \boxpath;
\draw[step =0.0843821,shift={(9.8359,2.86025)}](0.253146,-0.168764) \boxpath;
\draw[step =0.0843821,shift={(9.8359,2.86025)}](0.337528,-0.168764) \boxpath;
\draw[step =0.0843821,shift={(9.8359,2.86025)}](0.421911,-0.168764) \boxpath;
\draw[step =0.0843821,shift={(9.8359,2.86025)}](0.506293,-0.168764) \boxpath;
\draw[step =0.0843821,shift={(9.8359,2.86025)}](0,-0.253146) \boxpath;
\draw[step =0.0843821,shift={(9.8359,2.86025)}](0.0843821,-0.253146) \boxpath;
\draw[step =0.0843821,shift={(9.8359,2.86025)}](0.168764,-0.253146) \boxpath;
\draw[step =0.0843821,shift={(9.8359,2.86025)}](0.253146,-0.253146) \boxpath;
\draw[step =0.0843821,shift={(9.8359,2.86025)}](0.337528,-0.253146) \boxpath;
\draw[step =0.0843821,shift={(9.8359,2.86025)}](0,-0.337528) \boxpath;
\draw[step =0.0843821,shift={(9.8359,2.86025)}](0.0843821,-0.337528) \boxpath;
\draw[step =0.0843821,shift={(9.8359,2.86025)}](0.168764,-0.337528) \boxpath;
\draw[step =0.0843821,shift={(9.8359,2.86025)}](0,-0.421911) \boxpath;
\draw[step =0.0843821,shift={(9.8359,2.86025)}](0,-0.506293) \boxpath;
[6,1,4]
\filldraw [draw=gray, thick,fill = white, shift={(11,1.73205)}](0,0) -- (60:2cm) -- (120:2cm) --(0,0) node[green,fill opacity=1] at (0cm,1.6cm) {} node[black,fill opacity=1, draw opacity = 1] at (0cm,1.2cm) {} node[green,fill opacity=1] at (0cm,0.8cm) {};
\def\boxpath{-- +(-0.0773503,0) -- +(-0.0773503,-0.0773503) -- +(0,-0.0773503) -- cycle};
\draw[step =0.0773503,shift={(10.9359,3.27735)}](0,-0.0773503) \boxpath;
\draw[step =0.0773503,shift={(10.9359,3.27735)}](0.0773503,-0.0773503) \boxpath;
\draw[step =0.0773503,shift={(10.9359,3.27735)}](0.154701,-0.0773503) \boxpath;
\draw[step =0.0773503,shift={(10.9359,3.27735)}](0.232051,-0.0773503) \boxpath;
\draw[step =0.0773503,shift={(10.9359,3.27735)}](0.309401,-0.0773503) \boxpath;
\draw[step =0.0773503,shift={(10.9359,3.27735)}](0.386751,-0.0773503) \boxpath;
\draw[step =0.0773503,shift={(10.9359,3.27735)}](0.464102,-0.0773503) \boxpath;
\draw[step =0.0773503,shift={(10.9359,3.27735)}](0.541452,-0.0773503) \boxpath;
\draw[step =0.0773503,shift={(10.9359,3.27735)}](0.618802,-0.0773503) \boxpath;
\draw[step =0.0773503,shift={(10.9359,3.27735)}](0.696152,-0.0773503) \boxpath;
\draw[step =0.0773503,shift={(10.9359,3.27735)}](0,-0.154701) \boxpath;
\draw[step =0.0773503,shift={(10.9359,3.27735)}](0.0773503,-0.154701) \boxpath;
\draw[step =0.0773503,shift={(10.9359,3.27735)}](0.154701,-0.154701) \boxpath;
\draw[step =0.0773503,shift={(10.9359,3.27735)}](0.232051,-0.154701) \boxpath;
\draw[step =0.0773503,shift={(10.9359,3.27735)}](0.309401,-0.154701) \boxpath;
\draw[step =0.0773503,shift={(10.9359,3.27735)}](0.386751,-0.154701) \boxpath;
\draw[step =0.0773503,shift={(10.9359,3.27735)}](0.464102,-0.154701) \boxpath;
\draw[step =0.0773503,shift={(10.9359,3.27735)}](0.541452,-0.154701) \boxpath;
\draw[step =0.0773503,shift={(10.9359,3.27735)}](0,-0.232051) \boxpath;
\draw[step =0.0773503,shift={(10.9359,3.27735)}](0.0773503,-0.232051) \boxpath;
\draw[step =0.0773503,shift={(10.9359,3.27735)}](0.154701,-0.232051) \boxpath;
\draw[step =0.0773503,shift={(10.9359,3.27735)}](0.232051,-0.232051) \boxpath;
\draw[step =0.0773503,shift={(10.9359,3.27735)}](0.309401,-0.232051) \boxpath;
\draw[step =0.0773503,shift={(10.9359,3.27735)}](0.386751,-0.232051) \boxpath;
\draw[step =0.0773503,shift={(10.9359,3.27735)}](0,-0.309401) \boxpath;
\draw[step =0.0773503,shift={(10.9359,3.27735)}](0.0773503,-0.309401) \boxpath;
\draw[step =0.0773503,shift={(10.9359,3.27735)}](0.154701,-0.309401) \boxpath;
\draw[step =0.0773503,shift={(10.9359,3.27735)}](0.232051,-0.309401) \boxpath;
\draw[step =0.0773503,shift={(10.9359,3.27735)}](0,-0.386751) \boxpath;
\draw[step =0.0773503,shift={(10.9359,3.27735)}](0.0773503,-0.386751) \boxpath;
\draw[step =0.0773503,shift={(10.9359,3.27735)}](0,-0.464102) \boxpath;
\draw[step =0.0773503,shift={(10.9359,3.27735)}](0,-0.541452) \boxpath;
[2,2,0]
\filldraw [draw=gray,thick,fill opacity = .5, fill = pink!30,shift={(2.000000,3.4641)}](0,0) -- (60:2cm) -- (0:2cm) --(0,0) node[blue] at (0.8cm,0.4cm) {} node[black,fill opacity=1, draw opacity = 1] at (1cm,0.7cm) {} node[blue] at (0.8cm,1.2cm) {};
\def\boxpath{-- +(-0.154701,0) -- +(-0.154701,-0.154701) -- +(0,-0.154701) -- cycle};
\draw[step =0.154701,shift={(2.8359,4.5923)}](0,-0.154701) \boxpath;
\draw[step =0.154701,shift={(2.8359,4.5923)}](0.154701,-0.154701) \boxpath;
\draw[step =0.154701,shift={(2.8359,4.5923)}](0,-0.309401) \boxpath;
\draw[step =0.154701,shift={(2.8359,4.5923)}](0.154701,-0.309401) \boxpath;
\draw[step =0.154701,shift={(2.8359,4.5923)}](0,-0.464102) \boxpath;
\draw[step =0.154701,shift={(2.8359,4.5923)}](0,-0.618802) \boxpath;
[3,2,0]
\filldraw [draw=gray, thick,fill = white, shift={(4,3.4641)}](0,0) -- (60:2cm) -- (120:2cm) --(0,0) node[green,fill opacity=1] at (0cm,1.6cm) {} node[black,fill opacity=1, draw opacity = 1] at (0cm,1.2cm) {} node[green,fill opacity=1] at (0cm,0.8cm) {};
\def\boxpath{-- +(-0.154701,0) -- +(-0.154701,-0.154701) -- +(0,-0.154701) -- cycle};
\draw[step =0.154701,shift={(3.9359,5.24145)}](0,-0.154701) \boxpath;
\draw[step =0.154701,shift={(3.9359,5.24145)}](0.154701,-0.154701) \boxpath;
\draw[step =0.154701,shift={(3.9359,5.24145)}](0.309401,-0.154701) \boxpath;
\draw[step =0.154701,shift={(3.9359,5.24145)}](0,-0.309401) \boxpath;
\draw[step =0.154701,shift={(3.9359,5.24145)}](0.154701,-0.309401) \boxpath;
\draw[step =0.154701,shift={(3.9359,5.24145)}](0,-0.464102) \boxpath;
\draw[step =0.154701,shift={(3.9359,5.24145)}](0.154701,-0.464102) \boxpath;
\draw[step =0.154701,shift={(3.9359,5.24145)}](0,-0.618802) \boxpath;
\draw[step =0.154701,shift={(3.9359,5.24145)}](0,-0.773503) \boxpath;
[3,2,1]
\filldraw [draw=gray,thick,fill opacity = .5, fill = pink!30,shift={(4.000000,3.4641)}](0,0) -- (60:2cm) -- (0:2cm) --(0,0) node[blue] at (0.8cm,0.4cm) {} node[black,fill opacity=1, draw opacity = 1] at (1cm,0.7cm) {} node[blue] at (0.8cm,1.2cm) {};
\def\boxpath{-- +(-0.154701,0) -- +(-0.154701,-0.154701) -- +(0,-0.154701) -- cycle};
\draw[step =0.154701,shift={(4.8359,4.5923)}](0,-0.154701) \boxpath;
\draw[step =0.154701,shift={(4.8359,4.5923)}](0.154701,-0.154701) \boxpath;
\draw[step =0.154701,shift={(4.8359,4.5923)}](0.309401,-0.154701) \boxpath;
\draw[step =0.154701,shift={(4.8359,4.5923)}](0.464102,-0.154701) \boxpath;
\draw[step =0.154701,shift={(4.8359,4.5923)}](0,-0.309401) \boxpath;
\draw[step =0.154701,shift={(4.8359,4.5923)}](0.154701,-0.309401) \boxpath;
\draw[step =0.154701,shift={(4.8359,4.5923)}](0,-0.464102) \boxpath;
\draw[step =0.154701,shift={(4.8359,4.5923)}](0.154701,-0.464102) \boxpath;
\draw[step =0.154701,shift={(4.8359,4.5923)}](0,-0.618802) \boxpath;
\draw[step =0.154701,shift={(4.8359,4.5923)}](0,-0.773503) \boxpath;
[4,2,1]
\filldraw [draw=gray, thick,fill = white, shift={(6,3.4641)}](0,0) -- (60:2cm) -- (120:2cm) --(0,0) node[green,fill opacity=1] at (0cm,1.6cm) {} node[black,fill opacity=1, draw opacity = 1] at (0cm,1.2cm) {} node[green,fill opacity=1] at (0cm,0.8cm) {};
\def\boxpath{-- +(-0.116025,0) -- +(-0.116025,-0.116025) -- +(0,-0.116025) -- cycle};
\draw[step =0.116025,shift={(5.9359,5.1641)}](0,-0.116025) \boxpath;
\draw[step =0.116025,shift={(5.9359,5.1641)}](0.116025,-0.116025) \boxpath;
\draw[step =0.116025,shift={(5.9359,5.1641)}](0.232051,-0.116025) \boxpath;
\draw[step =0.116025,shift={(5.9359,5.1641)}](0.348076,-0.116025) \boxpath;
\draw[step =0.116025,shift={(5.9359,5.1641)}](0.464102,-0.116025) \boxpath;
\draw[step =0.116025,shift={(5.9359,5.1641)}](0,-0.232051) \boxpath;
\draw[step =0.116025,shift={(5.9359,5.1641)}](0.116025,-0.232051) \boxpath;
\draw[step =0.116025,shift={(5.9359,5.1641)}](0.232051,-0.232051) \boxpath;
\draw[step =0.116025,shift={(5.9359,5.1641)}](0,-0.348076) \boxpath;
\draw[step =0.116025,shift={(5.9359,5.1641)}](0.116025,-0.348076) \boxpath;
\draw[step =0.116025,shift={(5.9359,5.1641)}](0,-0.464102) \boxpath;
\draw[step =0.116025,shift={(5.9359,5.1641)}](0.116025,-0.464102) \boxpath;
\draw[step =0.116025,shift={(5.9359,5.1641)}](0,-0.580127) \boxpath;
\draw[step =0.116025,shift={(5.9359,5.1641)}](0,-0.696152) \boxpath;
[4,2,2]
\filldraw [draw=gray,thick,fill opacity = .5, fill = pink!30,shift={(6.000000,3.4641)}](0,0) -- (60:2cm) -- (0:2cm) --(0,0) node[blue] at (0.8cm,0.4cm) {} node[black,fill opacity=1, draw opacity = 1] at (1cm,0.7cm) {} node[blue] at (0.8cm,1.2cm) {};
\def\boxpath{-- +(-0.116025,0) -- +(-0.116025,-0.116025) -- +(0,-0.116025) -- cycle};
\draw[step =0.116025,shift={(6.8359,4.5923)}](0,-0.116025) \boxpath;
\draw[step =0.116025,shift={(6.8359,4.5923)}](0.116025,-0.116025) \boxpath;
\draw[step =0.116025,shift={(6.8359,4.5923)}](0.232051,-0.116025) \boxpath;
\draw[step =0.116025,shift={(6.8359,4.5923)}](0.348076,-0.116025) \boxpath;
\draw[step =0.116025,shift={(6.8359,4.5923)}](0.464102,-0.116025) \boxpath;
\draw[step =0.116025,shift={(6.8359,4.5923)}](0.580127,-0.116025) \boxpath;
\draw[step =0.116025,shift={(6.8359,4.5923)}](0,-0.232051) \boxpath;
\draw[step =0.116025,shift={(6.8359,4.5923)}](0.116025,-0.232051) \boxpath;
\draw[step =0.116025,shift={(6.8359,4.5923)}](0.232051,-0.232051) \boxpath;
\draw[step =0.116025,shift={(6.8359,4.5923)}](0.348076,-0.232051) \boxpath;
\draw[step =0.116025,shift={(6.8359,4.5923)}](0,-0.348076) \boxpath;
\draw[step =0.116025,shift={(6.8359,4.5923)}](0.116025,-0.348076) \boxpath;
\draw[step =0.116025,shift={(6.8359,4.5923)}](0,-0.464102) \boxpath;
\draw[step =0.116025,shift={(6.8359,4.5923)}](0.116025,-0.464102) \boxpath;
\draw[step =0.116025,shift={(6.8359,4.5923)}](0,-0.580127) \boxpath;
\draw[step =0.116025,shift={(6.8359,4.5923)}](0,-0.696152) \boxpath;
[5,2,2]
\filldraw [draw=gray, thick,fill = white, shift={(8,3.4641)}](0,0) -- (60:2cm) -- (120:2cm) --(0,0) node[green,fill opacity=1] at (0cm,1.6cm) {} node[black,fill opacity=1, draw opacity = 1] at (0cm,1.2cm) {} node[green,fill opacity=1] at (0cm,0.8cm) {};
\def\boxpath{-- +(-0.103134,0) -- +(-0.103134,-0.103134) -- +(0,-0.103134) -- cycle};
\draw[step =0.103134,shift={(7.9359,5.18989)}](0,-0.103134) \boxpath;
\draw[step =0.103134,shift={(7.9359,5.18989)}](0.103134,-0.103134) \boxpath;
\draw[step =0.103134,shift={(7.9359,5.18989)}](0.206267,-0.103134) \boxpath;
\draw[step =0.103134,shift={(7.9359,5.18989)}](0.309401,-0.103134) \boxpath;
\draw[step =0.103134,shift={(7.9359,5.18989)}](0.412535,-0.103134) \boxpath;
\draw[step =0.103134,shift={(7.9359,5.18989)}](0.515668,-0.103134) \boxpath;
\draw[step =0.103134,shift={(7.9359,5.18989)}](0.618802,-0.103134) \boxpath;
\draw[step =0.103134,shift={(7.9359,5.18989)}](0,-0.206267) \boxpath;
\draw[step =0.103134,shift={(7.9359,5.18989)}](0.103134,-0.206267) \boxpath;
\draw[step =0.103134,shift={(7.9359,5.18989)}](0.206267,-0.206267) \boxpath;
\draw[step =0.103134,shift={(7.9359,5.18989)}](0.309401,-0.206267) \boxpath;
\draw[step =0.103134,shift={(7.9359,5.18989)}](0.412535,-0.206267) \boxpath;
\draw[step =0.103134,shift={(7.9359,5.18989)}](0,-0.309401) \boxpath;
\draw[step =0.103134,shift={(7.9359,5.18989)}](0.103134,-0.309401) \boxpath;
\draw[step =0.103134,shift={(7.9359,5.18989)}](0.206267,-0.309401) \boxpath;
\draw[step =0.103134,shift={(7.9359,5.18989)}](0,-0.412535) \boxpath;
\draw[step =0.103134,shift={(7.9359,5.18989)}](0.103134,-0.412535) \boxpath;
\draw[step =0.103134,shift={(7.9359,5.18989)}](0,-0.515668) \boxpath;
\draw[step =0.103134,shift={(7.9359,5.18989)}](0.103134,-0.515668) \boxpath;
\draw[step =0.103134,shift={(7.9359,5.18989)}](0,-0.618802) \boxpath;
\draw[step =0.103134,shift={(7.9359,5.18989)}](0,-0.721936) \boxpath;
[5,2,3]
\filldraw [draw=gray,thick,fill opacity = .5, fill = white,shift={(8.000000,3.4641)}](0,0) -- (60:2cm) -- (0:2cm) --(0,0) node[blue] at (0.8cm,0.4cm) {} node[black,fill opacity=1, draw opacity = 1] at (1cm,0.7cm) {} node[blue] at (0.8cm,1.2cm) {};
\def\boxpath{-- +(-0.0928203,0) -- +(-0.0928203,-0.0928203) -- +(0,-0.0928203) -- cycle};
\draw[step =0.0928203,shift={(8.8359,4.5923)}](0,-0.0928203) \boxpath;
\draw[step =0.0928203,shift={(8.8359,4.5923)}](0.0928203,-0.0928203) \boxpath;
\draw[step =0.0928203,shift={(8.8359,4.5923)}](0.185641,-0.0928203) \boxpath;
\draw[step =0.0928203,shift={(8.8359,4.5923)}](0.278461,-0.0928203) \boxpath;
\draw[step =0.0928203,shift={(8.8359,4.5923)}](0.371281,-0.0928203) \boxpath;
\draw[step =0.0928203,shift={(8.8359,4.5923)}](0.464102,-0.0928203) \boxpath;
\draw[step =0.0928203,shift={(8.8359,4.5923)}](0.556922,-0.0928203) \boxpath;
\draw[step =0.0928203,shift={(8.8359,4.5923)}](0.649742,-0.0928203) \boxpath;
\draw[step =0.0928203,shift={(8.8359,4.5923)}](0,-0.185641) \boxpath;
\draw[step =0.0928203,shift={(8.8359,4.5923)}](0.0928203,-0.185641) \boxpath;
\draw[step =0.0928203,shift={(8.8359,4.5923)}](0.185641,-0.185641) \boxpath;
\draw[step =0.0928203,shift={(8.8359,4.5923)}](0.278461,-0.185641) \boxpath;
\draw[step =0.0928203,shift={(8.8359,4.5923)}](0.371281,-0.185641) \boxpath;
\draw[step =0.0928203,shift={(8.8359,4.5923)}](0.464102,-0.185641) \boxpath;
\draw[step =0.0928203,shift={(8.8359,4.5923)}](0,-0.278461) \boxpath;
\draw[step =0.0928203,shift={(8.8359,4.5923)}](0.0928203,-0.278461) \boxpath;
\draw[step =0.0928203,shift={(8.8359,4.5923)}](0.185641,-0.278461) \boxpath;
\draw[step =0.0928203,shift={(8.8359,4.5923)}](0.278461,-0.278461) \boxpath;
\draw[step =0.0928203,shift={(8.8359,4.5923)}](0,-0.371281) \boxpath;
\draw[step =0.0928203,shift={(8.8359,4.5923)}](0.0928203,-0.371281) \boxpath;
\draw[step =0.0928203,shift={(8.8359,4.5923)}](0,-0.464102) \boxpath;
\draw[step =0.0928203,shift={(8.8359,4.5923)}](0.0928203,-0.464102) \boxpath;
\draw[step =0.0928203,shift={(8.8359,4.5923)}](0,-0.556922) \boxpath;
\draw[step =0.0928203,shift={(8.8359,4.5923)}](0,-0.649742) \boxpath;
[6,2,3]
\filldraw [draw=gray, thick,fill = white, shift={(10,3.4641)}](0,0) -- (60:2cm) -- (120:2cm) --(0,0) node[green,fill opacity=1] at (0cm,1.6cm) {} node[black,fill opacity=1, draw opacity = 1] at (0cm,1.2cm) {} node[green,fill opacity=1] at (0cm,0.8cm) {};
\def\boxpath{-- +(-0.0843821,0) -- +(-0.0843821,-0.0843821) -- +(0,-0.0843821) -- cycle};
\draw[step =0.0843821,shift={(9.9359,5.14301)}](0,-0.0843821) \boxpath;
\draw[step =0.0843821,shift={(9.9359,5.14301)}](0.0843821,-0.0843821) \boxpath;
\draw[step =0.0843821,shift={(9.9359,5.14301)}](0.168764,-0.0843821) \boxpath;
\draw[step =0.0843821,shift={(9.9359,5.14301)}](0.253146,-0.0843821) \boxpath;
\draw[step =0.0843821,shift={(9.9359,5.14301)}](0.337528,-0.0843821) \boxpath;
\draw[step =0.0843821,shift={(9.9359,5.14301)}](0.421911,-0.0843821) \boxpath;
\draw[step =0.0843821,shift={(9.9359,5.14301)}](0.506293,-0.0843821) \boxpath;
\draw[step =0.0843821,shift={(9.9359,5.14301)}](0.590675,-0.0843821) \boxpath;
\draw[step =0.0843821,shift={(9.9359,5.14301)}](0.675057,-0.0843821) \boxpath;
\draw[step =0.0843821,shift={(9.9359,5.14301)}](0,-0.168764) \boxpath;
\draw[step =0.0843821,shift={(9.9359,5.14301)}](0.0843821,-0.168764) \boxpath;
\draw[step =0.0843821,shift={(9.9359,5.14301)}](0.168764,-0.168764) \boxpath;
\draw[step =0.0843821,shift={(9.9359,5.14301)}](0.253146,-0.168764) \boxpath;
\draw[step =0.0843821,shift={(9.9359,5.14301)}](0.337528,-0.168764) \boxpath;
\draw[step =0.0843821,shift={(9.9359,5.14301)}](0.421911,-0.168764) \boxpath;
\draw[step =0.0843821,shift={(9.9359,5.14301)}](0.506293,-0.168764) \boxpath;
\draw[step =0.0843821,shift={(9.9359,5.14301)}](0,-0.253146) \boxpath;
\draw[step =0.0843821,shift={(9.9359,5.14301)}](0.0843821,-0.253146) \boxpath;
\draw[step =0.0843821,shift={(9.9359,5.14301)}](0.168764,-0.253146) \boxpath;
\draw[step =0.0843821,shift={(9.9359,5.14301)}](0.253146,-0.253146) \boxpath;
\draw[step =0.0843821,shift={(9.9359,5.14301)}](0.337528,-0.253146) \boxpath;
\draw[step =0.0843821,shift={(9.9359,5.14301)}](0,-0.337528) \boxpath;
\draw[step =0.0843821,shift={(9.9359,5.14301)}](0.0843821,-0.337528) \boxpath;
\draw[step =0.0843821,shift={(9.9359,5.14301)}](0.168764,-0.337528) \boxpath;
\draw[step =0.0843821,shift={(9.9359,5.14301)}](0,-0.421911) \boxpath;
\draw[step =0.0843821,shift={(9.9359,5.14301)}](0.0843821,-0.421911) \boxpath;
\draw[step =0.0843821,shift={(9.9359,5.14301)}](0,-0.506293) \boxpath;
\draw[step =0.0843821,shift={(9.9359,5.14301)}](0.0843821,-0.506293) \boxpath;
\draw[step =0.0843821,shift={(9.9359,5.14301)}](0,-0.590675) \boxpath;
\draw[step =0.0843821,shift={(9.9359,5.14301)}](0,-0.675057) \boxpath;
\path [draw = red, ultra thick, draw opacity = 1, shift = {(0,0)}] (0,0) -- +(60:6);
\path [draw = red,ultra thick, draw opacity = 1] (0,0) -- +(10,0);
\path [draw = red,ultra thick,draw opacity = 1] (0,0) -- (0,0);

\path [draw = red, ultra thick, draw opacity = 1, shift = {(2,0)}] (0,0) -- +(60:6);
\path [draw = red,ultra thick, draw opacity = 1] (1,1.73205) -- +(10,0);
\path [draw = red,ultra thick,draw opacity = 1] (1,1.73205) -- (2,0);

\path [draw = red, ultra thick, draw opacity = 1, shift = {(4,0)},dashed] (0,0) -- +(60:6);
\path [draw = red,ultra thick, draw opacity = 1,dashed] (2,3.4641) -- +(10,0);
\path [draw = red,ultra thick,draw opacity = 1,dashed] (2,3.4641) -- (4,0);

\end{scope}

\end{tikzpicture}

\caption{The dominant chamber for $A_2$ with $3$-cores.The thick solid lines are hyperplanes from
  the $1$-Shi arrangement, the dashed ones are added to create the
  $2$-Shi arrangement.}
\label{fig:alcoveWCores}
\end{figure}
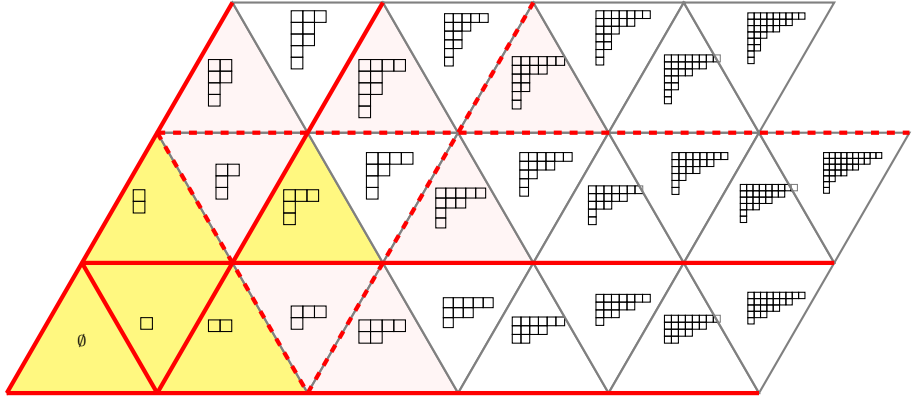

We'll just say a few more words about the geometry here. Richard's
pyramids have connected the core partitions to regions. We'll describe
Lascoux's \cite{lascoux1999} well-known bijection between $n$-cores
and certain elements of $\affSn$, and by extension, between $n$-cores
and alcoves in the dominant chamber. Please see Lapointe's and Morse's
paper \cite{lapointe-morse2005} for details. We describe the
bijection, as another way of seeing why core partitions pop up
here. An $n$-core partition may have several removable boxes of a
given residue or it may have several addable boxes of a given residue,
but it will never have both addable and removable boxes of the same
residue. Given an $n$-core partition $\lambda$ and the generators
$s_0,s_1,\ldots,s_{n-1}$ of $\affSn$, let $s_i(\lambda)$ be the
partition where all boxes of residue $i$ have been removed (added) if
there are removable (addable) boxes. Any $n$-core partition can be
expressed as $w(\emptyset)$. See Figure~\ref{fig:alcoveWCores} and
Example~\ref{ex:richards}. We associate the $n$-core $w(\emptyset)$
with the alcove $\Afund w$.

\omitt{We use Lapointe's and Morse's
\cite{lapointe-morse2005} reading of a word from the diagram of the
core partition. Given an $n$-core partition, fill the box with residue
$i$, $0\leq i\leq n-1$, with $s_i$. Now read the entries of the
diagram from bottome to top and from right to left and write the word
from left to right. We now have a reduced decomposition for some
$w\in\affSn$ which ends with $s_0$, ensuring that $w\Afund$ is in the
dominant chamber.}

We mention that Fishel and Vazirani mapped partitions
which are both $n$ and $nm+1$ cores to dominant regions in the $m$-Shi
arrangment of type $A_{n-1}$ in \cite{fishel-vazirani2010} using
abacus diagrams and the root lattice.

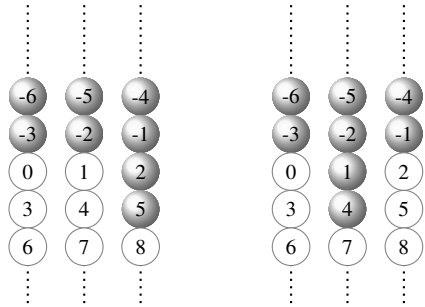
\begin{figure}[h]
\begin{tikzpicture}
\begin{scope}[fill opacity=1,every to/.style={draw,dotted},
 information text/.style={rounded corners,fill=red!10,inner sep=1ex}]
\tikzstyle{every node}=[font=\footnotesize]
\def\blankSpot{\filldraw[fill=white, draw=gray]}
\def\beadSpot{\shade[shading=ball, ball color = Gainsboro!30]}

\begin{scope}
\draw[thick,black,dotted,shift={(0,0)}] (0,1) to (0,5.0);
\draw[thick,black,dotted,shift={(0.75,0)}] (0,1) to (0,5.0);
\draw[thick,black,dotted,shift={(1.5,0)}] (0,1) to (0,5.0);

\beadSpot (0,3.75) circle (0.25) node {-6};
\beadSpot (0.75,3.75) circle (0.25) node {-5};
\beadSpot (1.5,3.75) circle (0.25) node {-4};

\beadSpot (0,3.25) circle (0.25) node {-3};
\beadSpot (0.75,3.25) circle (0.25) node {-2};
\beadSpot (1.5,3.25) circle (0.25) node {-1};

\blankSpot (0,2.75) circle (0.25) node {0};
\blankSpot (0.75,2.75) circle (0.25) node {1};
\beadSpot (1.5,2.75) circle (0.25) node {2};

\blankSpot (0,2.25) circle (0.25) node {3};
\blankSpot (0.75,2.25) circle (0.25) node {4};
\beadSpot (1.5,2.25) circle (0.25) node {5};

\blankSpot (0,1.75) circle (0.25) node {6};
\blankSpot (0.75,1.75) circle (0.25) node {7};
\blankSpot (1.5,1.75) circle (0.25) node {8};

\omitt{
\blankSpot (0,1.25) circle (0.25) node {15};
\blankSpot (0.75,1.25) circle (0.25) node {16};
\blankSpot (1.5,1.25) circle (0.25) node {17};
\blankSpot (2.25,1.25) circle (0.25) node {18};
\blankSpot (3,1.25) circle (0.25) node {19};
}

\end{scope}

\begin{scope}[shift={(3.5,0)}]
\draw[thick,black,dotted,shift={(0,0)}] (0,1) to (0,5.0);
\draw[thick,black,dotted,shift={(0.75,0)}] (0,1) to (0,5.0);
\draw[thick,black,dotted,shift={(1.5,0)}] (0,1) to (0,5.0);

\beadSpot (0,3.75) circle (0.25) node {-6};
\beadSpot (0.75,3.75) circle (0.25) node {-5};
\beadSpot (1.5,3.75) circle (0.25) node {-4};

\beadSpot (0,3.25) circle (0.25) node {-3};
\beadSpot (0.75,3.25) circle (0.25) node {-2};
\beadSpot (1.5,3.25) circle (0.25) node {-1};

\blankSpot (0,2.75) circle (0.25) node {0};
\beadSpot (0.75,2.75) circle (0.25) node {1};
\blankSpot (1.5,2.75) circle (0.25) node {2};

\blankSpot (0,2.25) circle (0.25) node {3};
\beadSpot (0.75,2.25) circle (0.25) node {4};
\blankSpot (1.5,2.25) circle (0.25) node {5};

\blankSpot (0,1.75) circle (0.25) node {6};
\blankSpot (0.75,1.75) circle (0.25) node {7};
\blankSpot (1.5,1.75) circle (0.25) node {8};

\omitt{
\blankSpot (0,1.25) circle (0.25) node {15};
\blankSpot (0.75,1.25) circle (0.25) node {16};
\blankSpot (1.5,1.25) circle (0.25) node {17};
\blankSpot (2.25,1.25) circle (0.25) node {18};
\blankSpot (3,1.25) circle (0.25) node {19};
}

\end{scope}

\end{scope}

\end{tikzpicture}

\caption{On the left, a $3$-abacus diagram for $(4,2)$. On the right,
  a diagram for $(3,1)$. The right abacus is the result of moving two
  beads in the left abacus, as described in
  Section~\ref{sec:algConn}.}
\label{fig:two3Cores}
\end{figure}

\begin{example}
  \label{ex:richards}
First, we construct the two pyramids for the partition $(4,2)$, one
each for $w=2$ and $w=3$. From Figure~\ref{fig:two3Cores}, we see that
$(p'_0,p'_1,p'_2)=(p_0,p_1,p_2)=(0,1,8)$. When $w=2$, the pyramid is
$(\acoord{0}{2}, \acoord{0}{1}, \acoord{1}{2})=(0,1,0)$ and when $w=3$, it is 
$(\acoord{0}{2}, \acoord{0}{1}, \acoord{1}{2})=(0,2,0)$. We also calculate the
coordinates/admissible sign type of the region containing $(4,2)$'s
alcove. For $m=1$ ($w=2$), $(\aprimecoord{0}{2}, \aprimecoord{0}{1}, \aprimecoord{1}{2})=(1,0,1)$ and for
$m=2$ ($w=3$), $(\aprimecoord{0}{2}, \aprimecoord{0}{1}, \aprimecoord{1}{2})=(2,0,2)$. Additionally,
$$s_0s_2s_1s_0(\emptyset)=s_0s_2s_1(\begin{ytableau}0\end{ytableau})=s_0s_2(\begin{ytableau}0&1\end{ytableau})=s_0(\begin{ytableau}0&1&2\\2\end{ytableau})=\begin{ytableau}0&1&2&0\\2&0\end{ytableau},$$ where the entries in the boxes are their residues mod $3$. We have placed $(4,2)$ in the alcove corresponding to
      $s_0s_2s_1s_0$. See Figure~\ref{fig:alcoveWCores}.

Consider the $1$-Shi region containing $(4,2)$'s alcove.  The minimal
alcove of this region is labeled with the $3$-core $(2)$. To reach
this region from $\Afund$, we must cross only one translate each of
$\Hak{\theta}{1}$ and $\Hak{\alpha_2}{0}$, and no translates of
$\Hak{\alpha_1}{0}.$ This is reflected in $\aprimecoord{0}{2}=1$, $\aprimecoord{1}{2}=1$, and
$\aprimecoord{0}{1}=0$ respectively. Next, we consider the $2$-Shi region
containing $(4,2)$'s alcove. This region, whose minimal element alcove
is labeled by $(4,2)$ itself, is separated from the fundamental alcove
by two translates for each of $\Hak{\theta}{1}$ and $\Hak{\alpha_2}{0}$,
and no translates of $\Hak{\alpha_1}{0}$, reflected in the $m=2$ pyramid
for $(4,2)$.

We repeat the calculations for the $3$-core $(3,1)$, whose $3$-abacus
is on the right in Figure~\ref{fig:two3Cores}.
We calculate
$$s_2s_1s_0(\emptyset)=s_2s_1(\begin{ytableau}0\end{ytableau})=s_2(\begin{ytableau}0&1\end{ytableau})=\begin{ytableau}0&1&2\\2\end{ytableau},$$
    so we have placed $(3,1)$ in the alcove corresponding to
    $s_2s_1s_0$.  The $w=2$ pyramid for $(3,1)$ is $(\acoord{0}{2},\acoord{0}{1},\acoord{1}{2})=(0,1,0)$, the same as the pyramid for $(4,2)$, indicating
    that its alcove will be in the same $1$-Shi region as the alcove
    labeled by $(4,2)$. Its $w=3$ pyramid is $(\acoord{0}{2},\acoord{0}{1},\acoord{1}{2})=(0,2,1)$, which is not the same as the $w=3$ pyramid for
    $(4,2)$, and indeed, their alcoves are in different $2$-Shi
    regions.

Lastly, we mention that the action of moving beads as described in
Section~\ref{sec:algConn} corresponds to flipping the alcove
containing the core of the original abacus over a hyperplane, to the
alcove containing the core obtained through the bead move.

  \end{example}


\subsection{Finite automata and reduced expressions} Headley used
the Shi arrangement to build an automaton which recognizes
reduced expressions. A {\em language} $\Lang$ is a subset of the set $B^*$
of words in a given finite alphabet $B$. For us, the alphabet will be
the set $S$ of a Coxeter group $W$ and the language will be reduced
expressions.

A {\em finite state automaton} is a finite directed graph, with one
vertex designated as the {\em initial state} $S_0$ and a subset of
vertices as final states, and with every edge labeled by an element of
$B$. We call the vertices states. A word $w\in\Lang$ is accepted by the
automaton if the sequence of edge labels along some directed path
starting at $S_0$ and ending at a final state is equal to $w$.

Headley was not the first nor the last to construct an automaton to
accept reduced words; see Bj\"orner and Brenti
\cite{bjorner-brenti2005}, Hohlweg, Nadeau, and Williams
\cite{hohlweg-nadeau-williams2016}, and Gunnells \cite{gunnells2010},
for instance. However, Headley realized that if $s_{i_2}s_{i_3}\cdots
s_{i_k}$ is reduced, then the fundamental alcove $\Afund$ and $\Afund
s_{i_2}s_{i_3}\cdots s_{i_k}$ lie on the same side of the hyperplane
fixed by $s_{i_1}$ if and only if $s_{i_1}s_{i_2}s_{i_3}\cdots
s_{i_k}$ is reduced. \omitt{If however $\Afund$ and
  $s_{i_2}s_{i_3}\cdots s_{i_k}\Afund$ are on opposite sides, then
  $s_{i_1}s_{i_2}s_{i_3}\cdots s_{i_k}$ is not reduced.} He built his
automaton on this observation. The key lemma is

\begin{lemma}\cite{headley1994} Let $W$ be an irreducible affine Weyl group
  with root system $\roots$. Let $R$ be a region of the Shi
  arrangement. If $R$ and $\Afund$ lie on the same side of the
  hyperplanes fixed by $s\in S$, then $Rs$ lies in a single
  region.
\end{lemma}

The states of his automaton are the regions of the Shi
arrangement. The fundamental alcove $\Afund$, which is also a Shi
region, is the initial state. All states are final. Let $s\in S$, and
let $R$ be a region. If $R$ and $\Afund$ are on the same
side of the hyperplane fixed by $s\in S$, then let $R'$ be the region
containing $Rs$ and place an arrow from $R$ to $R'$, labeled by
$s$. If $R$ and $\Afund$ are not on the same side of the hyperplane,
then $R$ has no outgoing arrow labeled by $s$. Headley not only showed
that the language accepted by this automaton is the set of all reduced
words, he showed that if $W$ is the affine symmetric group, the
automaton has the minimal number of states. There were actually two
automata in Headley's thesis. The first produced a nice generating
function, but it has more states.

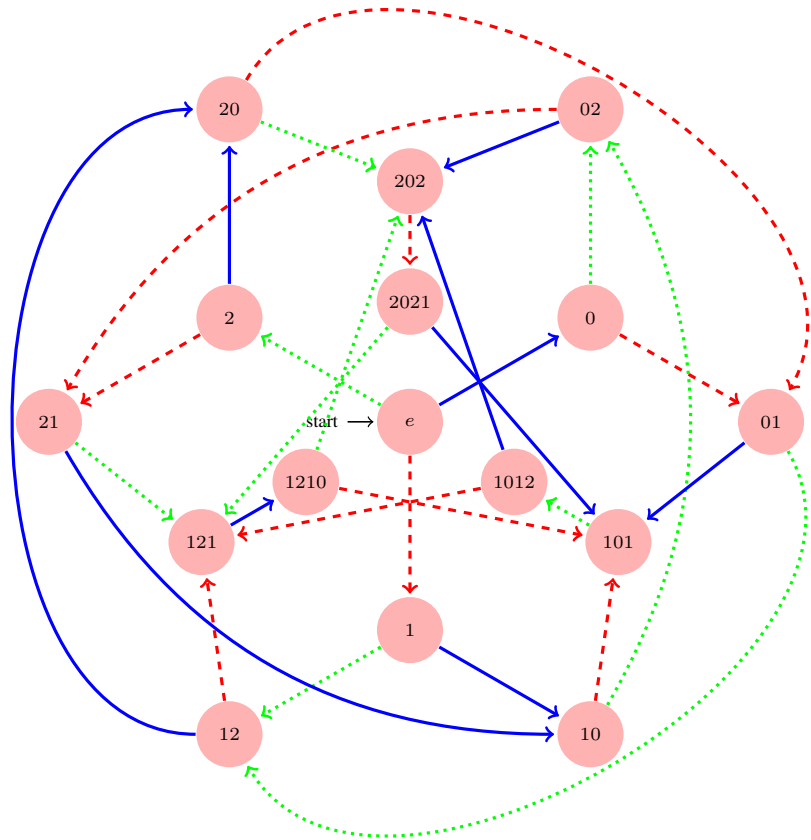
\begin{figure}[h]

\begin{center}
\begin{tikzpicture}[shorten >=1pt,auto,semithick,
every state/.style={fill=red!30,draw=none,text=black},scale=.8,font=\scriptsize]
  \def\a{2}
  \def\r{2}
  \def\d{60}
  \def\c{.8660254038}
  \node[initial,state] (e) at (0:0) {$e$};
\node[state] (0) at (30:2*\r*\c) {$0$};
\node[state] (2) at (150:2*\r*\c) {$2$};
\node[state] (1) at (270:2*\r*\c) {$1$};

\node[state] (101) at (330:2*\r) {$101$};
\node[state] (202) at (90:2*\r) {$202$};
\node[state] (121) at (210:2*\r) {$121$};

\node[state] (01) at (0*\d:3*\r) {$01$};
\node[state] (02) at (1*\d:3*\r) {$02$};
\node[state] (20) at (2*\d:3*\r) {$20$};
\node[state] (21) at (3*\d:3*\r) {$21$};
\node[state] (12) at (4*\d:3*\r) {$12$};
\node[state] (10) at (5*\d:3*\r) {$10$};

\node[state] (1012) at (330:1*\r) {$1012$};
\node[state] (2021) at (90:1*\r) {$2021$};
\node[state] (1210) at (210:1*\r) {$1210$};

\path[->,very thick,blue]
(e) edge  node {} (0)
(e) edge[dashed,red]  node {} (1)
(e) edge[dotted,green] node {} (2)
(1) edge node {} (10)
(1) edge[dotted,green] node {} (12)
(0) edge[dashed,red] node {} (01)
(0) edge[dotted,green] node {} (02)
(2) edge node {} (20)
(2) edge[dashed,red] node {} (21)
(21) edge[dotted,green] node {} (121)
(21) edge[bend right=30] node {} (10)
(12) edge[dashed,red] node {} (121)
(12) edge[bend left=90] node {} (20)
(10) edge[dashed,red] node {} (101)
(10) edge[dotted,green,bend right=30] node {} (02)
(01) edge node {} (101)
(01) edge[dotted,green,bend left=90] node {} (12)
(02) edge node {} (202)
(02) edge[bend right=30,dashed,red] node {} (21)
(20) edge[dotted,green] node {} (202)
(20) edge[dashed,red,bend left=90] node {} (01)
(121) edge node {} (1210)
(101) edge[dotted,green] node {} (1012)
(202) edge[dashed,red] node {} (2021)
(1210) edge[dotted,green] node {} (202)
(1210) edge[dashed,red] node {} (101)
(2021) edge[dotted,green] node {} (121)
(2021) edge node {} (101)
(1012) edge[dashed,red] node {} (121)
(1012) edge node {} (202);
\end{tikzpicture}

\end{center}

  \caption{Headley's automaton based on the Shi arrangement. It
    accepts reduced words in $\affS{3}$. Each state is labeled with
    the affine permutation corresponding to the minimal alcove of the
    region. We use $i$ instead of $s_i$. Solid arrows represent an
    edge labeled by $s_0$ and dashed (respectively dotted) edges
    represent edges labeled $s_1$ (respectively $s_2$).  See
    Example~\ref{ex:auto}.}
  \label{fig:auto}
  \end{figure}

\begin{example}
  \label{ex:auto}
  This example refers to Figure~\ref{fig:auto}. The path
  $$e\xrightarrow{0} 0\xrightarrow{2} 02\xrightarrow{1} 21\xrightarrow{0} 10\xrightarrow{1} 101\xrightarrow{2} 1012$$ represents the expression $s_0s_2s_1s_0s_1s_2$, which is reduced. Since $s_0s_1s_0s_1$ is not reduced, there is no path which starts at $e$ and follows solid, dashed, solid, dashed arrows. 
\end{example}

\subsection{More connections}
This connection is to the filters in $\pos$, not to the Shi
arrangement directly. In \cite{cellini-papi2000,cellini-papi2002,
  cellini-papi2004}, Cellini and Papi investigate $\ad$-nilpotent
ideals in a Borel subalgebra. They associate each ideal to a filter in
$\pos$ and also to an element of the affine Weyl group. \omitt{The
  elements $w$ which are associated to some ideal turn out to be the
  $w$ such that $w\Afund$ is a minimal alcove for a Shi region!
  \sftodo{Check this!!}} Very roughly speaking, they use the Cartan
decomposition $L=H\oplus N$ and $N= \bigoplus_{\alpha\in\pos}
L_{\alpha}$ and the definition of an $\ad$-nilpotent ideal as an ideal
contained in $N$ to define the
antichain $$\roots_I=\{\alpha\in\pos:L_{\alpha}\subseteq I\}$$ which
defines a filter. See also Suter \cite{suter2004}. Dong extends
Cellini and Papi's work from Borel subalgebras to parabolic
subalgebras in \cite{dong2013}. He uses deleted Shi arrangements,
which we don't address in this survey. Panyushev \cite{panyushev2004}
developes combinatorial aspects of the theory of $\ad$-nilpotent
ideals, giving a geometric interpretationfor the number of generators
of an ideal, for example.

Gunnells and Sommers study {\em Dynkin elements}, which we won't
define, in \cite{gunnells-sommers2003}. They define $N$-regions, which
turn out to be unions of Shi regions. A simplified version of their
main theorem is that if $x$ is the point of minimal Euclidean length
in the closure of an $N$-region, then $2x$ is a Dynkin element.


\section{Further}\label{sec:further}

We briefly mention a few recent results. In
\cite{gorsky-mazin-vazirani2016}, Gorsky, Mazin, and Vazirani developed
``rational slope'' versions of much of what has been discussed here. A
tuple $(b_1,\ldots,b_n)$ of nonnegative integers is called an
$M/n$-parking function if the Young diagram with row lengths equal to
$b_1,\ldots,b_n$ arranged in decreasing order fits above the diagonal
in an $n\times M$ rectangle. We've stated it a bit differently than in
Sections~\ref{sec:stanEnum} and \ref{sec:extShi}, but if we let
$M=n+1$ and $M=mn+1$ respectively and reverse the order of the tuple, we obtain the same
functions.
Gorsky, Mazin, and Vazirani defined {\em $M$-stable} permutations to
take the place of minimal permutations of Shi regions and generalized
the Pak-Stanley bijection, as well as another map defined by Anderson
\cite{anderson2002} in her study of core partitions. They conjectured
their generalization of the Pak-Stanley map is injective for all
relatively prime $M$ and $n$. In 2017, McCammond, Thomas, and Williams \cite{MTW2017} proved the conjectures in \cite{gorsky-mazin-vazirani2016}. Additionally, Gorsky, Mazin, and Vazirani connect their maps to
the combinatorics of $q,t$-Catalan polynomials. Sulzgruber
\cite{sulzgruber2015} built on \cite{gorsky-mazin-vazirani2016} by
finding the coordinates of the $M$-stable permutations, generalizing
\eqref{eq:mcoords}. Thiel \cite{thiel2014} extended their work to other
types, among other results.

As mentioned in Section~\ref{sec:enumeration} Hohlweg, Nadeau, and
Williams generalized the Shi arrangement to any Coxeter group, using
$n$-small roots, and then to indefinite Coxeter systems. They also
investigated automata.

\section{Themes we haven't included}
\label{sec:themes}
We give a short and incomplete list of topics we have not discussed.

\begin{enumerate}
\item The Shi arrangement is free.
  Either see original article by Athanasiadis \cite{athanasiadis1998b} or his excellent summary \cite{athanasiadis1998}. 
Abe, Suyama, and Tsujie \cite{abe-suyama-tsujie2017} show that the Ish arrangement is free.\omitt{See exercise in \cite{stanley2007} on something not being supersolvable.
Is original \cite{athanasiadis1998b}? See \cite{armstrong-rhoades2012}, p 11, Is \cite{athanasiadis1998} the survey?}
\item In graphical arrangements or deleted arrangements, some of the hyperplanes have been removed. We survey only the complete Shi arrangement.
\item We have no discussion of the connections to the torus $\check{Q}/(1+mh)\check{Q}$, where $\check{Q}$ is the coroot lattice of a root system, $(mh+1)\check{Q}$ is its dilate, and $h$ is the Coxeter number of the root system. See Athanasiadis \cite{athanasiadis2005} or Haiman \cite{haiman1994} for more information. 
\item The enumeration of bounded regions has nice results, which we
  have not discussed. See Athanasiadis and Tzanaki
  \cite{athanasiadis-tzanaki2006}, for example and Sommers \cite{sommers2005}.
\end{enumerate}

\section{Acknowlegements}
The author would like to thank H\'el\`ene Barcelo, Gizem Karaali, and
Rosa Orellana for allowing her to write an article for this AWM
series. She would like to thank Matthew Fayers, Sarah Mason, and
Jian-Yi Shi for comments on the manuscript, and Christos Athanasiadis,
Duncan Levear, Ant\'{o}nio Guedes de Oliveira, Brendon Rhoades, and
Nathan Williams for their help with references. She would like to
thank Patrick Headley for help with his thesis. The anonymous
referees' comments helped enormously to improve exposition.  She
thanks George Lusztig for correcting an error. This work was supported
by a grant from the Simons Foundation (\#359602, Susanna Fishel).

\omitt{
\section{Don't forget!!}
\begin{enumerate}
\item From \cite{athanasiadis1998} \cite{ram2003} invariant theory \cite{solomon-terao1998}\cite{ram2003} p 387 p 373
  
  \item Athanasiadis \cite{athanasiadis2010} recip for coeff of char poly? maybe with char poly subsec
\item More on what is known about regions inside of cells of other types.\cite{shi1988,shi1987c}.
\item In \cite{athanasiadis2004} indec. play role of min elts. p186, rep alcove establishe by shi 7.2 also p 186

 \item \cite{hughes-williams2014} labeling
 \item \cite{sommers2005} inverses form simplex, ideals
   \item \cite{panyushev2004} ad nilpotent, inverses form simples
\end{enumerate}
}

\bibliographystyle{amsalpha}
\label{sec:biblio}
\newcommand{\etalchar}[1]{$^{#1}$}
\providecommand{\bysame}{\leavevmode\hbox to3em{\hrulefill}\thinspace}
\providecommand{\MR}{\relax\ifhmode\unskip\space\fi MR }
\providecommand{\MRhref}[2]{%
  \href{http://www.ams.org/mathscinet-getitem?mr=#1}{#2}
}
\providecommand{\href}[2]{#2}

\end{document}